\documentclass[12pt]{amsart}

\usepackage[dvips]{graphics}

\usepackage{epsfig}
\usepackage{stmaryrd}

\usepackage{amsfonts}
\usepackage{amsmath, amsthm, amssymb, ulem, amscd} 

\newcommand{\stopthm}{\hfill$\square$\medskip}

\def\crn#1#2{{\vcenter{\vbox{
        \hbox{\kern#2pt \vrule width.#2pt height#1pt
           }
          \hrule height.#2pt}}}}
\def\intprod{\mathchoice\crn54\crn54\crn{3.75}3\crn{2.5}2}
\def\into{\mathbin{\intprod}}

\newcommand{\cP}{{\mathcal P}} 
\newcommand{\pa}{\partial}
\newcommand{\Up}{\Upsilon}

\newcommand{\ep}{\epsilon}
\newcommand{\al}{\alpha}
\newcommand{\be}{\beta}
\newcommand{\la}{\lambda}
\newcommand{\ka}{\kappa} 
\newcommand{\ga}{\gamma}  
\newcommand{\Ga}{{\Gamma}}
\newcommand{\Gat}{{\tilde \Gamma}}

\newcommand{\st}{{\tilde{\sigma}}}
\newcommand{\et}{{\tilde{\varepsilon}}}
\newcommand{\ttt}{{\tilde{\tau}}}
\newcommand{\ett}{{\tilde{\eta}}}
\newcommand{\mut}{{\tilde{\mu}}}
\newcommand{\cG}{{\mathcal G}}
\newcommand{\cM}{{\mathcal M}}
\newcommand{\cGt}{{\tilde {\mathcal G}}}
\newcommand{\cRt}{{\tilde {\mathcal R}}}
\newcommand{\cTt}{{\tilde {\mathcal T}}}
\newcommand{\Pt}{{\tilde {\Phi}}}
\newcommand{\Dt}{{\tilde {\Delta}}}  
\newcommand{\Lt}{{\tilde {\Lambda}}} 
\newcommand{\cO}{{\mathcal O}}
\newcommand{\cR}{{\mathcal R}}
\newcommand{\cV}{{\mathcal V}}
\newcommand{\cU}{{\mathcal U}}
  
\newcommand{\cT}{{\mathcal T}}
\newcommand{\cS}{{\mathcal S}}
\newcommand{\cN}{{\mathcal N}}
\newcommand{\cA}{{\mathcal A}}
\newcommand{\cE}{{\mathcal E}}
\newcommand{\cL}{{\mathcal L}}

\newcommand{\gt}{{\tilde g}}
\newcommand{\gb}{{\overline g}} 
\newcommand{\mb}{{\overline \mu}}  
\newcommand{\Wb}{{\overline W}} 
\newcommand{\htt}{{\tilde h}}

\newcommand{\Rt}{{\tilde R}}
\newcommand{\cH}{{\mathcal H}}
\newcommand{\Z}{{\mathbb Z}}
\newcommand{\R}{{\mathbb R}}

\newcommand{\N}{{\mathbb N}}

\newcommand{\C}{{\mathbb C}}

\newcommand{\nf}{\infty}
\newcommand{\nt}{\tilde \nabla}

\newcommand{\Ric}{\operatorname{Ric}}
\newcommand{\contr}{\operatorname{contr}}
\newcommand{\pcontr}{\operatorname{pcontr}}
\newcommand{\tf}{\operatorname{tf}}

\renewcommand{\tilde}{\widetilde} 
\newcommand{\wh}{\widehat} 

\theoremstyle{plain}
\newtheorem{theorem}{Theorem}[section]
\newtheorem{lemma}[theorem]{Lemma}
\newtheorem{proposition}[theorem]{Proposition}

\theoremstyle{definition}

\newtheorem{definition}[theorem]{Definition}

\theoremstyle{remark}
\newtheorem{remark}[theorem]{Remark}

\numberwithin{equation}{section}

\title{The Ambient Metric}

\author{Charles Fefferman}
\address{Department of Mathematics, Princeton University\\
  Princeton, NJ 08544}
\email{cf@math.princeton.edu}

\author{C. Robin Graham}
\address{Department of Mathematics, University of Washington,
Box 354350\\
Seattle, WA 98195-4350}
\email{robin@math.washington.edu}

\begin{document}

\maketitle

\tableofcontents

\thispagestyle{empty}

\section{Introduction}\label{intro}

Conformal geometry is the study of spaces in which one knows how to measure
infinitesimal angles but not lengths.  A conformal structure on a manifold 
is an equivalence class of Riemannian metrics, in which 
two metrics are identified if one is a positive smooth multiple of the
other.  The study of conformal geometry has a long and venerable history.   
From the beginning, conformal geometry has played an important role in
physical theories.   

A striking historical difference between conformal geometry compared
with Riemannian geometry is the scarcity of local invariants in the
conformal case.  Classically known conformally invariant tensors include
the Weyl conformal curvature tensor, which plays 
the role of the Riemann curvature tensor, its three-dimensional
analogue the Cotton tensor, and the Bach tensor in dimension four.
Further examples are not so easy to
come by.  By comparison, in the Riemannian case  
invariant tensors abound.  They can be easily constructed by covariant 
differentiation of the curvature tensor and tensorial operations.  The
situation is similar for 
other types of invariant objects, for example for differential operators.  
Historically, there are 
scattered examples of conformally invariant operators such as the
conformally invariant Laplacian and certain Dirac operators, whereas  it is 
easy to write down Riemannian invariant differential operators, of
arbitrary orders and between a wide variety of bundles.   

In Riemannian geometry, not only is it easy to write down invariants, it 
can be shown using Weyl's classical invariant theory for the orthogonal
group 
that all invariants arise by the covariant differentiation and tensorial
operations mentioned above.  In the case of scalar invariants, this
characterization as ``Weyl invariants'' has had important
application in the study of heat asymptotics: one can immediately write
down the form of coefficients in the expansion of heat kernels up to 
the determination of numerical coefficients.  

In \cite{FG}, we outlined a construction of a nondegenerate Lorentz metric
in $n+2$ dimensions associated to an $n$-dimensional conformal manifold,
which we called the ambient metric.  This association enables one to
construct conformal invariants in $n$ dimensions from pseudo-Riemannian
invariants in $n+2$ dimensions, and in particular shows that conformal
invariants are plentiful.  The construction of conformal invariants is
easiest and 
most effective for scalar invariants:  every scalar invariant of metrics in
$n+2$ dimensions immediately determines a scalar conformal invariant in $n$
dimensions (which may vanish, however).  For other types of invariants, for
example for differential operators, some effort is required to derive a  
conformal invariant from a pseudo-Riemannian invariant in two higher
dimensions, but in many cases this can be carried out and has led 
to important new examples.  

The ambient metric is homogeneous with respect to a family of dilations on
the $n+2$-dimensional space.  It is possible to mod out 
by these dilations and thereby obtain a metric in $n+1$ dimensions, 
also associated to the given conformal manifold in $n$ dimensions.  This
gives the ``Poincar\'e'' metric associated to the conformal manifold.  The   
Poincar\'e metric is complete and the conformal manifold forms its boundary
at infinity.

The construction of the ambient and Poincar\'e metrics associated to a
general conformal manifold is motivated by the
conformal geometry of the flat model, the sphere $S^n$, 
which is naturally described in terms of $n+2$-dimensional Minkowski space.
Let 
$$
Q(x)=\sum_{\al=1}^{n+1} \left(x^\al\right)^2 - \left(x^0\right)^2 
$$
be the standard Lorentz signature quadratic form on $\R^{n+2}$ and 
$$
\mathcal{N}=\left\{x\in \R^{n+2}\setminus \{0\}: Q(x)=0\right\}
$$ 
its null cone.  
The sphere $S^n$ can be identified with the space
of lines in $\mathcal{N}$, with projection $\pi:\mathcal{N}\rightarrow
S^n$.   Let  
$$
\gt = \sum_{\al=1}^{n+1} \left(dx^\al\right)^2 - \left(dx^0\right)^2 
$$
be the associated Minkowski metric on $\R^{n+2}$.  The conformal structure
on $S^n$
arises by restriction of $\gt$ to $\mathcal{N}$.  More specifically,    
for $x\in\mathcal{N}$ the restriction  
$\gt|_{T_x\mathcal{N}}$ is a degenerate quadratic form because it
annihilates the radial 
vector field $X=\sum_{I=0}^{n+1} x^I\pa_I\in T_x\mathcal{N}$.  So
$\gt|_{T_x\mathcal{N}}$       
determines an inner product on $T_x\mathcal{N}/\operatorname{span}{X}\cong 
T_{\pi(x)}S^n$.  As $x$ varies over a line in $\mathcal{N}$, the  
resulting inner products on $T_{\pi(x)}S^n$ vary only by scale and  
are the possible
values at $\pi(x)$ for a metric in the conformal class on $S^n$.  
The Lorentz group $O(n+1,1)$ acts linearly on $\R^{n+2}$ by isometries 
of 
$\gt$ preserving $\mathcal{N}$.  The induced action on lines in
$\mathcal{N}$ therefore preserves the conformal class of metrics and
realizes the group of conformal motions of $S^n$.        

If instead of restricting $\gt$ to $\mathcal{N}$, we 
restrict it to the hyperboloid 
$$
\mathcal{H}=\{x\in \R^{n+2}:\,Q(x)=-1 \},
$$
then we obtain the Poincar\'e metric associated to $S^n$.  Namely, 
$g_+:=\gt|_{T\mathcal{H}}$ is the hyperbolic metric of constant sectional 
curvature $-1$.  Under an appropriate identification of one sheet of
$\mathcal{H}$ with 
the unit ball in $\R^{n+1}$, $g_+$ can be realized as the Poincar\'e metric  
$$
g_+=4\left(1-|x|^2\right)^{-2} 
\sum_{\al=1}^{n+1}\left(dx^\al\right)^2, 
$$
and has the conformal structure on $S^n$ as conformal infinity.  The action
of
\linebreak  
$O(n+1,1)$ on $\R^{n+2}$ preserves $\mathcal{H}$.  The induced action  
on $\mathcal{H}$ is
clearly by isometries of $g_+$ and realizes the isometry group of
hyperbolic space.   

The ambient and Poincar\'e metrics associated to a general conformal
manifold are defined as solutions to certain systems of partial
differential equations with initial data determined by the conformal
structure.  Consider the ambient metric.  
A conformal class of metrics on a manifold $M$ determines and is determined
by its metric bundle $\cG$, an $\R_+$-bundle over $M$.  This is the
subbundle of symmetric 2-tensors whose sections are the  
metrics in the conformal class.  In the case $M=S^n$, 
$\cG$ can be identified with $\mathcal{N}$ (modulo $\pm 1$).  Regard 
$\cG$ as a hypersurface in $\cG\times \R$:  $\cG\cong\cG\times \{0\}\subset
\cG\times \R$.     
The conditions defining the ambient metric 
$\gt$ are that it be a Lorentz metric defined in a neighborhood of 
$\cG$ in $\cG\times \R$ which is homogeneous with respect to the natural
dilations on this space, that it satisfy an initial condition 
on the initial hypersurface $\cG$ determined by the conformal structure, 
and that it be 
Ricci-flat.  The Ricci-flat condition is the system of equations intended
to propagate the initial data off of the initial surface.     

This initial value problem is singular because the pullback of $\gt$ to 
the initial surface is degenerate.  However, for the applications to
the construction of conformal invariants, it is sufficient to have formal
power series solutions along the initial surface rather than actual
solutions in a neighborhood.  So we concern ourselves with the formal 
theory and do not discuss the interesting but more difficult problem of
solving the equations exactly.  

It turns out that the behavior of solutions of this system depends
decisively on the parity of the 
dimension.  When the dimension $n$ of the conformal 
manifold is odd, there exists a formal power series solution $\gt$ which is
Ricci-flat to infinite order, and it is unique up to   
diffeomorphism and up to terms vanishing to infinite order.  When $n\geq 4$
is even, there is a solution which is Ricci-flat to order $n/2-1$, 
uniquely determined up to diffeomorphism and up to terms vanishing to
higher order.  But at this order, the existence of smooth solutions   
is obstructed by a conformally invariant natural
trace-free symmetric 2-tensor, the ambient obstruction tensor.  When $n=4$, 
the obstruction tensor is the same as the classical Bach tensor.
When $n=2$, there is no obstruction, but uniqueness fails.  

It may seem contradictory that for $n$ odd, the solution of a second order
initial value  
problem can be formally determined to infinite order by only one piece of
Cauchy data: the 
initial condition determined by the conformal structure.  In fact, there
are indeed further formal solutions.  These correspond to the freedom of a
second piece of initial data at order $n/2$.  When $n$ is odd, $n/2$ is 
half-integral, and this freedom is  
removed by restricting to formal power series solutions.  It is  
crucially important for the applications to conformal invariants that we
are able to  
uniquely specify an infinite order solution in an invariant way.  On the
other hand, the additional formal solutions with nontrivial asymptotics at
order $n/2$ are also important; they necessarily arise in the global
formulation of the existence problem as a boundary value problem at
infinity for the Poincar\'e metric.  When $n$ is even, the
obstruction to 
the existence of formal power series solutions can be incorporated into 
log terms in the expansion, in which case there is again a 
formally undetermined term at order $n/2$.   

The formal theory described above was outlined in \cite{FG}, but the 
details were not given.  The first main goal of this monograph is to    
provide these details.  We give the full infinite-order formal theory,
including the freedom at order $n/2$ in all dimensions and the precise
description of the log terms when $n\geq 4$ is even.  This formal theory
for the ambient metric forms the content of Chapters~\ref{setup} and  
\ref{formalsect}.  The description of the solutions with freedom at order
$n/2$ and log terms extends and sharpens results of Kichenassamy \cite{K}.        
Convergence of the formal series determined by singular nonlinear initial
value problems of this type has been considered by several authors; these 
results imply that the formal series converge if the data are   
real-analytic.   

In Chapter~\ref{poincaresection}, we define Poincar\'e metrics:  they are
formal solutions to the equation $\Ric(g_+)=-ng_+$, and we show how   
Poincar\'e metrics are equivalent to ambient metrics satisfying an extra
condition which we call straight.  Then we use this equivalence to derive
the full formal theory for Poincar\'e metrics from that for ambient
metrics.  We discuss the
``projectively compact'' formulation of Poincar\'e metrics, modeled on the
Klein model of hyperbolic space, as well as the usual conformally compact
picture.  
As an application of the formal theory for Poincar\'e metrics, in 
Chapter~\ref{selfdualsect} we present a  
formal power series proof of a result of LeBrun \cite{LeB} asserting 
the existence and uniqueness of a real-analytic self-dual Einstein metric
in dimension 4 defined near the boundary 
with prescribed real-analytic conformal infinity.  

In Chapter~\ref{flat}, we analyze the ambient and Poincar\'e metrics for   
locally conformally flat manifolds and for conformal classes containing an 
Einstein metric.   The obstruction tensor vanishes for even
dimensional conformal structures of these types.  We show that for these
special conformal classes, 
there is a way to uniquely specify the 
formally undetermined term at order $n/2$ in an invariant way and thereby
obtain a unique ambient metric up to terms vanishing to infinite order   
and up to diffeomorphism, just like in odd dimensions.  We derive a   
formula of Skenderis and Solodukhin \cite{SS} for the ambient or
Poincar\'e metric in the locally 
conformally flat case which is in normal form relative to an arbitrary
metric in the conformal class, and prove a related unique continuation 
result for 
hyperbolic metrics in terms of data at conformal infinity.  The case $n=2$
is special for all of these considerations.  We also derive the form of  
the GJMS operators for an Einstein metric.

In \cite{FG}, we conjectured that when $n$ is odd, 
all scalar conformal invariants arise as Weyl invariants constructed  
from the ambient metric.      
The second main goal of this monograph is to prove this together with
an analogous result when $n$ is even.  
These results are contained in Theorems~\ref{vanishing},
\ref{oddlowweight}, and \ref{invariants}.  When $n$ is even, we  
restrict to   
invariants whose weight $w$ satisfies $-w\leq n$ because of the finite
order indeterminacy of the ambient metric:  Weyl invariants of higher
negative weight may involve derivatives of the ambient metric which are not 
determined.  
A particularly interesting phenomenon occurs in dimensions $n\equiv 0 \mod 
4$.  For all $n$ even, it is the case that all even   
(i.e. unchanged under orientation reversal) scalar conformal
invariants with $-w\leq n$ arise as Weyl invariants of the ambient metric.    
If $n\equiv 2 \mod 4$, this is also true for odd (i.e. changing sign under
orientation reversal) scalar conformal invariants with $-w\leq n$ (in 
fact, these all vanish).  But if $n\equiv 0 \mod 4$, there are    
odd invariants of weight $-n$ which     
are exceptional in the sense that they do not arise as Weyl invariants of
the ambient metric.  The set of such exceptional invariants of
weight $-n$ consists precisely of the nonzero elements of the vector space 
spanned by 
the Pontrjagin invariants whose integrals give the 
Pontrajgin numbers of a compact oriented $n$-dimensional manifold
(Theorem~\ref{oddlowweight}).   

The parabolic invariant theory needed to prove these results 
was developed in \cite{BEGr}, including the observation of the existence    
of exceptional invariants.  But substantial work is required to
reduce the theorems in Chapter~\ref{sinv} to the results of \cite{BEGr}. 
To understand this, we briefly review how Weyl's characterization of scalar
Riemannian invariants is proved.    

Recall that Weyl's theorem for even invariants states that every even
scalar Riemannian invariant is a linear   
combination of complete contractions of the form 
$
\contr\left(\nabla^{r_1}R\otimes \cdots \otimes\nabla^{r_L}R\right),
$
where the $r_i$ are nonnegative integers, $\nabla^{r}R$ denotes the $r$-th
covariant derivative of the curvature tensor, and $\contr$ 
denotes a metric contraction with respect to some pairing of all the
indices.  There are two main steps in the proof of Weyl's theorem.  The
first is to show that    
any scalar Riemannian invariant can be written as a polynomial in the
components of the covariant derivatives of the curvature tensor which is
invariant under the orthogonal group $O(n)$.  Since a Riemannian invariant
by definition is a polynomial in the Taylor coefficients of the metric in  
local coordinates whose value is independent of the choice of coordinates, 
one must pass from Taylor coefficients of the metric to covariant
derivatives of curvature.  This passage is carried out using geodesic
normal coordinates.  We refer to the result stating that the map from  
Taylor coefficients of the metric in geodesic normal coordinates to 
covariant derivatives of curvature is an $O(n)$-equivariant isomorphism as 
the jet isomorphism theorem for Riemannian geometry.  Once the jet
isomorphism theorem has been established, one is left with the algebraic 
problem of identifying the $O(n)$-invariant polynomials in the covariant
derivatives of curvature. This is solved by Weyl's classical invariant 
theory for the orthogonal group.

In the conformal case, the role of the covariant derivatives of the 
curvature tensor is played by the covariant derivatives of the curvature
tensor of 
the ambient metric.  These tensors are of course defined on the ambient
space.  But when evaluated on the initial surface, their components  
relative to a suitable frame define tensors on the base conformal   
manifold, which we call conformal curvature tensors.  For example, the
conformal curvature tensors defined by the curvature tensor of the ambient
metric itself (i.e. with no ambient covariant derivatives) are the
classical Weyl, Cotton, and Bach tensors (except that in dimension 4, the
Bach tensor does not arise as a conformal curvature tensor because of the  
indeterminacy of the ambient metric).  The covariant derivatives of
curvature of the
ambient metric satisfy identities and relations beyond those 
satisfied for general metrics owing to its homogeneity and Ricci-flatness.  
We derive these identities in Chapter~\ref{cct}.  We also derive the
transformation  
laws for the conformal curvature tensors under conformal change.   
Of all the conformal curvature tensors, 
only the Weyl tensor (and Cotton tensor in dimension 3) are conformally
invariant.  The transformation law of any other conformal curvature tensor  
involves only first derivatives of the conformal factor and ``earlier''
conformal curvature tensors.  These transformation laws may also be
interpreted in terms of tractors.  When $n$ is even, the definitions of the 
conformal curvature tensors and the identities which they satisfy are
restricted by the finite order indeterminacy of the ambient metric.  
The ambient obstruction tensor is not a conformal curvature tensor; it 
lies at the boundary of the range for which they are defined.  But it may
be regarded as the residue of an analytic continuation in the dimension of 
conformal curvature tensors in higher dimensions
(Proposition~\ref{bachgen}). 

Having understood the properties of the conformal curvature tensors, the
next step in the reduction of the theorems in Chapter~\ref{sinv} to the 
results of \cite{BEGr} is to formulate and prove a jet isomorphism theorem
for conformal geometry, in order to know that a scalar conformal invariant 
can be written in terms of conformal curvature tensors.  The  
Taylor expansion of the metric on the base manifold in geodesic 
normal coordinates can be further simplified since one now has the freedom
to change the metric by a conformal factor as well as by a diffeomorphism.  
This leads to a ``conformal normal form'' in which part of the base
curvature is 
normalized away to all orders.  Then the conformal jet isomorphism theorem
states that the map from the Taylor coefficients of a  
metric in conformal normal form to the space of all conformal curvature
tensors, realized as covariant derivatives of ambient curvature, is an 
isomorphism.  
Again, the spaces must be truncated at finite order in even dimensions.   
The proof of the conformal jet isomorphism theorem is much more involved
than in the Riemannian case; it is  
necessary to relate the normalization conditions in the conformal normal
form to the precise identities and relations satisfied by the ambient
covariant derivatives of curvature.  We carry this out in    
Chapter~\ref{jet} by making a direct algebraic study of these relations and
of the map from jets of normalized metrics to conformal curvature tensors.   
A more conceptual proof of the conformal jet isomorphism theorem due to the 
second author and K. Hirachi uses an ambient lift of the conformal
deformation complex and is outlined in \cite{Gr3}.      

The orthogonal group plays a central role in Riemannian geometry
because it is the isotropy group of a point in the group of isometries 
of the flat model $\R^n$.  The analogous group for  
conformal geometry is the isotropy subgroup $P\subset O(n+1,1)$ of the 
conformal group fixing a point in $S^n$, i.e. a null line.  Because of its
algebraic structure, 
$P$ is referred to as a parabolic subgroup of $O(n+1,1)$.      
Just as geodesic normal coordinates are determined up to the
action of $O(n)$ in the Riemannian case, the equivalent conformal 
normal forms for a metric at a given point are determined up to an  
action of $P$.  Since $P$ is a matrix group in
$n+2$ dimensions, there is  
a natural tensorial action of $P$ on the space of covariant derivatives of
ambient 
curvature, and the conformal transformation law for conformal curvature
tensors established in Chapter~\ref{cct} implies that the map from jets of
metrics in conformal normal form to 
conformal curvature tensors is $P$-equivariant.  

The jet isomorphism theorem reduces the study of conformal invariants to
the purely algebraic matter of understanding the $P$-invariants of the
space of covariant derivatives of ambient curvature.  This space is
nonlinear since the Ricci identity for commuting covariant derivatives is
nonlinear in curvature and its derivatives.  The results of 
\cite{BEGr} identify the $P$-invariants of the linearization of this space.
So the last steps, carried out in Chapter~\ref{sinv}, are to formulate the   
results about scalar invariants, to use the jet isomorphism theorem to
reduce these results to algebraic statements in invariant theory for $P$,
and finally to  
reduce the invariant theory for the actual nonlinear space to that for its
linearization.   
The treatment in Chapters~\ref{jet} and \ref{sinv} is inspired  
by, and to some   
degree follows, the treatment in \cite{F} in the case of CR geometry.   

Our work raises the obvious question of extending the theory to higher
orders in even dimensions.  This has recently been carried out by the 
second author and K. Hirachi.  An extension to all orders of the ambient
metric construction, jet isomorphism theorem, and invariant theory  has
been announced in \cite{GrH2}, \cite{Gr3} inspired by the work of Hirachi 
\cite{Hi} in the CR case.  The log terms in the expansion of an ambient
metric are modified by taking the log of a defining function homogeneous of
degree 2 rather than homogeneous of degree 0.  
This makes it possible to define the smooth part of an 
ambient metric with log terms in an invariant way.  The smooth part is
smooth and homogeneous but no 
longer Ricci flat to infinite order.  There is a family of such smooth
parts corresponding to 
different choices of the ambiguity at order $n/2$.  They can be used to
formulate a jet isomorphism theorem and to construct invariants, and the
main conclusion is that up to a linear combination of finitely many   
exceptional odd invariants 
in dimensions $n\equiv 0\mod 4$ which can be explicitly identified,    
all scalar conformal invariants arise from the ambient metric.  
An alternate development of a conformal invariant theory based on tractor 
calculus is given in general dimensions in \cite{Go1}.        

A sizeable literature concerning the ambient and Poincar\'e metrics has 
arisen since the publication of \cite{FG}.  
The subject has been greatly stimulated by its relevance in the study of
the AdS/CFT correspondence in physics.
We have tried to indicate some
of the most relevant references of which we are aware without attempting to
be exhaustive.   

A construction equivalent to the ambient metric was derived by
Haantjes and Schouten in 
\cite{HS}.  They obtained a version of the expansion for  
straight ambient metrics, to infinite order in odd dimensions and up to the 
obstruction in even dimensions.  In particular, 
they showed that there is an obstruction in even dimensions $n\geq 4$ and
calculated that it is the 
Bach tensor in dimension 4.  They observed that the obstruction vanishes
for conformally Einstein metrics and in this case derived the conformally    
invariant normalization uniquely specifying an infinite order ambient
metric in even dimensions.  They also obtained the   
infinite-order  expansion in the case of dimension 2, including the precise
description of the non-uniqueness of solutions.  Haantjes and Schouten did
not consider applications to conformal invariants and, unfortunately, it
seems that their work was largely forgotten.      

We are grateful to the National Science Foundation for support.  In
particular, the second author was partially supported by NSF grant 
\# DMS 0505701 during the preparation of this manuscript.   

Throughout, by smooth we will mean infinitely 
differentiable.  Manifolds are assumed to be smooth and second countable;
hence paracompact.  
Our setting is primarily algebraic, so we work with metrics of  
general signature.  In tensorial expressions, we denote by parentheses
$(ijk)$ symmetrization and by brackets $[ijk]$  skew-symmetrization  
over the enclosed indices.  

\section{Ambient Metrics}\label{setup} 

\setcounter{section}{2}

\medskip

Let $M$ be a smooth manifold of dimension $n \geq 2$ equipped  
with a conformal class $[g]$.  Here, $g$ is a smooth pseudo-Riemannian
metric of signature $(p , q)$ on $M$ and $[g]$ consists of all
metrics 
$$
\wh{g} = e^{2 \Up} g
$$ 
on $M$, where $\Up$ is any smooth real-valued function on $M$.

The space $\cG$ consists of all pairs $( h , x )$, where $x \in M$, and
$h$ is a symmetric bilinear form on $T_x M$ satisfying $h = s^2
g_x$ for some $s \in \R_+$.  Here and below, $g_x$ denotes the 
symmetric bilinear form on $T_xM$ induced by the metric $g$.
We write $\pi : \cG \rightarrow M$ for the projection map $(h , x )
\mapsto x$. Also, for $s \in \R_+$, we define the ``dilation''
$\delta_s : \cG \rightarrow \cG$ by setting $\delta_s ( h , x )
= ( s^2 h , x )$.  The space $\cG$, equipped with the projection
$\pi$ and the dilations $(\delta_s)_{s \in \R_+}$, 
is an $\R_+$-bundle.  We call it the {\it metric bundle} for $(M, [ g ])$.   
We denote by $T= \frac{d}{ds} \delta_s |_{s=1}$ the vector field on $\cG$ which is
the infinitesimal generator of the dilations $\delta_s$.

There is a tautological symmetric 2-tensor ${\bf g}_0$ on $\cG$, defined as
follows.  Let $z = ( h , x ) \in \cG$, and let 
$\pi_*: T\cG \rightarrow TM$ be the differential of the map
$\pi$.
Then, for tangent vectors $X, Y \in T_z \cG$, we define ${\bf g}_0 ( X , Y ) =
h( \pi_* X , \pi_* Y )$.  
The 2-tensor ${\bf g}_0$ is homogeneous of degree 2 with respect to the 
dilations $\delta_s$.  That is, $\delta_s^* {\bf g}_0 = s^2 {\bf g}_0$.
One checks easily that the $\R_+$-bundle $\cG$, the maps
$\delta_s$ and $\pi$, and the tautological 2-tensor ${\bf g}_0$ on
$\cG$, all depend only on the conformal class $[g]$, and are independent
of the choice of the representative $g$.
However, once we fix a representative $g$, we obtain a trivialization of
the bundle $\cG$.  In fact, we identify
$$
( t , x ) \in \R_+ \times M \ \mbox{ with } \ (t^2 g_x , x ) \in \cG \,.
$$
In terms of this identification, the dilations $\delta_s$, the
projection $\pi$, the vector field $T$, and the tautological 2-tensor ${\bf g}_0$
are given by 
$$
\delta_s\,:\,( t , x ) \mapsto (st , x ),\qquad
\pi \, :\,(t , x ) \mapsto x ,\qquad 
T   = t \partial_t, \qquad {\bf g}_0 = t^2\pi^*g.
$$
The metric $g$ can be regarded as a section of the bundle $\cG$.  The image
of this section is the submanifold of $\cG$ given by $t=1$.  
The choice of $g$ also determines a horizontal subspace
$\cH_z \subset T_z\cG$ for each $z\in \cG$, namely $\cH_z =
\operatorname{ker}(dt)_z$.    
In terms of another representative $\wh{g} = e^{2 \Up} g$ of the same
conformal class, we obtain another trivialization of $\cG$, by
identifying 
$$ 
( \wh{t} , x ) \in \R_+ \times M \ \mbox{ with } \ ( {\wh{t}}^2  
{\wh{g}}_{x} , x ) \in \cG \,. 
$$
The two trivializations are then related by the formula 
\begin{equation}\label{triv}
\wh{t}  =  e^{- \Up ( x )} t \mbox {  for  } \ (t , x ) \in \R_+ \times
M. 
\end{equation}

If $(x^1\cdots, x^n)$ are local coordinates on an open set $U$ in $M$, and
if $g$ is given in these coordinates as $g = g_{ij}(x)dx^idx^j$, 
then $(t,x^1\cdots, x^n)$ are local coordinates on $\pi^{-1}(U)$, and ${\bf g}_0$ 
is given by 
$$
{\bf g}_0= t^2 g_{ij}(x)dx^idx^j.  
$$
The horizontal subspace $\cH_z$ is the span of 
$\{\partial_{x^1},\cdots, \partial_{x^n}\}$.  

Consider now the space $\cG \times \R$.  
We write points of $\cG \times \R$ as $(z , \rho)$, with $z \in \cG$, $\rho
\in \R$. 
The dilations $\delta_s$ extend to
$\cG \times \R$ acting in the first factor alone, and we denote these
dilations also by $\delta_s$.  The infinitesimal dilation $T$ also extends
to $\cG \times \R$.
We imbed $\cG$ into $\cG \times \R$ by $\iota : z \mapsto ( z , 0 )$ 
for $z\in \cG$.  Note that $\iota$  commutes with dilations. 
If $g$ is a representative for the
conformal structure with associated fiber coordinate $t$, and if
$(x^1,\cdots, x^n)$ are local coordinates on $M$ as above, then 
$(t, x^1,\cdots, x^n,\rho)$ are local coordinates on $\cG\times \R$.  
We use $0$ to label the $t$-component, $\nf$ to label the
$\rho$-component, lower case Latin letters for $M$, and
capital Latin letters for $\cG\times \R$.  Even without
a choice of coordinates on $M$, we can use $0$, $i$, $\nf$ as 
labels for the components relative to the identification   
$\cG\times \R \simeq \R_+\times M \times \R$ induced by the choice of 
representative metric $g$.  Such an interpretation is coordinate-free and 
global on $M$.  
\begin{definition}\label{ambientdef}
A {\it pre-ambient space} for $(M , [g])$, where $[ g ]$ is a
conformal class of signature $(p , q )$ on $M$, is a pair 
$(\cGt, \gt )$, where:
\begin{enumerate}
\item[(1)] $\cGt$ is a dilation-invariant open neighborhood of $\cG 
\times \{0\}$ in $\cG \times \R$;
\item[(2)] $\gt$ is a smooth metric of signature $(p + 1 , q + 1)$ on
$\cGt$;
\item[(3)] $\gt$ is homogeneous of degree 2 on $\cGt$ (i.e., 
$\delta_s^* \gt =s^2 \gt$, for $s\in \R_+)$; 
\item[(4)] The pullback $\iota^* \gt$ is the tautological tensor 
${\bf g}_0$ on $\cG$.
\end{enumerate}
If $(\cGt,\gt)$ is a pre-ambient space, the metric $\gt$ is called a
{\it pre-ambient metric}.
If the dimension $n$ of $M$ is odd or $n=2$, 
then a pre-ambient space $(\cGt,\gt)$
is called an {\it ambient space} for $(M ,[g] )$ provided we have
\begin{enumerate}
\item[(5)] $\Ric(\gt)$ vanishes to infinite order at every point of  
$\cG \times \{ 0 \}$.  
\end{enumerate}
\end{definition}
\noindent
We prepare to define ambient spaces in the even-dimensional case.
Let $S_{IJ}$ be a symmetric 2-tensor field on an  
open neighborhood of $\cG \times \{ 0 \}$ in $\cG \times  \R$.   
For $m \geq 0$, we write $S_{IJ} = O^+_{IJ}( \rho^m)$ if:
\begin{enumerate}
\item[(i)]
$S_{IJ} = O(\rho^m)$; and
\item[(ii)]
For each point $z\in \cG$, the symmetric 2-tensor $(\iota^*(\rho^{-m}S))(z)$  
is of the form 
$\pi^*s$ for some symmetric 2-tensor $s$ at $x=\pi(z)\in M$ satisfying 
$\operatorname{tr}_{g_x}s = 0$.  The symmetric 2-tensor $s$ is allowed to
depend on $z$, not just on $x$.  
\end{enumerate}
In terms of components relative to a choice of representative metric $g$,
$S_{IJ}= O_{IJ}^+(\rho^m)$ if and only if  
all components satisfy $S_{IJ}= O(\rho^m)$ and if in addition one has that 
$S_{00}$, $S_{0i}$ and $g^{ij}S_{ij}$ are $O(\rho^{m+1})$.  The condition
$S_{IJ}=O^+_{IJ}(\rho^m)$ is easily seen to be preserved by diffeomorphisms  
$\phi$ on a neighborhood of $\cG \times \{ 0 \}$ in $\cG \times \R$
satisfying $\phi|_{\cG \times \{ 0 \}}$ = identity.

Now suppose $( \cGt , \gt)$ is a pre-ambient space for
$(M , [ g ])$, with $n = \dim M$ even and $n\geq 4$.  We say that  
$( \cGt,\gt)$ is an   
{\it ambient space} for $(M , [ g ])$, provided we have 
\begin{enumerate}
\item[($5'$)] $\Ric(\gt)= O_{IJ}^+ (\rho^{n/2 - 1} )$.
\end{enumerate}
If $(\cGt,\gt)$ is an ambient space, the metric $\gt$ is called an 
{\it ambient metric}.

Next, we define a notion of ambient equivalence for pre-ambient spaces. 
\begin{definition}
Let $(\cGt_1,\gt_1)$ and $(\cGt_2, \gt_2)$ 
be two pre-ambient
spaces for $(M , [ g ])$.  We say that $(\cGt_1,\gt_1)$ and 
$(\cGt_2, \gt_2)$ 
are {\it ambient-equivalent} if there exist open sets
$\cU_1 \subset \cGt_1$, $\cU_2 \subset 
\cGt_2$ and a 
diffeomorphism $\phi: \cU_1\rightarrow \cU_2$, with the 
following properties:
\begin{enumerate}
\item[(1)] $\cU_1$ and $\cU_2$ both contain $\cG \times
  \{ 0 \}$;
\item[(2)] $\cU_1$ and $\cU_2$ are dilation-invariant
  and $\phi$ commutes with dilations;
\item[(3)] The restriction of $\phi$ to $\cG \times \{ 0 \}$ is the
identity map; 
\item[(4)] If $n = \dim M$ is odd, then $\gt_1 - \phi^*  
\gt_2$ vanishes to infinite order at every point of $\cG \times
\{ 0 \}$.
\item[($4'$)] If $n = \dim M$ is even, then $\gt_1 - \phi^*  
\gt_2 = O_{IJ}^+ ( \rho^{n/2 })$.
\end{enumerate}
\end{definition}
\noindent
It is easily seen that ambient-equivalence is an equivalence
relation.

One of the main results of this monograph is the following. 

\begin{theorem}\label{main}
Let $(M, [ g ])$ be a smooth manifold of
dimension $n\geq 2$, equipped with a conformal class.  
Then there exists an ambient space for $(M , [g])$. 
Moreover, any two ambient spaces for $(M , [g])$ are ambient-equivalent.  
\end{theorem}
\noindent
It is clear when $n$ is odd that a pre-ambient space is an ambient space
provided it is ambient-equivalent to an ambient space.  The kinds of
arguments we use in Chapter~\ref{formalsect} can be used to show that   
this is also true if $n$ is even and $n\geq 4$.  
Thus for $n>2$, an ambient space for $(M,[g])$ is determined 
precisely up to ambient-equivalence.  This is clearly not true when 
$n=2$:  changing the metric at high finite order generally affects the
infinite-order   
vanishing of the Ricci curvature.  The uniqueness of ambient metrics when
$n=2$ will be clarified in Theorem~\ref{2dim}. 

Our proof of Theorem~\ref{main} will establish an additional important
property 
of ambient metrics.  In Chapter~\ref{formalsect} we will prove the
following two Propositions.
\begin{proposition}\label{straightequiv}
Let $(\cGt,\gt)$ be a pre-ambient space for $(M,[g])$.  There is a
dilation-invariant open set  
$\cU\subset \cGt$ containing $\cG \times \{0\}$ such that the
following three conditions are equivalent.
\begin{enumerate}
\item 
$\tilde{\nabla}T = Id$ on $\cU$. 
\item
$2T\into \gt = d(\|T\|^2)$ on $\cU$.  
\item
For each $p\in \cU$, the parametrized dilation orbit
$s\mapsto \delta_sp$ is a geodesic for $\gt$.
\end{enumerate}
In (1), $\tilde{\nabla}$ denotes the covariant derivative with respect to the
Levi-Civita connection of $\gt$.  So $\tilde\nabla T$ is a $(1,1)$-tensor
on $\cU$, and the requirement is that it be the identity endomorphism at
each point.  In (2), $\|T\|^2 = \gt(T,T)$.  
\end{proposition}
\begin{definition}
A pre-ambient space $(\cGt,\gt)$ for $(M,[g])$ will be said to be {\it
straight} if the equivalent properties of Proposition~\ref{straightequiv}
hold 
with $\cU=\cGt$.  In this case, the pre-ambient metric $\gt$ is also said
to be straight.    
\end{definition}
\noindent
Note that if $(\cGt,\gt)$ is a straight pre-ambient space for $(M,[g])$ and   
$\phi$ is a diffeomorphism of a dilation-invariant open neighborhood $\cU$
of $\cG\times \{0\}$ into $\cGt$ which commutes with dilations and
satisfies that $\phi|_{\cG\times \{0\}}$ is the identity map, then the
pre-ambient space $(\cU,\phi^*\gt)$ is also  straight. 
\begin{proposition}\label{straightprop}
Let $(M, [ g ])$ be a smooth manifold of 
dimension $n> 2$ equipped with a conformal class.
Then there exists a straight ambient space for $(M , [g])$. 
Moreover, if $\gt$ is any ambient metric for $(M,[g])$, 
there is a straight ambient metric $\gt'$ such that if $n$ is odd, then 
$\gt -\gt'$ vanishes to infinite order at $\cG\times \{0\}$, while if $n$
is even, then $\gt -\gt' = O^+_{IJ}(\rho^{n/2})$.  
\end{proposition}
\noindent
Because of Proposition~\ref{straightprop}, one can usually 
restrict attention to straight ambient spaces.  
Observe that the second statement of Proposition~\ref{straightprop} follows
from the first statement, the uniqueness up to ambient-equivalence in
Theorem~\ref{main}, and the 
diffeomorphism-invariance of the straightness condition.  

In this chapter, we will begin the proof of Theorem~\ref{main} by using a
diffeomorphism 
to bring a pre-ambient metric into a normal form relative to a choice of
representative of the conformal class.  Chapter~\ref{formalsect}
will analyze metrics in that normal form and complete the proof
of Theorem~\ref{main}.  We start by formulating the normal form condition.  

\begin{definition}\label{normalformdef}
A pre-ambient space $( \cGt , \gt )$ for 
$(M , [g])$ is said to be in {\it normal form} 
relative to a representative
metric $g$ if the following three conditions hold:
\begin{enumerate}
\item[(1)] For each fixed $z \in \cG$, the set of all $\rho \in
\R$ such that $(z , \rho ) \in \cGt$ is an open interval $I_z$  
containing $0$.
\item[(2)] For each $z \in \cG$, the parametrized curve $I_z
\ni \rho \mapsto ( z , \rho )$ is a geodesic for the metric
$\gt$.  
\item[(3)] Let us write $(t,x,\rho)$ for a point of $\R_+ \times M \times
  \R \simeq \cG \times 
\R$ under the identification induced by $g$, as discussed above.  Then, at
  each point $(t , x , 0) \in \cG \times \{ 0 \}$, the  
metric tensor $\gt$ takes the form
\begin{equation}\label{initialform}
\gt = \,  {\bf g}_0  + 2 t\,  dt\,
d \rho\, .
\end{equation}
\end{enumerate}
\end{definition}

The main result proved in this chapter is the following.
\begin{proposition}\label{normalform}
Let $( M , [ g ])$ be a  smooth manifold
equipped with a conformal class, let $g$ be a representative of the
conformal class, and let $( \cGt , \gt )$ be a
pre-ambient space for $(M , [g])$.  Then there exists a dilation-invariant
open set  
$\cU\subset \cG\times \R$ containing $\cG\times \{0\}$ 
on which there is a unique diffeomorphism $\phi$ from $\cU$ into $\cGt$, such that  
$\phi$ commutes with dilations, $\phi|_{{\cG} \times \{ 0 \}}$ is 
the identity map, and such that the pre-ambient space $(\cU,\phi^* \gt)$ 
is in normal form relative to $g$.  
\end{proposition}
\noindent
Thus, once we have picked a representative $g$ of the conformal class
$[g]$, we can uniquely place any given pre-ambient metric into normal form   
by a diffeomorphism $\phi$.  In 
Proposition~\ref{normalform}, note that 
$(\cU,\phi^*\gt)$ is ambient-equivalent to $(\cGt,\gt)$.    

In Chapter~\ref{formalsect} below, we will establish the following result.  
\begin{theorem}\label{formal}
Let $M$ be a  smooth manifold of dimension $n\geq 2$ and $g$ a
smooth metric on $M$.
\begin{enumerate}
\item[(A)] There exists an ambient space $(\cGt,\gt)$ for $(M , [ g ])$ 
which is in normal form relative to $g$.  
\item[(B)] Suppose that $(\cGt_1,\gt_1)$ and $(\cGt_2,\gt_2)$ are 
two ambient spaces for $(M , [ g ])$, both of which are in normal form
relative to $g$.  
If $n $ is odd, then $\gt_1 - \gt_2$ 
vanishes to infinite order at every point of $\cG \times \{ 0 \}$.  
\noindent
If $n $ is even, then $\gt_1 - \gt_2
= O^+_{IJ}(\rho^{n/2})$.
\end{enumerate}
\end{theorem}
\noindent
{\it Proof of Theorem~\ref{main} using Proposition~\ref{normalform} and
Theorem~\ref{formal}.} 
Given $(M , [g])$, we pick a representative $g$ and invoke
Theorem~\ref{formal}(A).  Thus, there exists an ambient space for $(M ,
[ g ])$.  For the uniqueness, let $(\cGt_1,\gt_1)$ and $(\cGt_2, \gt_2)$ be
ambient spaces for 
$(M , [g])$.  Again, we pick a representative $g$.  Applying 
Proposition~\ref{normalform}, we
find that $(\cGt_1,\gt_1)$ is ambient-equivalent to an ambient space
in normal form relative to $g$.  Similarly, $(\cGt_2,\gt_2)$ 
is ambient-equivalent to an ambient space also in normal form relative to
$g$.  Theorem~\ref{formal}(B) shows that these two ambient spaces in 
normal form are ambient-equivalent.  
Consequently, $(\cGt_1,\gt_1)$ is ambient-equivalent to $(\cGt_2,\gt_2)$. 
\stopthm

The rest of this chapter is devoted to the proof of
Proposition~\ref{normalform}.  We first formulate a notion which will play
a key role in the proof.
Let $(M,[g])$ be a conformal manifold, let $(\cGt,\gt)$ be a pre-ambient
space for $(M,[g])$, and let $g$ be a metric in the conformal class.  
Recall that $g$ determines the fiber coordinate $t:\cG \rightarrow \R_+$
and the  
horizontal subbundle $\cH = \operatorname{ker}(dt)\subset T\cG$.  
For $z\in \cG$, we may view $\cH_z$ as a subspace of $T_{(z,0)}(\cG\times  
\R)$ via the inclusion $\iota:\cG\rightarrow \cG\times \R$.   
For $z\in\cG$, we say that a vector $V\in T_{(z,0)}(\cG\times \R)$ is a  
$g$-transversal for $\gt$ at $(z,0)$ if it satisfies the following
conditions.  
\begin{equation}\label{Vequations}
\begin{split}
&\gt ( V, T ) = t^2\\
&\gt ( V , X) = 0 \mbox{  for all  } X \in \cH_{z}\\
&\gt ( V , V ) = 0.
\end{split}
\end{equation}
In the first line, $t$ denotes the fiber coordinate for the point $z$.  
The motivation for this definition is the observation that 
$V=\pa_{\rho}$ is a $g$-transversal for $\gt$ if $\gt$ 
satisfies condition (3) of Definition~\ref{normalformdef}.  In general we
have the following elementary result.
\begin{lemma}\label{gtransverse}
For each $z\in \cG$,
there exists one and only one $g$-transversal $V_z$ for $\gt$ at $(z,0)$.  
Moreover, $V_z$ is transverse to $\cG\times \{0\}$, $V_z$ depends smoothly
on $z$, and $V_z$ is dilation-invariant in the sense that
$(\delta_s)_* \, V_z = V_{\delta_s(z)}$
for $s\in \R_+$ and $z \in \cG$.  
\end{lemma}
\begin{proof}
Recall the identification $\cG\times \R \simeq \R_+\times M\times \R$ 
determined by $g$.  In terms of the coordinates $(t,x^i,\rho)$ induced by a
choice of local coordinates on $M$, we can express 
$V = V^0\pa_t +V^i\pa_{x^i} +V^\nf \pa_\rho$.  We can also express $\gt$ in 
these coordinates.  By condition (4) of
Definition~\ref{ambientdef}, at a point $(t,x,0)$, $\gt$ takes the
form  
$$
\gt = t^2g_{ij}(x)dx^idx^j + 2\gt_{0\nf}dtd\rho +2\gt_{j\nf}dx^jd\rho 
+\gt_{\nf\nf}(d\rho)^2  
$$
where $\gt_{0\nf}$, $\gt_{j\nf}$ and $\gt_{\nf\nf}$ depend on $(t,x)$.
Nondegeneracy of $\gt$ implies that $\gt_{0\nf}\neq 0$.  
The conditions \eqref{Vequations} defining a $g$-transversal become  
$$
\gt_{0\nf}V^\nf= t, \qquad t^2g_{ij}V^i + \gt_{j\nf}V^\nf = 0,
$$
$$
2\gt_{0\nf}V^0V^\nf +t^2 g_{ij}V^iV^j +2\gt_{i\nf}V^iV^\nf 
+ \gt_{\nf\nf}(V^\nf)^2 = 0. 
$$
It is clear that these equations can be successively solved uniquely for 
$V^\nf$, $V^i$, $V^0$, and the other conclusions of
Lemma~\ref{gtransverse} follow easily from the smoothness and homogeneity 
properties defining a pre-ambient space.  
\end{proof}
\noindent
{\it Proof of Proposition~\ref{normalform}.}  For $z\in \cG$, let $V_z$ be 
the $g$-transversal for $\gt$ at $(z,0)$ given by 
Lemma~\ref{gtransverse}.  Let $\la \mapsto \phi ( z, \la ) \in \cGt$ be a 
(parametrized) geodesic for $\gt$, with initial conditions 
\begin{equation}\label{initialgeod}
\phi ( z, 0 ) \, = \, ( z, 0) \, \qquad\qquad
\pa_\la \phi ( z, \la )|_{\la = 0} \, = \,
V_{z} \, . 
\end{equation}
Since $\gt$ needn't be geodesically complete, $\phi (z , \la)$ is
defined only for $(z, \la)$ in an open neighborhood $\cU_0$ of
$\cG \times \{ 0 \}$ in $\cG \times \R$.  Since
$\gt$ and $V_z$ are homogeneous with respect to the dilations,
we may take $\cU_0$ to be
dilation-invariant. Thus, $\phi : \cU_0 \rightarrow \cGt$ is a smooth
map, commuting with dilations, and satisfying \eqref{initialgeod}.

Since $V_z$ is transverse to $\cG\times \{ 0 \}$,
it follows that $\cU_1 = \{ ( z, \la ) \in \cU_0: \det
\phi^\prime ( z, \la ) \ne 0 \}$ is a dilation-invariant open
neighborhood of $\cG \times \{ 0 \}$ in $\cU_0$.  Thus,
$\phi$ is a local diffeomorphism from $\cU_1$ into $\cGt$, 
commuting with dilations.
Moreover, by definition of $\phi$, we have:
\begin{enumerate}\label{geodinterval}
\item[] Let $z \in \cG$ and let $I$ be an interval
containing $0$.
Assume that $(z, \la ) \in \cU_1$ for all $\la \in I$.  Then $I
\ni \la \mapsto \phi ( z, \la )$ is a geodesic for $\gt$, with
initial conditions \eqref{initialgeod}.
\end{enumerate}
The map $\phi$ need not be globally one-to-one on $\cU_1$.  However,
arguing as in the proof of the Tubular Neighborhood Theorem (see,
e.g. \cite{L}), one concludes that there exists a 
dilation-invariant open neighborhood $\cU_2$ of $\cG \times
\{ 0 \}$ in $\cU_1$, such that $\phi|_{\cU_2}$ is globally one-to-one. 
Thus, $\phi$ is a diffeomorphism from $\cU_2$ to a dilation-invariant open
subset of $\cGt$ containing 
$\cG\times \{ 0 \}$.   

Next, we define $\cU = \{ ( z, \la ) \in \cU_2 :(z,
\mu ) \in \cU_2$ for all $\mu \in \R$ for which $| \mu | \leq | \la | \}$.  
Thus, $\cU$ is a dilation-invariant open neighborhood of $\cG
\times \{ 0 \}$ in $\cU_2$.  Moreover, 
for each fixed $z \in \cG$, 
$
\{\la \in \R: (z, \la ) \in \cU\}
$ 
is an open interval $I_{z}$ containing $0$. 
It follows that for each fixed $z \in \cG$, the  
parametrized curve $I_{z} \ni \la \mapsto \phi ( z, \la)$
is a geodesic for the metric $\gt$.

Since $( \cGt, \gt )$ is a pre-ambient space for
$(M , [g])$, so is $( \cU, \phi^*\gt)$. 
For each fixed $z \in \cG$, the 
parametrized curve $I_{z} \ni \la \mapsto ( z, \la)$ is a
geodesic for $\phi^*\gt$.  
{From} the facts that $V$ satisfies \eqref{Vequations} and $\phi$ satisfies 
\eqref{initialgeod}, it follows that 
under the identification 
$\R_+ \times M \times \R \simeq \cG \times 
\R$ induced by $g$, we have at $\la = 0$:
\[
\begin{split}
&(\phi^*\gt) (\pa_{\la} , T ) = t^2,\\
&(\phi^*\gt) ( \pa_\la , X) = 0 \mbox{  for  } X \in TM, \\  
&(\phi^*\gt) ( \pa_\la ,\pa_\la) = 0.  
\end{split}
\]
Together with property (4) of Defintion~\ref{ambientdef} of the pre-ambient
space $( \cGt , \gt)$, these equations show that  
$\phi^*\gt =  {\bf g}_0  + 2t \, dt \, d \la $ when $\la = 0$. 
This establishes the existence part of Proposition~\ref{normalform}.  
The uniqueness follows from the fact that the above construction of 
$\phi$ is forced.  If $\phi$ is any diffeomorphism with the required
properties, then at $\rho = 0$, $\phi_*(\pa_\la)$ is a $g$-transversal for
$\gt$, so must be $V$.  Then for $z\in \cG$, the curve 
$I_{z} \ni \la \mapsto \phi (z, \la )$ must be the unique geodesic
satisfying the initial conditions \eqref{initialgeod}.  These requirements 
uniquely determine $\phi$ on $\cU$.  
\stopthm

\section{Formal theory}\label{formalsect}
The first goal of this chapter is to prove Theorem~\ref{formal} for $n>2$.  
We begin with the following Lemma.
\begin{lemma}\label{explicitnormal}
Let $(\cGt,\gt)$ be a pre-ambient space for $(M,[g])$, where $\cGt$
has the property that for  
each $z\in \cG$, the set of all $\rho \in \R$ such that $(z,\rho)\in 
\cGt$ is an open interval $I_{z}$ containing $0$.  Let $g$ be a metric in
the conformal class, with associated identification 
$\R_+\times M\times \R \simeq \cG\times \R$.  Then $(\cGt,\gt)$ is in
normal form relative to $g$ if and only if one has on $\cGt$:
\begin{equation}\label{geocond}
\gt_{0\nf} = t,\qquad \gt_{i\nf} = 0,\qquad \gt_{\nf\nf}=0.
\end{equation}
\end{lemma}
\begin{proof}
Since a pre-ambient metric satisfies $\iota^*\gt={\bf g}_0$, if $\gt$ 
satisfies \eqref{geocond}, then it has the form \eqref{initialform} at
$\rho =0$.  Thus we must show that for $\gt$ satisfying \eqref{initialform}
at $\rho =0$, the condition that the 
lines $\rho\ni I_{(t,x)}\rightarrow (t,x,\rho)$ are geodesics for 
$\gt$ is equivalent to \eqref{geocond}.  Now the $\rho$-lines are geodesics 
if and only if 
$\Gat_{\nf\nf I} = 0$, where
$\Gat_{IJK} = \gt_{KL}\Gat_{IJ}^L$ and $\Gat_{IJ}^L$ are the usual
Christoffel symbols for $\gt$.  Taking $I=\nf$ gives $\pa_\rho\gt_{\nf\nf}
= 0$, which combined with $\gt_{\nf\nf}|_{\rho  = 0} 
= 0$ from \eqref{initialform} yields $\gt_{\nf\nf}=0$.  Now taking 
$I=i$ and $I=0$ and using \eqref{initialform} gives $\gt_{\nf i}=0$ and
$\gt_{\nf 0}=t$.   
\end{proof}

The case $n=2$ is exceptional for Theorem~\ref{formal}.  We give the proof
for $n>2$ now; a sharpened version for $n=2$ will be given in
Theorem~\ref{2dim}.   
\medskip

\noindent
{\it Proof of Theorem~\ref{formal} for $n>2$.} Given $M$ and a metric $g$
on $M$,  
we must construct a smooth metric $\gt$ on a suitable neighborhood $\cGt$
of $\R_+\times M\times \{0\}$ with the properties:   
\begin{enumerate}
\item[(1)]
$\delta_s^*\gt = s^2\gt,\quad s>0$
\item[(2)]
$\gt = t^2 g(x) +2t\,dt\,d\rho$ when $\rho = 0$ 
\item[(3)]
For each $(t,x)$, the curve $\rho\rightarrow (t,x,\rho)$ is a
geodesic for $\gt$  
\item[(4)]
If $n$ is odd, then $\Ric(\gt)=0$ to infinite order at $\rho = 0$   
\item[$(4')$]
If $n$ is even, then $\Ric(\gt) = O^+_{IJ}(\rho^{n/2-1})$.   
\end{enumerate}
Also we must show that if $n$ is odd, 
then such a metric is uniquely determined to infinite order at $\rho = 0$,
while if $n$ even, then it is determined modulo $O^+_{IJ}(\rho^{n/2})$. 

Lemma~\ref{explicitnormal} enables us to replace (3) above with
\eqref{geocond}.   
Thus the components $\gt_{I\nf}$ are determined.  We consider the remaining
components of $\gt$ as unknowns, subject to the homogeneity conditions
determined by (1) above.  Then (2) can be interpreted as initial conditions
and the equation $\Ric(\gt)=0$ as a system of partial differential
equations to be solved formally.  

Set $\gt_{00}=a$, $\gt_{0i} = tb_i$ and $\gt_{ij} = t^2g_{ij}$, where 
all of $a$, $b_i$, $g_{ij}$ are functions of $(x,\rho)$.  Condition (2)
gives at $\rho = 0$:  $a=0$, $b_i=0$ and $g_{ij}$ is the given metric. 
In order to determine the first derivatives of $a, b_i, g_{ij}$ at
$\rho = 0$, we calculate at $\rho = 0$ the components
$\Rt_{IJ}=\Ric_{IJ}(\gt)$ for  
$I,J \neq \nf$.  This is straightforward but tedious.  We have 
\begin{equation}\label{inverse}
\gt^{IJ} = 
\left(
\begin{matrix}
0&0&t^{-1}\\
0&t^{-2}g^{ij}&-t^{-2}b^i\\
t^{-1}&-t^{-2}b^j&t^{-2}(b_kb^k -a)
\end{matrix}
\right)
\end{equation}
and in particular at $\rho = 0$ we have
$$
\gt^{IJ} = 
\left(
\begin{matrix}
0&0&t^{-1}\\
0&t^{-2}g^{ij}&0\\
t^{-1}&0&0
\end{matrix}
\right).
$$
The Christoffel symbols are given by:

\begin{gather}\label{c}
\begin{gathered}
2\Gat_{IJ0} = 
\left(
\begin{matrix}
0&\pa_ja&\pa_\rho a\\
\pa_ia&t(\pa_jb_i+\pa_ib_j-2g_{ij})&t\pa_\rho b_i\\
\pa_\rho a&t\pa_\rho b_j&0
\end{matrix}
\right)\\
2\Gat_{IJk} = 
\left(
\begin{matrix}
2b_k-\pa_k a&t(\pa_jb_k-\pa_kb_j+2g_{jk})&t\pa_\rho b_k\\
t(\pa_ib_k-\pa_kb_i+2g_{ik})&2t^2\Ga_{ijk}&t^2\pa_\rho g_{ik}\\ 
t\pa_\rho b_k&t^2\pa_\rho g_{jk}&0
\end{matrix}
\right)\\
2\Gat_{IJ\nf} = 
\left(
\begin{matrix}
2-\pa_\rho a&-t\pa_\rho b_j&0\\
-t\pa_\rho b_i&-t^2\pa_\rho g_{ij}&0\\
0&0&0
\end{matrix}
\right),
\end{gathered}
\end{gather}
where in the second equation the $\Ga_{ijk}$ refers to the Christoffel
symbol of the metric $g_{ij}$ with $\rho$ fixed.
The Ricci curvature is given by:
\begin{equation}\label{Ricci}
\begin{split}
\Rt_{IJ} =  \tfrac12 \gt^{KL}
\left(\pa^2_{IL}\gt_{JK}\right.&\left. +\pa^2_{JK}\gt_{IL}
-\pa^2_{KL}\gt_{IJ}-\pa^2_{IJ}\gt_{KL}\right) \\
& +\gt^{KL}\gt^{PQ}\left(\Gat_{ILP}\Gat_{JKQ}-\Gat_{IJP}\Gat_{KLQ}\right). 
\end{split}
\end{equation}
Computing this using the above gives at $\rho = 0$:
\begin{gather}\label{firstorder}
\begin{gathered}
\Rt_{00} = \frac{n}{2t^2} \left(2-\pa_\rho a\right),\\
\Rt_{i0}=\frac1{2t}\left(\pa^2_{i\rho}a - n\pa_\rho b_i\right),\\
\Rt_{ij}=\tfrac12 \left[(\pa_\rho a -n)\pa_\rho g_{ij}
-g^{kl}\pa_\rho g_{kl} g_{ij} 
+\nabla_j\pa_\rho b_i +\nabla_i\pa_\rho b_j -\pa_\rho b_i \pa_\rho
b_j\right] +R_{ij},
\end{gathered}
\end{gather}
where in the last equation $R_{ij}$ refers to the Ricci curvature of the
initial metric and $\nabla_j$ denotes the covariant derivative 
with respect to its Levi-Civita connection.  
Setting these to $0$ successively shows that the vanishing of these
components of $\Rt_{IJ}$ at $\rho=0$ is equivalent to the conditions:
\begin{gather}\label{initial}
\begin{gathered}
a = 2\rho +O(\rho^2),  \qquad b_i = O(\rho^2),\\
g_{ij}(x,\rho )= g_{ij}(x) + 2P_{ij}(x)\rho +O(\rho^2),
\end{gathered}
\end{gather}
where 
\begin{equation}\label{Ptensor}
P_{ij}=(n-2)^{-1}\left(R_{ij} - \frac{R}{2(n-1)}g_{ij}\right). 
\end{equation}

Next we carry out an inductive perturbation calculation for higher orders.
Suppose for some $m\geq 2$ that $\gt^{(m-1)}_{IJ}$ is a metric satisfying 
\eqref{geocond}, \eqref{initial}.
Set $\gt^{(m)}_{IJ} = \gt^{(m-1)}_{IJ} 
+ \Phi_{IJ}$, where 
\begin{equation}\label{phi}
\Phi_{IJ}= \rho^m
\left(
\begin{matrix}
\phi_{00}&t\phi_{0j}&0\\
t\phi_{i0}&t^2\phi_{ij}&0\\
0&0&0
\end{matrix}
\right)
\end{equation}
and the $\phi_{IJ}$ are functions of $(x,\rho)$.
{From} \eqref{Ricci} it follows that
\begin{equation}\label{perturb}
\begin{split}
\Rt^{(m)}_{IJ}&=\Rt^{(m-1)}_{IJ} 
+\tfrac12 \gt^{KL}
\left(\pa^2_{IL}\Phi_{JK}+\pa^2_{JK}\Phi_{IL}
-\pa^2_{KL}\Phi_{IJ}-\pa^2_{IJ}\Phi_{KL}\right)\\  
+&\gt^{KL}\gt^{PQ}\left(\Gat_{ILP}\Ga^\Phi_{JKQ} + \Ga_{ILP}^\Phi\Gat_{JKQ} 
-\Gat_{IJP}\Ga_{KLQ}^\Phi - \Ga_{IJP}^\Phi\Gat_{KLQ}\right)+O(\rho^m),  
\end{split}
\end{equation}
where $\gt^{AB}$ and $\Gat_{ABC}$ refer to the metric $\gt^{(m)}_{IJ}$,  
and $2\Ga^\Phi_{IJK} = \pa_J \Phi_{IK}+\pa_I\Phi_{JK} - \pa_K\Phi_{IJ}$.  
These are given modulo $O(\rho^m)$ by:
\begin{gather}\label{cP}
\begin{gathered}
2\Ga^\Phi_{IJ0} = 
\left(
\begin{matrix}
0&0&\pa_\rho \Phi_{00}\\
0&0&\pa_\rho \Phi_{i0}\\
\pa_\rho \Phi_{00}&\pa_\rho \Phi_{0j}&0
\end{matrix}
\right)\\
2\Ga^\Phi_{IJk} = 
\left(
\begin{matrix}
0&0&\pa_\rho \Phi_{0k}\\
0&0&\pa_\rho \Phi_{ik}\\ 
\pa_\rho \Phi_{0k}&\pa_\rho \Phi_{jk}&0
\end{matrix}
\right)\\
2\Ga^\Phi_{IJ\nf} = 
\left(
\begin{matrix}
-\pa_\rho \Phi_{00}&-\pa_\rho \Phi_{0j}&0\\
-\pa_\rho \Phi_{i0}&-\pa_\rho \Phi_{ij}&0\\
0&0&0
\end{matrix}
\right).
\end{gathered}
\end{gather}
On the right hand side of \eqref{perturb},
one can take for $\gt^{AB}$ and $\Gat_{ABC}$ the quantities obtained by
substituting  
\eqref{initial} into \eqref{inverse}, \eqref{c}. 
Calculating, one finds 
\begin{equation}\label{perturbricci}
\begin{split}
t^2\Rt^{(m)}_{00} &= t^2\Rt^{(m-1)}_{00} +
m(m-1-\tfrac{n}{2})\rho^{m-1}\phi_{00} +O(\rho^m),\\ 
t\Rt^{(m)}_{0i} &= t\Rt^{(m-1)}_{0i} +
m(m-1-\tfrac{n}{2})\rho^{m-1}\phi_{0i}  
+ \tfrac{m}{2}\rho^{m-1}\pa_i\phi_{00}+O(\rho^m),\\
\Rt^{(m)}_{ij} &= \Rt^{(m-1)}_{ij} +
m \rho^{m-1}\left[(m-\tfrac{n}{2})\phi_{ij}-\tfrac12
  g^{kl}\phi_{kl}g_{ij}\right.\\   
&\qquad\qquad \left. +\tfrac12(\nabla_j\phi_{0i} +
\nabla_i\phi_{0j}) +P_{ij}\phi_{00}\right] +O(\rho^m)\\  
t\Rt^{(m)}_{0\nf} &= t\Rt^{(m-1)}_{0\nf} +
m(m-1)\rho^{m-2}\phi_{00}/2 +O(\rho^{m-1})\\
\Rt^{(m)}_{i\nf} &= \Rt^{(m-1)}_{i\nf} +
m(m-1)\rho^{m-2}\phi_{0i}/2 +O(\rho^{m-1})\\
\Rt^{(m)}_{\nf\nf} &= \Rt^{(m-1)}_{\nf\nf} -
m(m-1)\rho^{m-2}g^{kl}\phi_{kl}/2 +O(\rho^{m-1}).
\end{split}
\end{equation}

We first consider only the components $\Rt_{IJ}$ with $I$, $J\neq \nf$.
Suppose inductively that $\gt^{(m-1)}$ has been determined so that
$\Ric_{IJ}(\gt^{(m-1)}) = O(\rho^{m-1})$ for 
$I,J\neq \nf$, and that $\gt^{(m-1)}$ is uniquely determined modulo
$O(\rho^m)$ by this condition and \eqref{geocond}.  Define $\gt^{(m)}$ as 
above.
If $n$ is odd or if $n$ is even and $m\leq n/2$, then the
coefficient $m-1-n/2$ appearing in the first two formulae of 
\eqref{perturbricci} does not vanish, so one can  
uniquely choose $\phi_{00}$ and $\phi_{0i}$ at $\rho = 0$ to make
$\Rt_{00}^{(m)}$ and $\Rt_{0i}^{(m)}$ be $O(\rho^m)$.  The map  
$\phi_{ij} \rightarrow (m-n/2)\phi_{ij} - \frac12 g^{kl}\phi_{kl} g_{ij}$ is
bijective on symmetric 2-tensors unless $m=n/2$ or $m=n$, so except for
these $m$ one can similarly make $\Rt_{ij}^{(m)}=O(\rho^m)$.  
Thus the induction proceeds up to $m<n/2$ for $n$ even and up to 
$m<n$ for $n$ odd.  Consider the next value of $m$ in each case.  
For $n$ even and $m=n/2$, one can
uniquely determine $\phi_{00}$ and $\phi_{0i}$ at $\rho = 0$ to make 
$\Rt_{00}$, $\Rt_{0i} = O(\rho^{n/2})$ just as before.  In addition, one
can choose 
$g^{ij}\phi_{ij}$ to guarantee that $g^{ij}\Rt_{ij}=O(\rho^{n/2})$. 
Thus for $n$ even, we deduce that $\gt_{IJ} \mod O_{IJ}^+(\rho^{n/2})$ 
is uniquely determined by the condition $\Rt_{IJ} = O_{IJ}^+(\rho^{n/2-1})$  
for $I,\,J \neq \nf$.  
For $n$ odd and $m=n$, one can again uniquely 
determine $\phi_{00}$ and $\phi_{0i}$ at $\rho = 0$ to make 
$\Rt_{00}$, $\Rt_{0i} = O(\rho^n)$, and now can uniquely determine 
the trace-free part of $\phi_{ij}$ to guarantee that 
$\Rt_{ij} = \lambda g_{ij} \rho^{n-1} \mod O(\rho^n)$ for some function
$\lambda$.  

In order to analyze 
the remaining components $\Rt_{I\nf}$ and to complete the analysis above in
the case $m=n$ when 
$n$ is odd, we consider the contracted Bianchi identity. 
The Ricci curvature of $\gt$ satisfies the Bianchi identity
$\gt^{JK}\nt_I \Rt_{JK} = 2\gt^{JK}\nt_J\Rt_{IK}$.  Writing this in terms
of coordinate derivatives gives
\begin{equation}\label{bianchi}
2\gt^{JK}\pa_J\Rt_{IK} -\gt^{JK}\pa_I\Rt_{JK} - 
2\gt^{JK}\gt^{PQ}\Gat_{JKP}\Rt_{QI} =0.
\end{equation}
Suppose for some $m\geq 2$ that $\Rt_{IJ}=O(\rho^{m-1})$ for 
$I,\,J \neq \nf$ and $\Rt_{I\nf}=O(\rho^{m-2})$.  Write out
\eqref{bianchi} for 
$I=0,i,\nf$, using \eqref{inverse}, \eqref{c}, \eqref{initial} and 
the homogeneity of the components $\Rt_{IJ}$.  Calculating mod
$O(\rho^{m-1})$, one obtains:
\begin{gather}\label{bcomp}
\begin{gathered}
(n-2-2\rho\pa_{\rho})\Rt_{0\nf} + t\pa_{\rho}\Rt_{00}= O(\rho^{m-1}),\\ 
(n-2-2\rho\pa_{\rho})\Rt_{i\nf} - t\pa_i\Rt_{0\nf} +t\pa_\rho \Rt_{i0}= 
O(\rho^{m-1}),\\
(n-2-\rho\pa_{\rho})\Rt_{\nf\nf} +g^{jk}\nabla_j\Rt_{k\nf} 
+tP_k{}^k\Rt_{0\nf} -\tfrac12 g^{jk}\pa_{\rho}\Rt_{jk}= O(\rho^{m-1}).
\end{gathered}
\end{gather}

Suppose first that $n$ is even.  Let $\gt_{IJ}$ be the metric 
determined above mod $O_{IJ}^+(\rho^{n/2})$ by the requirement $\Rt_{IJ} = 
O_{IJ}^+(\rho^{n/2-1})$ for $I,\,J \neq \nf$.  We show by induction on $m$
that $\Rt_{I\nf} = O(\rho^{m-1})$ for $1\leq m \leq n/2$.  
The statement is clearly true 
for $m=1$.  Suppose it holds for $m-1$ and write 
$\Rt_{I\nf} = \gamma_I \rho^{m-2}$.  The hypotheses for \eqref{bcomp} are
satisfied.  The first 
equation of \eqref{bcomp} gives $(n+2-2m)\gamma_0 = O(\rho)$, so 
$\Rt_{0\nf} = O(\rho^{m-1})$.  The second equation of \eqref{bcomp} then
gives $(n+2-2m)\gamma_i = O(\rho)$, so $\Rt_{i\nf} = O(\rho^{m-1})$.  
The last equation then gives $(n-m)\gamma_\nf = O(\rho)$, so
$\Rt_{\nf\nf} = O(\rho^{m-1})$, completing the induction.  

Hence, if $n$ is even, we have uniquely determined 
$\gt_{IJ} \mod O_{IJ}^+(\rho^{n/2})$ so that $\Ric(\gt) = 
O_{IJ}^+(\rho^{n/2-1})$.  A finite order Taylor polynomial for 
$\gt$ will be nondegenerate on a neighborhood $\cGt$ of 
$\R_+\times M\times \{0\}$ satisfying the required properties.  
This concludes the proof of Theorem~\ref{formal} for $n$ even.  

If $n$ is odd, let $\gt_{IJ}$ be the metric determined above 
by the requirements $\Rt_{00}$, $\Rt_{0i} = O(\rho^n)$
and $\Rt_{ij} = \lambda g_{ij} \rho^{n-1} + O(\rho^n)$.  Then $\gt_{IJ}$
is uniquely determined up to $O(\rho^{n+1})$ except that $\gt_{ij}$ has an 
additional indeterminacy of the form $c g_{ij}\rho^n$  for some function
$c$.  As for $n$ even, consider the induction based on \eqref{bcomp}.
The constants $n+2-2m$ never vanish for $n$ odd and $m$ integral, and one
may proceed with 
the induction for all three equations to conclude that
$\Rt_{I\nf} = O(\rho^{n-2})$.  The hypotheses of \eqref{bcomp} now hold
with $m=n$.  The first two equations give $\Rt_{0\nf}$, 
$\Rt_{i\nf}= O(\rho^{n-1})$.  In the third equation, the coefficient of
$\gamma_\nf \rho^{n-2}$ is 
now 0, so the third equation reduces to $\lambda \rho^{n-2} =
O(\rho^{n-1})$.  We conclude that $\lambda = O(\rho)$, i.e. 
$\Rt_{ij} = O(\rho^n)$.  However, there is still the indeterminacy of 
$c g_{ij}\rho^n$ in
$\gt_{ij}$ and we do not yet know that $\Rt_{\nf\nf}=O(\rho^{n-1})$.  
These can be dealt with simultaneously by 
observing directly that one can uniquely choose $c$ at $\rho = 0$
to make $\Rt_{\nf\nf}=O(\rho^{n-1})$.  Namely, set
$\gt'_{IJ}= \gt_{IJ} +\Phi_{IJ}$, where $\Phi_{IJ}$ is given by \eqref{phi} 
with $m=n$ and $\phi_{00}=\phi_{0i}=0$, $\phi_{ij}=cg_{ij}$.  According to 
the last formula of \eqref{perturbricci}, one has
$
\Rt'_{\nf\nf} = \Rt_{\nf\nf} -n^2(n-1)c\rho^{n-2}/2 +O(\rho^{n-1}).
$
Therefore, $c \mod O(\rho)$ is uniquely determined by the requirement  
$\Rt'_{\nf\nf} = O(\rho^{n-1})$.  Removing the $'$, we thus have  
$\gt_{IJ}$ uniquely determined $\mod O(\rho^{n+1})$ by the conditions
$\Rt_{IJ}=O(\rho^n)$ for $I$, $J \neq \nf$ and $\Rt_{I\nf}=O(\rho^{n-1})$.   
We can now proceed inductively to all higher orders with no problems:
\eqref{perturbricci} shows that the requirement 
$\Rt_{IJ}=0$ to infinite order for $I$, $J \neq \nf$ uniquely determines 
$\gt_{IJ}$ to infinite order, and \eqref{bcomp} shows that the so
determined $\gt_{IJ}$ also satisfies $\Rt_{I\nf}=0$ to infinite order. 

Summarizing, the components $\gt_{I\nf}$ are given by \eqref{geocond}.  We 
have determined all derivatives of the components $\gt_{IJ}$   
for $I$, $J\neq \nf$ at $\rho = 0$
uniquely to ensure that all components of $\Ric(\gt)$ vanish to infinite
order.  By Borel's Theorem, we can find 
a homogeneous symmetric 2-tensor on a   
neighborhood of $\rho = 0$ with the prescribed Taylor expansion.  We can
then choose a dilation-invariant subneighborhood $\cGt$    
satisfying condition (1) of Definition~\ref{normalformdef} on 
which this tensor is nondegenerate with signature $(p+1,q+1)$.
\stopthm

Next we show that the metric $\gt_{IJ}$ of Theorem~\ref{formal} takes a   
special form.
\begin{lemma}\label{0comp}
Let $n\geq 2$.  If $\gt$ has the form
\begin{equation}\label{spform}
\gt_{IJ} = 
\left(
\begin{matrix}
2\rho&0&t\\
0&t^2g_{ij}&0\\
t&0&0
\end{matrix}
\right),
\end{equation}
where $g_{ij}=g_{ij}(x,\rho)$ is a one-parameter family of metrics on
$M$, then the Ricci curvature of $\gt$ satisfies 
$\Rt_{0I}=0$.
\end{lemma}
\begin{proof}
For $\gt$ of the form \eqref{spform}, the Christoffel symbols \eqref{c}
become:  
\begin{gather}\label{cn}
\begin{gathered}
\Gat_{IJ0} = 
\left(
\begin{matrix}
0&0&1\\
0&-tg_{ij}&0\\
1&0&0
\end{matrix}
\right)\\
\Gat_{IJk} = 
\left(
\begin{matrix}
0&tg_{jk}&0\\
tg_{ik}&t^2\Ga_{ijk}&\frac12 t^2g_{ik}'\\ 
0&\frac12 t^2g_{jk}'&0
\end{matrix}
\right)\\
\Gat_{IJ\nf} = 
\left(
\begin{matrix}
0&0&0\\
0&-\frac12t^2g_{ij}'&0\\
0&0&0
\end{matrix}
\right).
\end{gathered}
\end{gather}
Here $'$ denotes $\pa_\rho$.  
Lemma~\ref{0comp} is a straightforward    
computation from \eqref{Ricci} using \eqref{inverse} and \eqref{cn}.  
Alternate derivations are given below in the comments after 
Proposition~\ref{curvrelation} and in the proof of 
Proposition~\ref{Tcontract}.
\end{proof}
For future reference we record the
raised index version of \eqref{cn}:
\begin{gather}\label{cnr}
\begin{gathered}
\Gat_{IJ}^0 = 
\left(
\begin{matrix}
0&0&0\\
0&-\frac12 tg_{ij}'&0\\
0&0&0
\end{matrix}
\right)\\
\Gat_{IJ}^k = 
\left(
\begin{matrix}
0&t^{-1}\delta_j{}^k&0\\
t^{-1}\delta_i{}^k&\Ga_{ij}^k&\frac12 g^{kl}g_{il}'\\ 
0&\frac12 g^{kl}g_{jl}'&0 
\end{matrix}
\right)\\
\Gat_{IJ}^\nf = 
\left(
\begin{matrix}
0&0&t^{-1}\\
0&-g_{ij} +\rho g_{ij}'&0\\
t^{-1}&0&0
\end{matrix}
\right).
\end{gathered}
\end{gather}

\begin{proposition}\label{simpler}  Let $n>2$.  
The metric $\gt_{IJ}$ in Theorem~\ref{formal} satisfies
$\gt_{00}=2\rho$ and $\gt_{0i}=0$, to infinite order for $n$ odd, and
modulo $O(\rho^{n/2+1})$ for $n$ even.
\end{proposition}
\begin{proof}
Consider the inductive construction of the metric $\gt$ in the proof of  
Theorem~\ref{formal}.  At each step of the induction, the perturbation
terms $\phi_{00}$ and $\phi_{0i}$ at $\rho = 0$ are determined from the
first two equations of \eqref{perturbricci}.  If $\gt^{(m-1)}$ is of the
form \eqref{spform}, then by Lemma~\ref{0comp} we have 
$\Rt^{(m-1)}_{00}$, $\Rt^{(m-1)}_{0i} = 0$, so we obtain 
$\phi_{00}$, $\phi_{0i} = O(\rho)$.  It therefore follows by induction that 
$\gt_{00}=2\rho$ and $\gt_{0i}=0$ to all orders for which these components
are determined. 
\end{proof}

The next result shows that the special form given by 
Proposition~\ref{simpler} can be reinterpreted in terms of the conditions
of Proposition~\ref{straightequiv}.

\begin{proposition}\label{tgeodesics} Suppose $n\geq 2$.  
Let the pre-ambient space $(\cGt,\gt)$ be in normal form relative to a
representative metric $g$.  The following conditions are equivalent. 
\begin{enumerate}
\item
$\gt_{00} = 2\rho$ and $\gt_{0i} = 0$.
\item
For each $p\in \cGt$, the dilation orbit $s\rightarrow \delta_s p$ is a
geodesic for $\gt$.
\item
$2T\into \gt = d(\|T\|^2)$
\item
The infinitesimal dilation field $T$ satisfies ${\tilde \nabla}T =Id$. 
\end{enumerate}
\end{proposition}
\begin{proof}
By Lemma~\ref{explicitnormal}, we have \eqref{geocond}.  The computations
leading to \eqref{c} are valid for all $n\geq 2$, so 
the Christoffel symbols are given by \eqref{c}.    
Consider a dilation orbit $s\rightarrow \delta_s p$ for
$p\in \cGt$, given in components by $s\rightarrow (st,x,\rho)$. 
Its tangent vector is a constant multiple of $\pa_t$.  These orbits are
therefore geodesics for $\gt$ if and only if 
$\Gat_{00I}=0$.  From \eqref{c} and the initial normalization  
\eqref{initialform}, it is easily seen that the condition $\Gat_{00I}=0$ is
equivalent to  
$\gt_{00}=2\rho$, $\gt_{0i}=0$.  Therefore (1) is equivalent to (2).  
If $\gt$ is any pre-ambient metric, then (3)
states that $2t\gt_{I0} = \pa_I(t^2\gt_{00})$, while the condition 
$\Gat_{00I}=0$ can be written $2\pa_0\gt_{0I}=\pa_I\gt_{00}$.  Both of
these are automatic for $I=0$ and are easily seen to be equivalent for
$I=i$, $\infty$ upon using the 
homogeneity of the components $\gt_{0I}$.  Thus (2) and (3) are equivalent.   
As for (4), if ${\tilde \nabla}T =Id$, then  
${\tilde \nabla}_T T =T$, which is a restatement of the 
geodesic condition.  On the other hand, if $\gt$ has the form 
\eqref{spform}, then \eqref{cnr} gives ${\tilde \nabla}_I T^J =  
\delta^J{}_I$.    
\end{proof}
\noindent
{\it Proof of Proposition~\ref{straightequiv}.}
Choose a representative metric $g$ and invoke Proposition~\ref{normalform}.   
If $\cU$ and $\phi$ are as in Proposition~\ref{normalform}, then the
pre-ambient space $(\cU,\phi^*\gt)$ is in normal form relative to $g$.  
By Proposition~\ref{tgeodesics}, the conditions (1)-(3) of 
Proposition~\ref{straightequiv} are equivalent for the metric $\phi^*\gt$
on $\cU$.  But all of these conditions are invariant under $\phi$, so they 
are also equivalent for $\gt$ on $\phi(\cU)$.   
\stopthm

\noindent
{\it Proof of Proposition~\ref{straightprop}}.  As pointed out after the
statement of Proposition~\ref{straightprop}, we only need to prove the 
existence of a straight ambient space for $(M,[g])$.  Choose a
representative metric $g$.  According to Proposition~\ref{simpler}, we may
as well take the ambient metric given by Theorem~\ref{formal} to be of 
the form \eqref{spform}.  Then Proposition~\ref{tgeodesics} shows that upon
choosing $\cGt$ suitably, the ambient space $(\cGt,\gt)$ is straight.   
\stopthm

Let us return to consider the formal determination of $\gt_{IJ}$.  
According to Proposition~\ref{simpler} and Lemma~\ref{0comp}, $\gt$ may be
taken of the form 
\eqref{spform}, and for any such $\gt$ the 
equations $\Rt_{0I}=0$ already hold.  Therefore 
$\gt_{ij}=t^2g_{ij}(x,\rho)$ can be regarded as the only ``unknown''
component.  One can 
calculate the remaining components of $\Rt_{IJ}$ 
to obtain explicit equations for $g_{ij}(x,\rho)$.  Again calculating from 
\eqref{Ricci}, one finds
\begin{gather}\label{evol}
\begin{gathered}
\Rt_{ij} = \rho g_{ij}'' -\rho g^{kl}g_{ik}'g_{jl}' +\tfrac12 \rho
g^{kl}g_{kl}'g_{ij}' 
- (\tfrac{n}{2}-1) g_{ij}' - \tfrac12 g^{kl}g_{kl}'g_{ij} +R_{ij} \\ 
\Rt_{i\nf} = \tfrac12 g^{kl}(\nabla_k g_{il}' - \nabla_i g_{kl}')\\
\Rt_{\nf\nf}= -\tfrac12g^{kl}g_{kl}'' +\tfrac14
g^{kl}g^{pq}g_{kp}'g_{lq}'. 
\end{gathered}
\end{gather}
Here $R_{ij}$ and $\nabla$ denote the Ricci curvature and Levi-Civita
connection  
of $g_{ij}(x,\rho)$ with $\rho$ fixed.  The Taylor expansion of 
$g_{ij}(x,\rho)$ can be determined by successively 
differentiating and evaluating at $\rho = 0$ the equations obtained by  
setting these expressions to 0.  For example, simply evaluating the first
equation at $\rho = 0$ recovers the fact that $g_{ij}'|_{\rho = 0} =
2P_{ij}$, which we obtained in \eqref{initial}.  According to the proof of  
Theorem~\ref{formal}, for $n$ even the first equation of \eqref{evol}
determines the derivatives $\pa_{\rho}^m g_{ij}$ 
for $m<n/2$ and also $g^{ij}\pa_{\rho}^{n/2}g_{ij}$, and then the
second and third equations automatically hold mod $O(\rho^{n/2 -1})$.  
In practice, it is easier to calculate the traces
$g^{ij}\pa_{\rho}^mg_{ij}$ for $m\leq n/2$ 
using the last equation rather than the first.  For 
$n$ odd, the first equation determines the derivatives $\pa_{\rho}^m
g_{ij}$ for $m<n$ as well as the trace-free part of $\pa_{\rho}^n g_{ij}$.
The trace part of the first equation at order $n$ is automatically true.
The value of $g^{ij}\pa_{\rho}^n g_{ij}$ is determined by the third 
equation, and all higher derivatives are then determined by the first
equation.  

If the initial metric is Einstein, one can identify explicitly the solution 
$g_{ij}(x,\rho)$.  It is straightforward to check that if the initial
metric satisfies $R_{ij} = 2\lambda (n-1) g_{ij}$, then 
$g_{ij}(x,\rho) = (1+\lambda \rho)^2 g_{ij}(x)$ solves \eqref{evol}.  

In general it is feasible to carry out the first few iterations by hand.
One finds at $\rho = 0$:
\begin{equation}\label{deriv}
\begin{split}
\tfrac12 (n-4)g_{ij}'' &= -B_{ij} +(n-4)P_i{}^kP_{jk},\qquad
  n\neq 4\\
\tfrac12 (n-4)(n-6)g_{ij}''' &=  
 B_{ij},_k{}^k - 2W_{kijl}B^{kl} -4(n-6)P_{k(i}B_{j)}{}^k-4P_k{}^kB_{ij}\\
  +4(n-4)&P^{kl}C_{(ij)k},_l
 -2(n-4)C^k{}_i{}^lC_{ljk}+(n-4)C_i{}^{kl}C_{jkl}\\  
  +2(n-4)&P^k{}_k,_lC_{(ij)}{}^l
-2(n-4)W_{kijl}P^k{}_mP^{ml},\qquad n\neq 4,6,
\end{split}
\end{equation}
where $W_{ijkl}$ is the Weyl tensor,
$$
C_{ijk} = P_{ij},_k-P_{ik},_j
$$
is the Cotton tensor, and 
$$
B_{ij} = C_{ijk},{}^k - P^{kl}W_{kijl}
$$
is the Bach tensor.
The traces are given by: 
\begin{equation}\label{traces}
\begin{split}
 g^{ij}g_{ij}'' &= 2P_{ij}P^{ij}\\
 (n-4)g^{ij}g_{ij}''' &=  
-8P_{ij}B^{ij},  \qquad n\neq 4.
\end{split}
\end{equation}
The first equation of \eqref{traces} also holds for $n=4$ and the second
for $n=6$. 

The part of the derivatives depending linearly on curvature can be
calculated for all orders.  Differentiating the last equation of
\eqref{evol} shows that  
for $m\geq 2$, the trace $g^{ij}\pa_\rho^m g_{ij}|_{\rho = 0}$ has
vanishing linear part.  An easy
induction using the derivative of Ricci curvature  
\begin{equation}\label{derivricci}
R'_{ij} = \tfrac12(g'_{ik},_j{}^k + g'_{jk},_i{}^k -
g'_{ij},_k{}^k - g'_k{}^k,_{ij})
\end{equation}
and the Bianchi identity $P_{ik},{}^k = P_k{}^k,_i$
shows that at $\rho = 0$:
\begin{equation}\label{prinpart}
(4-n)(6-n)\cdots(2m-n)\pa_\rho^m g_{ij} = 2
\left (\Delta^{m-1}P_{ij}-\Delta^{m-2}P_k{}^k{},_{ij} \right ) 
+ \mbox{ lots }
\end{equation}
for $m\geq 2$ (and $m<n/2$ for $n$ even).  Here lots denotes quadratic and 
higher terms involving 
fewer derivatives of curvature and our sign convention is 
$\Delta = \nabla^k\nabla_k$.

A further observation can be made concerning the derivatives of 
$g_{ij}$ at $\rho = 0$:  each of them can be expressed in terms only of
Ricci 
curvature and its covariant derivatives.  Terms involving Weyl curvature
and its derivatives need not appear.
\begin{proposition}\label{onlyricci}
Each derivative $\pa_\rho^m g_{ij}|_{\rho = 0}$ can be expressed as a
linear combination of contractions of Ricci
curvature and covariant derivatives of Ricci curvature for the initial
metric $g$.  This holds for all $m$ for which these expressions are
determined:  for $m\geq 1$ for $n$ odd and for $1\leq m<n/2$
for $n$ even, and also for $g^{ij}\pa_\rho^{n/2} g_{ij}|_{\rho = 0}$ for
$n$ even.   
\end{proposition}
\begin{proof}
The proof is by induction on $m$.  We already know that $g_{ij}'|_{\rho = 0} =
2P_{ij}$.  Consider the inductive determination of $\pa_\rho^m
g_{ij}|_{\rho = 0}$ for $m\geq 2$ by taking the
equation obtained by setting the first expression of \eqref{evol} equal to
0, applying $\pa_{\rho}^{m-1}$, and setting $\rho = 0$.  The first, fourth
and 
fifth terms involve $\pa_\rho^m g_{ij}|_{\rho = 0}$.  The second and third
terms give rise to combinations of contractions of previously determined
derivatives which are of the desired form by the induction hypothesis.  The
fifth term also generates contractions of this form in addition to the term
involving $\pa_\rho^m g_{ij}$.  The result therefore 
follows if $\pa_\rho^{m-1} R_{ij}|_{\rho = 0}$ can be written only in terms 
of Ricci curvature and its covariant derivatives of the original
metric.  The  
first derivative is given by \eqref{derivricci}.  In differentiating
\eqref{derivricci} again, the derivative can fall on $g'$, on $g^{-1}$, or
on the 
connection.  Differentiating a Christoffel symbol shows that for a 
1-form $\eta_i$ which depends on $\rho$, one has 
$(\eta_i,{}_j)' = \eta'{}_i,{}_j-\Gamma'^k_{ij} \eta_k$, 
where $\Gamma'^k_{ij}$ is the
tensor $\Gamma'^k_{ij} = \frac12 g^{kl}(g'_{il},{}_j +g'_{jl},{}_i 
-g'_{ij},{}_l)$.  There is an analogous formula for the covariant
derivative of a tensor of higher rank.  Iterating such formulae together
with the Leibnitz formula and the formula for $g^{ij}{}'$, then setting $\rho
= 0$ and applying the induction hypothesis, it is clear 
that $\pa_\rho^{m-1} R_{ij}|_{\rho = 0}$ has the desired form.

If $n$ is odd, the trace $g^{ij}\pa_{\rho}^n g_{ij}|_{\rho = 0}$ is
determined from the third line of \eqref{evol} rather than from the first.  
But it is clear that differentiating the third line and using the induction
hypothesis gives rise to an expression of the desired form.
\end{proof}

We remark that it is a consequence of Proposition~\ref{onlyricci} that 
objects constructed solely out of the tensors 
$\pa_\rho^m g_{ij}|_{\rho = 0}$ also can be written in terms of Ricci
curvature and its covariant derivatives.  Two examples are the 
``conformally invariant powers of the Laplacian'' of \cite{GJMS} and
Branson's $Q$-curvature (\cite{Br}).  It is easily seen from the
construction in \cite{GJMS} that the coefficients of the conformally
invariant natural operators constructed there can be written in terms 
of the $\pa_\rho^m g_{ij}|_{\rho = 0}$; hence by
Proposition~\ref{onlyricci} in terms of the Ricci curvature and its
derivatives.  For $Q$-curvature, the same GJMS construction can be used to
establish the result following Branson's original definition.
Alternatively and more directly, this follows from the characterization of
$Q$-curvature given in \cite{FH}.  

The equations defining an ambient metric in normal form possess a symmetry
under reflection in $\rho$.  Let $R:\R_+\times M\times \R\rightarrow
\R_+\times M\times \R$ be $R(t,x,\rho)=(t,x,-\rho)$.  
If $\gt$ is an ambient metric for $(M,[g])$ in
normal form relative to a representative $g$, then $R^*\gt$ is also an
ambient metric for $(M,[g])$ but is not in normal form because condition
(3) of Definition~\ref{normalformdef} does not hold.  However, the
following Proposition is easily verified.
\begin{proposition}\label{reflection}
If $\gt$ is an ambient metric for $(M,[g])$ in normal form relative to $g$,
then $-R^*\gt$ is an ambient metric for $(M,[-g])$ in normal form relative
to $-g$.  
\end{proposition}

\noindent
Recall that these ambient metrics are unique up to the order specified in
Theorem~\ref{formal}(B).  

Next we prove a sharpened version of Theorem~\ref{formal} for $n=2$.  
We denote by $\tf$ the trace-free part with respect to $g$.  
\begin{theorem}\label{2dim}
Let $M$ be a smooth manifold of dimension 2.  Let $g$ be a smooth metric  
on $M$ and $h$ be smooth symmetric 2-tensor on $M$ satisfying
$g^{ij}h_{ij}=0$. Then there is an ambient metric $\gt$   
for $(M,[g])$ in normal form relative to $g$ which satisfies 
$\tf\left(\pa_\rho \gt_{ij}|_{\rho = 0}\right) = t^2 h_{ij}$.  
These conditions 
uniquely determine $\gt$ to infinite order at $\rho = 0$.   
The metric $\gt$ 
is straight to infinite order if and only if $h_{ij},{}^j = \frac12 R,_i$.      

\end{theorem}
\begin{proof}
The proof begins the same way as the proof of Theorem~\ref{formal} for
$n>2$.  The components $\gt_{I\nf}$ 
are determined by Lemma~\ref{explicitnormal}.  The computations leading to 
\eqref{firstorder} remain valid.  The vanishing of the first two lines
of \eqref{firstorder} is again equivalent to $\pa_\rho a|_{\rho = 0} =2$  
and $\pa_\rho b_i|_{\rho = 0} = 0$.  However 
when $n=2$, the coefficient $(\pa_\rho a -n)$ vanishes in the third line of 
\eqref{firstorder}.  For $n=2$, one always has $2R_{ij}= Rg_{ij}$, 
so vanishing of the third line is therefore equivalent to  
\begin{equation}\label{tracen2}
g^{kl}\pa_\rho g_{kl}|_{\rho = 0}=R.
\end{equation}
Thus the trace of 
$\pa_\rho g_{ij}|_{\rho = 0}$ is determined, but the trace-free part 
remains undetermined by the Einstein condition.  We observe that this
information is   
already enough to prove part (B) of Theorem~\ref{formal} when $n=2$:   
we have now shown that a solution $\gt$ is uniquely determined modulo 
$O^+_{IJ}(\rho)$.  
The prescription
$\tf(\pa_\rho g_{ij})|_{\rho = 0} = h_{ij}$ fixes the ambiguity 
in $\pa_\rho g_{ij}|_{\rho = 0}$.  It is convenient to define $P_{ij}$ by 
$$
2P_{ij} = h_{ij} +\tfrac12 R g_{ij}
$$
so that $g_{ij}(x,\rho)$ is still given by the second line of
\eqref{initial}.  

We will consider the inductive determination of the Taylor expansion of
$\gt$ for higher orders as in the proof of Theorem~\ref{formal} for $n>2$.
The next order is 
the tricky one.  For definiteness, we define $\gt^{(1)}$ to be given by 
\eqref{spform} with $g_{ij}(x,\rho)$ given by the second line of
\eqref{initial} (say with the $O(\rho^2)$ term set to $0$); i.e., we fix
the $O(\rho^2)$ indeterminacy in $a$ and $b_i$  
in \eqref{initial} to be $0$.  With notation as in the proof of
Theorem~\ref{formal} above, Lemma~\ref{0comp} implies that
$\Rt^{(1)}_{0I}=0$.  Also, $\Rt^{(1)}_{i\nf}$ is given by the middle line
of \eqref{evol}.  In particular, 
\begin{equation}\label{riccomp}
\Rt^{(1)}_{i\nf}|_{\rho = 0} = P_{ij},{}^j- P_j{}^j,_i 
= \tfrac12 \left( h_{ij},{}^j - \tfrac12 R,_i\right).
\end{equation} 
 
Set $\gt^{(2)}_{IJ} = \gt^{(1)}_{IJ}+ \Phi_{IJ}$
with $\Phi_{IJ}$ given by \eqref{phi} with $m=2$.  
The perturbed Ricci tensor is given by \eqref{perturbricci}.  
The first line says $\Rt^{(2)}_{00} = O(\rho^2)$.  The fourth line shows
that $\Rt^{(2)}_{0\nf} = O(\rho)$ if and only if $\phi_{00}=O(\rho)$, so we
require $\phi_{00}|_{\rho = 0}= 0$.  The second line then gives  
$\Rt^{(2)}_{0i} = O(\rho^2)$.  The fifth line shows that 
$\phi_{0i}|_{\rho  = 0}$ can be uniquely chosen so that 
$\Rt^{(2)}_{i\nf} = O(\rho)$, and \eqref{riccomp} shows that the so
determined $\phi_{0i}|_{\rho  = 0}$ vanishes if and only if 
$h_{ij},{}^j = \frac12 R,_i$.   The third line shows that the trace-free
part of $\phi_{ij}|_{\rho = 0}$ can be uniquely chosen so that  
\begin{equation}\label{modtrace}
\Rt^{(2)}_{ij} = \la\rho g_{ij} + O(\rho^2)
\end{equation}
for some function $\la$ on $M$.  
The last line shows that the trace of $\phi_{ij}|_{\rho = 0}$ can be
uniquely determined so that $\Rt^{(2)}_{\nf\nf} = O(\rho)$.  Thus all 
$\phi_{IJ}|_{\rho = 0}$ have been determined, and all components of 
$\Rt^{(2)}_{IJ}$ vanish to the desired orders except for \eqref{modtrace}. 
Consider now the last line of \eqref{bcomp} applied to $\Rt^{(2)}$ with
$m=2$.  It reduces to $g^{ij}\pa_\rho \Rt^{(2)}_{ij} = O(\rho)$.  Thus
$\la=0$ as desired.  Therefore we have shown that one can uniquely determine 
$\Phi_{IJ} \mod O(\rho^3)$ to make $\Rt^{(1)}_{IJ} = O(\rho^2)$ for 
$I$, $J\neq \nf$ and $\Rt^{(1)}_{I\nf} = O(\rho)$.  Moreover, the metric 
$\gt^{(2)}_{IJ}$ so determined is of the form \eqref{spform} $\mod
O(\rho^3)$ if and only if $h_{ij},{}^j = \frac12 R,_i$.  

Consider now the induction for higher $m$.  The argument proceeds as
in the proof of Theorem~\ref{formal} above.  The relevant
coefficients in the first three lines of \eqref{perturbricci} never vanish
for 
$n=2$ and $m\geq 3$.  Thus the conditions $\Rt^{(m)}_{IJ} = O(\rho^m)$ for
$I$, $J \neq \nf$ uniquely determine $\Phi_{IJ} \mod O(\rho^{m+1})$.  Then 
the three lines of \eqref{bcomp} successively show that 
$\Rt^{(m)}_{I\nf} = O(\rho^{m-1})$.  Thus the induction continues to all
orders.  
 
If $\gt$ is straight to infinite order, then $\gt_{0i} = 0$ to infinite 
order, so by the determination of $\phi_{0i}|_{\rho = 0}$ when $m=2$ noted
above, we must have  
$h_{ij},{}^j = \frac12 R,_i$.  Conversely, if $h_{ij},{}^j = 
\frac12 R,_i$, then we have $\phi_{0i}|_{\rho = 0}$ when $m=2$.  The
argument of Proposition~\ref{simpler} then shows that  
$\gt_{00}=2\rho$ and $\gt_{0i}=0$ to all higher orders.  
\end{proof}
Theorem~\ref{formal} for $n=2$ is a consequence of Theorem~\ref{2dim} and
its proof.  Part (A) follows upon choosing any $h_{ij}$ satisfying
$h_i{}^i = 0$.  We already noted that part (B) holds at the end of the 
first paragraph in the proof of Theorem~\ref{2dim} above. 

We remark that in the straight case, i.e. when $h_{ij},{}^j = \frac12
R,_i$, the solution $\gt$ in Theorem~\ref{2dim} can be written explicitly;
see Chapter~\ref{flat}.  

Observe in Theorem~\ref{2dim} that the metric $\gt$ may be put into normal
form relative to another metric $\wh{g}$ in the conformal class, giving
rise to another trace-free tensor $\wh{h}$.  Since the
straightness condition is invariant under diffeomorphisms, $h$ satisfies 
$h_{ij},{}^j=\frac12 R,_i$ if and only if $\wh{h}$ satisfies
$\wh{h}_{ij},{}^j=\frac12 \wh{R},_i$, where in the latter equation the
covariant derivative and scalar curvature are that of $\wh{g}$.  Note also
that if $g$ has constant scalar curvature,  
then $h_{ij}=0$ satisfies $h_i{}^i=0$ and $h_{ij},{}^j = \frac12
R,_i$.  
In the case of definite signature, the uniformization theorem
implies that every conformal class $(M,[g])$ on any 2-manifold contains a 
metric of constant 
scalar curvature.  Thus it follows that for any definite signature metric
$g$, there exists a trace-free $h$ satisfying $h_{ij},{}^j=\frac12 R,_i$.  

When $n\geq 4$ is even, the existence of formal power series solutions for
the ambient metric at order $n/2$ is 
in general obstructed.  The obstruction can be identified as a conformally
invariant natural tensor generalizing the Bach tensor in dimension 4, which
we call the ambient obstruction tensor and denote $\cO_{ij}$.  We next
define the obstruction tensor and establish its
basic properties.  Suppose that $n\geq 4$ is even and that $\gt$ is 
an ambient metric for $(M,[g])$.  By Theorem~\ref{main}, $\gt$ is 
uniquely determined modulo 
$O_{IJ}^+(\rho^{n/2})$ up to a homogeneous diffeomorphism of $\cGt$ which 
restricts to the identity on $\cG$.  Set $Q =  \| T \|^2 = \gt(T,T)$, where
as usual $T$ denotes the infinitesimal dilation.  
Then $Q$ is a defining function for $\cG\times \{0\}\subset \cGt$  
invariantly associated to $\gt$, which is homogeneous of degree 2. 
(To see that $Q$ is a defining function, one can put $\gt$ into normal
form, whereupon Proposition~\ref{simpler} shows that $Q= 2\rho t^2 \mod
O(\rho^{n/2 +1})$.)   We  
identify $\cG$ with $\cG\times \{0\}$ via the inclusion $\iota$.
Now $\Ric(\gt)=O^+_{IJ}(\rho^{n/2-1})$, so 
$(Q^{1-n/2}\Ric{\gt})|_{T\cG}$ is a tensor field on $\cG$, 
homogeneous of degree $2-n$, which annihilates $T$.  It therefore defines a 
symmetric 2-tensor-density on $M$ of weight $2-n$, which is trace-free.
If $g$ is a metric in the conformal class, evaluating this tensor-density
at the image of $g$ viewed as a section of $\cG$ defines a 2-tensor on $M$
which we denote by 
$(Q^{1-n/2} \Ric{\gt})|_g$.  We define the obstruction tensor of $g$ to
be  
\begin{equation}\label{obsdef}
\cO = c_n (Q^{1-n/2} \Ric{\gt})|_g, 
\qquad c_n = (-1)^{n/2-1} \frac{2^{n-2} (n/2-1)!^2}{n-2}.
\end{equation}
For $\gt$ in normal form relative to $g$, this reduces to 
$$
\cO_{ij} = 2^{1-n/2}c_n (\rho^{1-n/2}\Rt_{ij})|_{\rho=0}.
$$

\begin{theorem}\label{obstruction}
Let $n\geq 4$ be even.  The obstruction tensor $\cO_{ij}$ of $g$ is
independent of 
the choice of ambient metric $\gt$ and has the following properties:   
\begin{enumerate}
\item
$\cO$ is a natural tensor invariant of the metric $g$;
i.e. in local coordinates the components of $\cO$ are given by
universal polynomials in the components of $g$, $g^{-1}$ and the curvature 
tensor of $g$ and its covariant derivatives, and can be written just in
terms of the Ricci curvature and its covariant derivatives.
The expression for $\cO_{ij}$  takes the form 
\begin{equation}\label{Oform}
\begin{split}
\cO_{ij} &= \Delta^{n/2-2}\left (P_{ij},_k{}^k-P_k{}^k,_{ij}\right )
+ lots\\
&=(3-n)^{-1}\Delta^{n/2-2}W_{kijl},{}^{kl} + lots,
\end{split}
\end{equation}
where  
$\Delta = \nabla^i\nabla_i$
and lots denotes quadratic and higher terms in curvature involving fewer 
derivatives.  
\item
One has
$$
\cO_i{}^i=0 \qquad\qquad\qquad \cO_{ij},{}^j = 0.
$$
\item
$\cO_{ij}$ is conformally invariant of weight $2-n$; i.e. if $0<\Omega \in 
  C^{\infty}(M)$ and $\wh{g}_{ij} =
\Omega^2 g_{ij}$, then $\wh{\cO}_{ij} = \Omega^{2-n}\cO_{ij}$.   
\item
If $g_{ij}$ is conformal to an Einstein metric, then $\cO_{ij}=0$.  
\end{enumerate}
\end{theorem}
\begin{proof}
We can assume that $\gt$ is in normal form relative to $g$.
Then 
$\gt$ is unique up to addition of $\Phi_{IJ}$ of the form \eqref{phi} with 
$m=n/2$, where $\phi_{00}$, $\phi_{0i}$, and $g^{ij}\phi_{ij}$ all vanish
at $\rho = 0$.  The independence of $(Q^{1-n/2} \Ric{\gt})|_g$ on $\gt$ 
is then an immediate consequence of \eqref{perturbricci}.  

According to 
Proposition~\ref{simpler}, we may take $\gt$ to satisfy $\gt_{00}=2\rho$,  
$\gt_{0i}=0$.  Then $\cO_{ij}$ may be obtained by setting
$\Rt_{ij} = c_n^{-1}(2\rho)^{n/2-1} \cO_{ij} \mod O(\rho^{n/2})$ in
\eqref{evol}, applying $\pa_{\rho}^{n/2 -1}|_{\rho = 0}$, and taking the 
trace-free part.
This shows that $\cO_{ij}$ is a natural tensor and
Proposition~\ref{onlyricci} shows that it can be written just in terms of
Ricci and its derivatives.  Its linear part may be
calculated using \eqref{derivricci} and \eqref{prinpart} to be given by  
the first line of \eqref{Oform}.  The second line follows from the fact
that 
$
W_{kijl},{}^{kl} = (3-n)(P_{ij},_k{}^k - P_{ik},_j{}^k).
$

We have already observed that $\cO_{ij}$ is trace-free.  Its conformal
invariance follows from its definition in terms of a
conformally invariant tensor-density.  If $g_{ij}$ is Einstein, we have
previously noted that there is a solution for $\gt$ to all orders, so
$\cO_{ij}=0$.  

It only remains to establish that $\cO_{ij},{}^j = 0$.  This follows from
the Bianchi identity as follows.  Recall that the hypotheses for
\eqref{bcomp} were that $\Rt_{IJ}=O(\rho^{m-1})$ for 
$I,\,J \neq \nf$ and $\Rt_{I\nf}=O(\rho^{m-2})$.  Our metric $\gt$
satisfies these with $m=n/2+1$ except that 
$c_n\Rt_{ij}= (2\rho)^{n/2-1} \cO_{ij} \mod O(\rho^{n/2})$.
If one recalculates the middle line of \eqref{bcomp} allowing the
possibility that $\Rt_{ij} = O(\rho^{m-2})$ but all other components
vanish as before, one finds that there are two extra terms:
$g^{jk}\nabla_j\Rt_{ik} -\frac12 g^{jk}\nabla_i \Rt_{jk}$.  For $m=n/2 +1$,  
the coefficient of $\Rt_{i\nf}$ vanishes, so we obtain
$$
t(\pa_\rho \Rt_{i0}- \pa_i\Rt_{0\nf}) + c'\rho^{n/2-1}\cO_{ij},{}^j  
=O(\rho^{n/2})
$$
for some nonzero constant $c'$.  We may as well take $\gt$ to have 
$\gt_{00} = 2\rho$, $\gt_{i0} = 0$, in which case the first two terms
vanish, giving the desired conclusion.
\end{proof}

For $n=4,\,6$ one can calculate $\cO_{ij}$ by hand by carrying out the
computation indicated in the proof of Theorem~\ref{obstruction}.  One
obtains the tensor which obstructs the validity of \eqref{deriv}:
$\cO_{ij} = B_{ij}$ for $n=4$ and 
\[
\begin{split}
\cO_{ij} = B_{ij},_k{}^k - 2W_{kijl}&B^{kl} -4P_k{}^kB_{ij} +8P^{kl}C_{(ij)k},_l 
-4C^k{}_i{}^lC_{ljk} \\
&+2C_i{}^{kl}C_{jkl} +4P^k{}_k,_lC_{(ij)}{}^l
-4W_{kijl}P^k{}_mP^{ml}
\end{split}
\]
for $n=6$.  

When the obstruction tensor is nonzero, there are no formal power series
solutions for $\gt$ beyond order $n/2$.  However, one can still continue
the solution to higher orders by introducing log terms.  In this
case one 
is obliged to introduce an indeterminacy and the solution is no longer
determined solely by the initial metric.  We have already seen this
indeterminacy phenomenon when $n=2$ in Theorem~\ref{2dim}, although for
$n=2$ there 
is no obstruction and consequently there are no log terms.   When $n$ is  
odd, there are solutions with expansions involving half-integral powers
which also have an indeterminacy at order $n/2$.  

We broaden our terminology to encompass such metrics.  Recall that in  
Definition~\ref{ambientdef}, we required an ambient metric to be smooth.  
We now define a {\it generalized ambient metric}
to be a metric $\gt$ satisfying all the conditions of
Definition~\ref{ambientdef}, except that the smoothness condition is
relaxed to the requirement that $\gt$ be 
$C^{\infty}(\cGt \setminus \{\rho =0\}) \cap C^1(\cGt)$, and in all
dimensions we require $\Ric({\gt})$ to vanish to infinite order along 
$\cG\times \{0\}$ (i.e. all derivatives of all components of $\Ric({\gt})$ 
extend continuously to $\rho=0$ and vanish there).
Definition~\ref{normalformdef} and 
Lemma~\ref{explicitnormal} concerning metrics in normal form extend to
generalized ambient metrics. 
Proposition~\ref{straightequiv} and the notion of straightness also extend 
to generalized ambient metrics.

We now discuss the existence and uniqueness of generalized ambient
metrics in normal form, beginning with the case $n$ odd.

\begin{theorem}\label{odd}
Let $M$ be a smooth manifold of odd dimension $n$.  
Suppose given a smooth metric $g$ 
and a smooth symmetric 2-tensor $h$ on $M$ satisfying   
$g^{ij}h_{ij} = 0$.  Then there exists a generalized ambient metric 
$\gt$ which is in normal form relative to $g$ and whose restriction to 
either $\cGt\cap \{\rho > 0\}$ or $\cGt\cap \{\rho < 0\}$ 
has the form 
$$
\gt_{IJ} = \psi^{(0)}_{IJ} + \psi^{(1)}_{IJ}|\rho|^{n/2} 
$$
where $\psi^{(0)}_{IJ}$, $\psi^{(1)}_{IJ}$ extend smoothly up to $\rho =0$ 
and $\tf\left(\psi^{(1)}_{ij}|_{\rho =0}\right) = t^2 h_{ij}$.  The Taylor
expansions of 
the $\psi^{(0)}_{IJ}$ and $\psi^{(1)}_{IJ}$ are uniquely determined 
to infinite order by these conditions, and the solution satisfies 
$g^{ij}\psi^{(1)}_{ij}|_{\rho =0}=0$.  The metric $\gt$ is
straight to 
infinite order if and only if $h_{ij},{}^j = 0$.    
\end{theorem}
\begin{proof}
We construct $\gt$ separately on $\cGt\cap \{\rho > 0\}$ and on $\cGt\cap
\{\rho < 0\}$.  
First consider $\{\rho >0\}$.    

Return to the inductive construction of $\gt$ in Theorem~\ref{formal}.  
If we pause in that construction at $m=(n-1)/2$, we have $\gt^{((n-1)/2)}$ 
determined uniquely mod $O(\rho^{(n+1)/2})$ by \eqref{geocond} and the 
condition $\Ric_{IJ}(\gt^{((n-1)/2)}) = O(\rho^{(n-1)/2})$ for $I$,~$J\neq 
\nf$.  
As in Theorem~\ref{simpler}, we may as well
choose $\gt^{((n-1)/2)}$ to be of the form \eqref{spform}, in which case we
have $\Rt_{00}^{((n-1)/2)}=0$, $\Rt_{0i}^{((n-1)/2)}=0$, and of course 
$\Rt_{ij}^{((n-1)/2)}$ is smooth, i.e. it has no half-integral powers in
its expansion.   
Previously we considered only formal power series solutions, so we
next modified $\gt$ at order $(n+1)/2$.  Now we 
instead modify $\gt$ at order $n/2$:  set
$\gt^{(n/2)}_{IJ} = \gt^{((n-1)/2)}_{IJ} + \Phi_{IJ}$, where $\Phi_{IJ}$ is
of the form 
\eqref{phi} with $m=n/2$ and each $\phi_{IJ}$ is asymptotic to
a formal power series in $\sqrt \rho$.  It follows from \eqref{Ricci}
that each component of the Ricci curvature of such a metric is also
asymptotic to a formal power series in $\sqrt \rho$.  (For $n=3$, the
series for the components $\Rt_{I\nf}$ may in principle also contain a
$\rho^{-1/2}$ term arising from the $\rho^{3/2}$ term in $\gt$; see below.)
Now \eqref{perturbricci} still 
holds except the $\Rt^{(m-1)}_{IJ}$ which appear on the right
hand side are now replaced by $\Rt^{(m-1/2)}_{IJ}$ and the error
terms are shifted by $\frac12$:  the error terms
in the first three lines are $O(\rho^{m-1/2})$ and those in the last three
lines are $O(\rho^{m-3/2})$.  Vanishing of the first three lines of
\eqref{perturbricci} gives at $\rho = 0$:  $\phi_{00} = 0$, 
$\phi_{0i} = 0$,  
$g^{ij}\phi_{ij} = 0$, but the trace-free part of $\phi_{ij}$ may be chosen 
arbitrarily.  We define $\gt^{(n/2)}$ by taking $\phi_{ij}=h_{ij}$ and 
$\gt^{(n/2)}$ of the form \eqref{spform}, so that
$\Rt_{00}^{(n/2)}=0$, $\Rt_{0i}^{(n/2)}=0$.  Also, 
$\Rt_{ij}^{(n/2)}$ is 
asymptotic to a formal power series in $\sqrt \rho$, which by 
construction satisfies $\Rt_{ij}^{(n/2)} = O(\rho^{(n-1)/2})$.  
We now 
modify $g^{(n/2)}$ by addition of a term \eqref{phi} with $m=(n+1)/2$ to
obtain $\gt^{((n+1)/2)}$.  
Once again \eqref{perturbricci} holds with the shifted error terms and
superscripts on $\Rt_{IJ}$.
None of the relevant constants which appear in \eqref{perturbricci} vanish
for this value of $m$, so we deduce that $\phi_{00}$, $\phi_{0i}$, and 
$\phi_{ij}$ are all detemined at $\rho = 0$ and once again $\phi_{00}$ and 
$\phi_{0i}$ may be taken to be identically 0.  Now
$\Rt_{00}^{((n+1)/2)}=0$, $\Rt_{0i}^{((n+1)/2)}=0$, and
$\Rt_{ij}^{((n+1)/2)}$ is 
asymptotic to a formal power series in $\sqrt \rho$ satisfying 
$\Rt_{ij}^{((n+1)/2)} = O(\rho^{n/2})$. 

Before proceeding to the next value of $m$, consider the components 
$\Rt_{I\nf}^{((n+1)/2)}$.  We have $\Rt_{0\nf}^{((n+1)/2)} = 0$ by
Lemma~\ref{0comp}.  The 
components $\Rt^{((n+1)/2)}_{i\nf}$ and $\Rt^{((n+1)/2)}_{\nf\nf}$ 
are given by formal power series in $\sqrt \rho$.  (When $n=3$, at first
glance it appears
from the third line of \eqref{evol} that the $\rho^{3/2}$ term in $g_{ij}$ 
generates a $\rho^{-1/2}$ term in $\Rt^{((n+1)/2)}_{\nf\nf}$.  However 
this term has coefficient 0 because $g^{ij}h_{ij}=0$.)  
Now it is easily
checked that \eqref{bcomp} holds for $m \in \frac12 \mathbb Z$ 
and for the $\gt$ that we are considering with expansions in 
$\sqrt \rho$, still 
under the same hypotheses:  $\Rt_{IJ}=O(\rho^{m-1})$ for 
$I,\,J \neq \nf$ and $\Rt_{I\nf}=O(\rho^{m-2})$.  (When $n=3$ one
modification is required:  the error term in the last line is 
$O(\rho^{m-1})+O(\rho^{1/2} \Rt_{0\nf})$, owing to the $\rho^{3/2}$
term in the expansion of $g_{ij}$.)  
If we proceed by induction on the order of vanishing of 
$\Rt^{((n+1)/2)}_{i\nf}$ and $\Rt^{((n+1)/2)}_{\nf\nf}$ in 
\eqref{bcomp} similarly to the proof of Theorem~\ref{formal}, 
now using $\Rt_{0I}^{((n+1)/2)}=0$ and $\Rt_{ij}^{((n+1)/2)} =
O(\rho^{n/2})$, we 
find  $\Rt_{i\nf}^{((n+1)/2)}$, $\Rt_{\nf\nf}^{((n+1)/2)} =
O(\rho^{n/2-1})$.  We make one further observation about 
$\Rt_{i\nf}^{((n+1)/2)}$; namely, 
$(\rho^{1-n/2}\Rt_{i\nf}^{((n+1)/2)})|_{\rho = 0}$ is a constant multiple
of $h_{ij},{}^j$.  To see this, note that $\gt^{((n+1)/2)}$ is of the form 
\eqref{spform} where $g_{ij} = \eta_{ij}^{(0)} + \eta_{ij}^{(1)}\rho^{n/2}$ 
for $\eta_{ij}^{(0)}$, $\eta_{ij}^{(1)}$ smooth and 
$\eta_{ij}^{(1)}|_{\rho = 0} = h_{ij}$.  The middle line of 
\eqref{evol} together with the fact that $g^{ij}h_{ij}=0$ then imply that 
the $\rho^{n/2-1}$ coefficient in the expansion of
$\Rt_{i\nf}^{((n+1)/2)}$ is a multiple of $h_{ij},{}^j$.

Return now to the inductive construction of $\gt$.  We next define
$\gt^{(n/2+1)}_{IJ}=\gt^{((n+1)/2)}_{IJ} +\Phi_{IJ}$ with $\Phi_{IJ}$ given
by \eqref{phi} with $m=n/2+1$.  The first line of \eqref{perturbricci}
(with shifted error) tells us that $\Rt^{(n/2+1)}_{00}=O(\rho^{(n+1)/2})$ 
independent of the choice of $\phi_{00}$.  However, the fourth line of
\eqref{perturbricci} tells us that we must choose $\phi_{00} = 0$ at 
$\rho = 0$ in order to make $\Rt^{(n/2+1)}_{0\nf} = O(\rho^{(n-1)/2})$.  
Now the second line of \eqref{perturbricci} says that 
$\Rt^{(n/2+1)}_{0i}=O(\rho^{(n+1)/2})$ independent of the choice of
$\phi_{0i}$.  However, the fifth line determines $\phi_{0i}|_{\rho = 0}$ by
the requirement that $\Rt^{(n/2+1)}_{i\nf}=O(\rho^{(n-1)/2})$.  The third
line then determines $\phi_{ij}|_{\rho = 0}$ by the requirement that 
$\Rt_{ij}=O(\rho^{(n+1)/2})$.    Taking $m=(n+1)/2$ in \eqref{bcomp}, we
already know that the 
first two lines hold, and the third line tells us that 
$\Rt^{(n/2+1)}_{\nf\nf}=O(\rho^{(n-1)/2})$.
Thus the $\phi_{IJ}|_{\rho =0}$ have been uniquely determined and we have  
$\Rt^{(n/2+1)}_{IJ}=O(\rho^{(n+1)/2})$ for $I$, $J\neq \nf$ and 
$\Rt^{(n/2+1)}_{I\nf}=O(\rho^{(n-1)/2})$.  In this determination, we found 
$\phi_{00}|_{\rho = 0} = 0$ and $\phi_{0i}|_{\rho = 0}$ was determined by
the fifth line of \eqref{perturbricci}.  By the observation noted above
that 
$(\rho^{1-n/2}\Rt_{i\nf}^{((n+1)/2)})|_{\rho = 0}$ is a constant multiple
of $h_{ij},{}^j$, we deduce that $\phi_{0i}|_{\rho = 0}$ is a constant
multiple of $h_{ij},{}^j$, and in particular $\phi_{0i}|_{\rho = 0}$
vanishes if and only if $h_{ij},{}^j = 0$.  

Now we modify $\gt$ to higher orders successively by induction increasing
$m$ by $1/2$ each step.  The induction statement is that $\gt^{(m)}_{IJ}$
is uniquely determined mod $O(\rho^{m+1/2})$ by the requirement that  
$\Rt^{(m)}_{IJ} = O(\rho^{m-1/2})$ for $I$, $J\neq \nf$.  We have
established above the case $m=n/2+1$.  Up through $m=n-1/2$, the induction
step follows from the first 
three lines of \eqref{perturbricci} just as in the proof of
Theorem~\ref{formal}.  Again just as in the proof of Theorem~\ref{formal},
we then deduce that $\Rt^{n-1/2}_{I\nf} = O(\rho^{n-2})$ using
\eqref{bcomp}.  

For $m=n$ we encounter the vanishing of the coefficient of the trace in 
the third line of \eqref{perturbricci}.  The same reasoning used in the
proof of Theorem~\ref{formal} applies here:  the fact that we are
increasing $m$ by $1/2$ rather than 1 at each step plays no role in that
analysis.  The induction for higher $m$ then proceeds as usual, including 
the induction based on \eqref{bcomp} for the $\Rt_{I\nf}$ components.  
Thus it follows that there is a unique series for $\gt$ of the desired form 
for which all components of $\Ric(\gt)$ vanish to infinite order.

Finally we observe that, just as in the proof of Proposition~\ref{simpler}, 
once we get beyond $m=n/2+1$ the special form \eqref{spform} is preserved
in the induction.  So if $h_{ij},{}^j=0$, then the solution has
$\gt_{00}=2\rho$, $\gt_{0i}=0$ to infinite order.

We construct $\gt$ for $\{\rho <0\}$  
by using the reflection $R$ as in Proposition~\ref{reflection}.  
If $\gt$ is a solution on $\{\rho >0\}$ with data $-g$, $-h$, then
$-R^*\gt$  
is a solution on $\{\rho <0\}$ with data $g$, $h$ and vice versa.
That the solutions match to first order at $\rho =0$ 
can be checked using \eqref{initial}.  In fact, 
Proposition~\ref{reflection} and the uniqueness of the expansion imply 
that the solution obtained by reflection is $C^{(n-1)/2}$ across $\rho =0$.     
\end{proof}

The analogue of Theorem~\ref{odd} for $n$ even is the following.
\begin{theorem}\label{even}
Let $M$ be a smooth manifold of even dimension $n\geq 4$. 
If $g$ is a smooth metric and
$h$ a smooth symmetric 2-tensor on $M$ satisfying   
$g^{ij}h_{ij} = 0$, then there exists a generalized ambient metric $\gt$
in normal form relative to $g$, which has an expansion of the form 
\begin{equation}\label{logexpand}
\gt_{IJ} \sim \sum_{N= 0}^{\infty} \gt^{(N)}_{IJ}(\rho^{n/2}\log|\rho|)^N    
\end{equation}
where the $\gt^{(N)}_{IJ}$ are smooth on $\cGt$ and 
$\tf\left(\partial_\rho^{n/2} \gt^{(0)}_{ij}\right) = 
t^2 h_{ij}$ at    
$\rho =0$.   
The Taylor expansions of the $\gt^{(N)}_{IJ}$ are uniquely determined to 
infinite order by these conditions.  
The solution $\gt$ is smooth 
(i.e. $\gt^{(N)}$ vanishes to infinite order for $N\geq 1$) if and only if  
the obstruction tensor $\cO_{ij}$ vanishes on $M$. There is   
a natural pseudo-Riemannian invariant 1-form $D_i$ such that the solution    
$\gt$ is straight to infinite order if
and only if $h_{ij},{}^j = D_i$.     
\end{theorem}
\begin{remark}
Observe that the corresponding result when $n=2$ (Theorem~\ref{2dim}) 
takes precisely the same form, except that there are no log terms and the
obstruction tensor always vanishes.  Theorem~\ref{2dim} shows that when
$n=2$, one has $D_i = \frac12 R,_i$.   
For general even $n$, $D_i$ is given by \eqref{Dformula} below, and is the
same tensor which appears in Theorem 3.5 of \cite{GrH2}.  When $n=4$, 
one has 
$$
D_i=4P^{jk}P_{ij,k}-3 P^{jk}P_{jk,i}+ 2 P_i{}^jP^k{}_{k,j}.
$$
\end{remark}
\begin{remark}
We do not know whether for any $(M,g)$ the equation $h_{ij},{}^j = D_i$  
admits a solution.  We have already observed after the proof of
Theorem~\ref{2dim} that the existence of a solution is a conformally
invariant condition and that there is always a solution if $g$ is definite
and $n=2$.  For $n\geq 3$, the equation $h_{ij},{}^j = f_i$ always admits a
local solution for any smooth 1-form $f_i$.  This  
follows from the fact that at each point of the cotangent bundle minus the
zero section, the symbol of the operator
$\operatorname{div}:\odot^2_0T^*M\rightarrow T^*M$ is surjective, so that
there 
is a right parametrix (see Theorem 19.5.2 of \cite{Ho}).  This also 
implies that on a compact manifold, the range of $\operatorname{div}$ 
has finite codimension.  For $g$ definite and $M$ compact, the 
range of $\operatorname{div}$ is  
the $L^2$ orthogonal complement of the space of global conformal Killing
fields. Using this and the explicit formula for $D_i$ above, one can show
that there exists a solution to $h_{ij},{}^j = D_i$ if 
$n=4$, $g$ is definite, and $M$ is compact.
\end{remark}

We prepare to prove Theorem~\ref{even}.  
Denote by $\cA$ the space of formal asymptotic expansions of scalar
functions $f$ on $\R_+\times M \times \R$ of the form  
$$
f\sim \sum_{N\geq 0} f^{(N)}(\rho^{n/2}\log |\rho|)^N
$$
where each $f^{(N)}$ is smooth and homogeneous of degree 0 in $t$.
It is easily checked that $\cA$
is an algebra, that $\cA$ is preserved by $\pa_{x^i}$ and $\rho
\pa_{\rho}$, and that $f^{-1}\in \cA$ if $f \in \cA$ and $f \neq 0$ at 
$\rho =0$.  The metric $\gt$ which we construct actually has a more refined 
expansion than \eqref{logexpand}.  
Let $\cM$ denote the space of formal asymptotic expansions of metrics $\gt$
on $\R_+\times M \times \R$ of the form    
\begin{equation}\label{refined}
\gt_{IJ}=
\begin{pmatrix}
2\rho + \al &tb_j &t\\
tb_i & t^2 g_{ij}& 0\\
t & 0 & 0
\end{pmatrix}
\end{equation}
with $\al \in \rho^{n/2+2}\cA$, $b_i\in \rho^{n/2+1}\cA$, $g_{ij}\in \cA$.
We will show that there is a unique expansion $\gt\in \cM$ satisfying  
that $g_{ij}$ is the given representative at $\rho =0$, that 
$\tf(\pa_{\rho}^{n/2}g_{ij}^{(0)}) = h_{ij}$ at $\rho =0$,  
and that $\Ric(\gt)=0$.   

For $\gt_{IJ}$ of the form \eqref{refined}, the inverse $\gt^{IJ}$ and
Christoffel symbols $\Gat_{IJK}$ are given by \eqref{inverse} and \eqref{c}
with $a=2\rho +\al$.  Observe that for all $IJ$ we have
$t^{2-\#(IJ)}\gt^{IJ}\in \cA$, where 
$\#(IJ)$ denotes the number of zeros in the list $IJ$.  Also, for the
off-anti-diagonal elements there is an improvement:  
$t^{2-\#(IJ)}\gt^{IJ}\in \rho \cA$ unless both $IJ$ are between $1$ and 
$n$ or one is $0$ and the other $\infty$.  Similarly, for the Christoffel
symbols we have:   
$\rho\Gat_{IJK}\in t^{2-\#(IJK)}\cA$ for all components, and 
$\Gat_{IJK}\in t^{2-\#(IJK)}\cA$ unless two of $IJK$ are between $1$ and
$n$ and the third is $\infty$.

\begin{proposition}\label{formRt}
If $\gt\in \cM$, then $t^2\Rt_{00}\in \rho^{n/2+1}\cA$, 
$t\Rt_{0j},t\Rt_{0\infty}\in \rho^{n/2}\cA$, 
$\rho\Rt_{ij}$, $\rho\Rt_{i\infty}\in \cA$,
$\rho^2\Rt_{\infty\infty}\in \cA$.
\end{proposition}
\begin{proof}
We first derive a formula for the components $\Rt_{0J}$ 
for $\gt$ of the form \eqref{refined}.  
Set $E_{IJ} = T_{[I,J]}$.  The fact that $\cL_T \gt_{IJ} = 2\gt_{IJ}$  
implies that $T_{(I,J)} = \gt_{IJ}$, so 
$T_{I,J} = \gt_{IJ} + E_{IJ}$.  Thus $T_{I,JK} = E_{IJ,K}$.  
Now 
$t\Rt_{0IJK} = T^L\Rt_{LIJK} = 2T_{I,[JK]}= 2E_{I[J,K]}$,
so $t\Rt_{0J}=t\gt^{IK}\Rt_{0IJK}=\gt^{IK}E_{IJ,K}$.
Expand the covariant derivative to obtain 
\begin{equation}\label{Ric0}
t\Rt_{0J}=\gt^{IK}\pa_K E_{IJ}-\gt^{IK}\gt^{PQ}\Gat_{IKP}E_{QJ}
-\gt^{IK}\gt^{PQ}\Gat_{JKP}E_{IQ}.
\end{equation}

The components of $T_I$ are given explicitly by   
$T_I = t\gt_{I0} = (2\rho t +\al t,t^2b_i,t^2)$.  The 1-form
$2\rho t dt + t^2d\rho$ is $d(\rho t^2)$, so these terms can be ignored
when calculating $E_{IJ} = T_{[I,J]}=\pa_{[J}T_{I]}$.  Thus we have
$$
2E_{IJ}=
\begin{pmatrix}
0&t(\pa_j\al -2b_j)&t\pa_\rho \al\\
t(2b_i- \pa_i\al)&t^2(\pa_jb_i
-\pa_ib_j)&t^2\pa_\rho b_i\\
-t\pa_\rho \al&-t^2\pa_\rho b_i&0 
\end{pmatrix}.
$$
We will use that this is of the form
\begin{equation}\label{Eform}
E_{IJ}=
\begin{pmatrix}
0&t\rho^{n/2+1}\cA  &t\rho^{n/2+1}\cA  \\
t\rho^{n/2+1}\cA  &t^2\rho^{n/2+1}\cA &t^2\rho^{n/2}\cA  \\ 
t\rho^{n/2+1}\cA  &t^2\rho^{n/2}\cA &0 \\ 
\end{pmatrix}.
\end{equation}

The conclusion of Proposition~\ref{formRt} for the components $\Rt_{0J}$ follows 
upon expanding the contractions in \eqref{Ric0} and using \eqref{Eform} and
our knowledge of $\gt^{IJ}$ and $\Gat_{IJK}$.  We indicate  
the details for $\Rt_{00}$; the cases $\Rt_{0j}$ and  
$\Rt_{0\infty}$ are similar.  In all cases, one knows ahead of time that
the powers of $t$ work out correctly. 

Setting $J=0$ in \eqref{Ric0} gives 
$$
t^2\Rt_{00}=t\gt^{IK}\pa_K E_{I0}-t\gt^{IK}\gt^{PQ}\Gat_{IKP}E_{Q0}
-t\gt^{IK}\gt^{PQ}\Gat_{0KP}E_{IQ}.
$$
It follows using $E_{I0}\in t^{1-\#(I)}\rho^{n/2+1}\cA$ and  
$t^{2-\#(IK)}\gt^{IK}\in \cA$ that when the first term on the right hand side
is expanded, all terms are in $\rho^{n/2+1}\cA$ except possibly those with 
$K=\infty$.  Since $E_{00}=0$, only the terms with $K=\infty$ and $I\neq 0$
need be considered.  But $t^{2-\#(I)}\gt^{I\infty} \in \rho \cA$ for
$I\neq 0$, so the first term is in $\rho^{n/2+1}\cA$.   

A term in the expansion of the second term is in $\rho^{n/2+1}\cA$ unless
$\Gat_{IKP}\notin t^{2-\#(IKP)}\cA$.  This gives that
two of $IKP$ must be between $1$ and $n$ and the
third must be $\infty$.  Each of the pairs $IK$, $PQ$ must be $0\infty$ or 
with both indices between $1$ and $n$.  $Q$ cannot be $0$.  There are no
such possibilities, so the second term is in $\rho^{n/2+1}\cA$.   

For the third term, we have $\Gat_{0KP}\in t^{1-\#(KP)}\cA$.  So for a term 
not to be in $\rho^{n/2+1}\cA$, it must be that $E_{IQ}$ is not in  
$t^{2-\#(IQ)}\rho^{n/2+1}\cA$, which gives that one of $IQ$ is between $1$
and $n$ and the other is $\infty$.  Again each of the pairs $IK$, $PQ$ must be
$0\infty$ or with both indices between $1$ and $n$.  So one of $KP$ must be
$0$ and the other between $1$ and $n$.  However, \eqref{c} shows that 
$\Gat_{0k0}$, $\Gat_{00k}\in \rho\cA$, so the third term is in
$\rho^{n/2+1}\cA$.  

For the components $\Rt_{IJ}$ in which neither $I$ nor $J$ is $0$, we use
\eqref{Ricci}.  It is straightforward to check using the
observations above about $\gt^{IJ}$ and $\Gat_{IJK}$ that for each of the
possibilities $IJ = ij$, $i\infty$, $\infty\infty$, each term on the right
hand side of 
\eqref{Ricci} is in $\cA$ when multiplied by the indicated power of
$\rho$.  
\end{proof}

In order to carry out the inductive perturbation analysis for
Theorem~\ref{even}, we need to extend \eqref{perturbricci} to the case
where the perturbations involve log terms.  If $0\leq m\in \Z$, we will say
that an expansion is $O^m$ if it can be written in the form
$\sum_{N\geq 0}u^{(N)}(\log|\rho|)^N$, where each $u^{(N)}$ is smooth, 
homogeneous of some degree in $t$, and $O(\rho^m)$.  Set $\cA^m 
= O^m\cap \cA$.     
The same calculations that gave
\eqref{perturbricci} give the following.
\begin{proposition}\label{perturblog}
Let $\gt$ have the form \eqref{refined} with $\al$, $b_i\in \cA^2$,  
$g_{ij}(x,\rho)=g_{ij}(x) +2P_{ij}\rho +\cA^2$.  Set 
$\gt_{IJ}' = \gt_{IJ} + \Phi_{IJ}$, where 
\begin{equation}\label{philog}
\Phi_{IJ}= 
\left(
\begin{matrix}
\phi_{00}&t\phi_{0j}&0\\
t\phi_{i0}&t^2\phi_{ij}&0\\
0&0&0
\end{matrix}
\right)
\end{equation}
with the $\phi_{IJ}\in \cA^m$, $m\geq 2$.  Then 
\begin{equation}\label{perturbriccilog} 
\begin{split}
t^2\Rt'_{00} &= t^2\Rt_{00} + 
(\rho\pa_\rho^2-\tfrac{n}{2}\pa_\rho)\phi_{00} +O^m\\
t\Rt'_{0i} &= t\Rt_{0i} +
(\rho\pa_\rho^2-\tfrac{n}{2}\pa_\rho)\phi_{0i}  
+ \tfrac12 \pa^2_{i\rho}\phi_{00} +O^m\\ 
\Rt'_{ij} &= \Rt_{ij} +
\left[\rho\pa_\rho^2+(1-\tfrac{n}{2})\pa_\rho\right]\phi_{ij}-\tfrac12  
g^{kl}\pa_\rho\phi_{kl}g_{ij}\\  
&\qquad\qquad +\tfrac12 (\nabla_j\pa_\rho\phi_{0i} + 
\nabla_i\pa_\rho\phi_{0j}) +P_{ij}\pa_\rho\phi_{00} +O^m\\  
t\Rt'_{0\nf} &= t\Rt_{0\nf} +
\tfrac12 \pa^2_\rho\phi_{00} +O^{m-1}\\
\Rt'_{i\nf} &= \Rt_{i\nf} +
\tfrac12 \pa^2_\rho\phi_{0i} +O^{m-1}\\
\Rt'_{\nf\nf} &= \Rt_{\nf\nf} - 
\tfrac12 g^{kl}\pa^2_\rho\phi_{kl} +O^{m-1} 
\end{split}
\end{equation}
\end{proposition}

We also need the analogue of \eqref{bcomp} for expansions with logs.  The
same reasoning as for \eqref{bcomp} shows that if $\gt\in\cM$ satisfies
for some $m\geq 2$ that $\Rt_{IJ}=O^{m-1}$ for $I$, $J\neq \infty$ and 
$\Rt_{I\infty}=O^{m-2}$, then \eqref{bcomp} holds except that the error
terms are all $O^m$ rather than $O(\rho^m)$.  

\medskip
\noindent
{\it Proof of Theorem~\ref{even}}. 
The analysis of Theorem~\ref{formal} and Proposition~\ref{simpler} leaves
us with a smooth $g_{ij}(x,\rho)$ determined modulo $O(\rho^{n/2})$ so that  
the metric $\gt$ defined by \eqref{spform} 
satisfies $\Rt_{0I}=0$, $\Rt_{ij}$, $\Rt_{j\infty}=O(\rho^{n/2-1})$,   
$\Rt_{\infty\infty}=O(\rho^{n/2-2})$.  (We will save the determination of 
the trace of the $\rho^{n/2}$ term in the expansion of $g_{ij}$ for the
next step of the induction.)  Each coefficient in the Taylor 
expansion of $g_{ij}$ through order $n/2-1$ is a natural tensorial
invariant of the initial metric $g$.  Let us write 
$(\frac{n}{2}-1)!\,\Rt_{ij}=\rho^{n/2-1}r_{ij}$ and 
$(\frac{n}{2}-1)!\,\Rt_{j\infty}=\rho^{n/2-1}r_{j\infty}$, 
where $r_{ij}$ and $r_{j\infty}$ are smooth and independent of $t$.  
The obstruction tensor is given by
$$
(\tfrac{n}{2}-1)!\,\cO_{ij}= 2^{1-n/2}c_n\tf{(r_{ij})}|_{\rho=0} 
$$
with $c_n$ as in \eqref{obsdef}.  
The values at $\rho =0$ of $r_{j\infty}$ and the trace $g^{kl}r_{kl}$
do depend on the choice of $g_{ij}$ at order $n/2$.  However, if we fix 
$g_{ij}$ to be its finite Taylor polynomial of order $n/2-1$, then these 
values can be expressed as natural tensorial invariants of the initial
metric $g$.

For specificity, fix $g_{ij}$ to be this Taylor polynomial.
Define $\gt'_{IJ}$ as in Proposition~\ref{perturblog} with $m=n/2$. 
It is clear from \eqref{perturbriccilog} that only the 
$\phi_{IJ} \mod \cA^{n/2+1}$ can affect the $\Rt_{IJ}$ at the next order.
So we may as well take  
\[
\begin{split}
&(\tfrac{n}{2})!\,\,\phi_{00}=\rho^{n/2}\left( \ka^{(1)}\log|\rho| +
  \ka^{(0)}\right)\\ 
&(\tfrac{n}{2})!\,\,\phi_{0i}=\rho^{n/2}\left( \mu_i^{(1)}\log|\rho| +
  \mu_i^{(0)}\right)\\  
&(\tfrac{n}{2})!\,\,\phi_{ij}=\rho^{n/2}\left( \la^{(1)}_{ij}\log|\rho| 
  +\la^{(0)}_{ij}\right)   
\end{split}
\]
with coefficients $\ka^{(N)}$, $\mu_i^{(N)}$, $\la_{ij}^{(N)}$, $N= 0,1$,
smooth and independent of $t$.  

The first equation of \eqref{perturbriccilog} and the requirement 
$\Rt_{00}'=O^{n/2}$ give
$(\rho\pa_\rho^2-\tfrac{n}{2}\pa_\rho)\phi_{00} =O^{n/2}$, which is easily 
seen to imply that $\phi_{00}=O^{n/2+1}$.  Similarly the second equation
of \eqref{perturbriccilog} and the requirement  
$\Rt_{0i}'=O^{n/2}$ give $\phi_{0i}=O^{n/2+1}$.  The third equation and the
requirement $\Rt_{ij}'= O^{n/2}$ give
\begin{equation}\label{ijcomp}
r_{ij} -\tfrac12 g^{kl}\la_{kl}^{(1)}g_{ij}\log|\rho| 
-\left(
\tfrac12 g^{kl}\la_{kl}^{(0)}+\tfrac{1}{n}g^{kl}\la_{kl}^{(1)}\right)g_{ij}  
+\la_{ij}^{(1)}=O^1.
\end{equation}
Clearly we must have $g^{kl}\la_{kl}^{(1)}=0$ at $\rho=0$.  Taking the
trace-free part then shows that 
$\tf(r_{ij})+\la^{(1)}_{ij}=0$ at $\rho=0$, so  
$$
c_n\la^{(1)}_{ij}|_{\rho =0}=-2^{n/2-1}(\tfrac{n}{2}-1)!\,\cO_{ij}.   
$$
Now taking the trace in \eqref{ijcomp} gives 
$g^{kl}r_{kl}=\frac{n}{2}g^{kl}\la_{kl}^{(0)}$ at $\rho =0$.  These
determinations are necessary and sufficient for $\Rt_{ij}'=O^{n/2}$.  The
trace-free part of $\la_{ij}^{(0)}$ is undetermined by the Einstein
condition, but is fixed by
the choice of $h_{ij}$:  $\tf(\la_{ij}^{(0)}|_{\rho = 0}) = h_{ij}$.   
Thus 
$$
\la_{ij}^{(0)}|_{\rho =0}= h_{ij}+\frac{2}{n^2}r_k{}^kg_{ij}.    
$$

Now $\gt'_{IJ}$ is determined mod $\cA^{n/2+1}$.  We fix the $\cA^{n/2+1}$ 
indeterminacy in $\gt_{00}'$, $\gt_{0i}'$ by taking $\phi_{00}$,
$\phi_{0i}=0$, so that $\gt_{IJ}'$ has the form \eqref{spform}.  Then    
$\Rt_{0I}'=0$, $\Rt_{ij}'=O^{n/2}$, $\Rt_{i\infty}'=O^{n/2-1}$,  
$\Rt_{\infty\infty}'=O^{n/2-2}$.  Also we know that
$\rho^2\Rt_{\infty\infty}'\in \cA$ by Proposition~\ref{formRt}.  Substituting
this information into the last line of \eqref{bcomp} with $m=n/2$ shows
that $\Rt_{\infty\infty}'=O^{n/2-1}$.  For reference in the next step, we
will need to know the leading term of $\Rt_{i\infty}'$.  This component is
given explicitly by the middle line of \eqref{evol} (where, however, $'$
denotes $\pa_\rho$).  Replacing $g_{ij}$ by 
$g_{ij}+\phi_{ij}$ in \eqref{evol} and recalling 
$\cO_j{}^j=0$, $\cO_{ij},{}^j=0$, one finds that 
\begin{equation}\label{iinfty}
(\tfrac{n}{2}-1)!\,\Rt_{i\infty}'=\rho^{n/2-1}\left(
r_{i\infty}+\frac12 h_{ij},{}^j- \frac{n-1}{n^2}r_j{}^j,{}_i\right) 
+O^{n/2}.
\end{equation}
This completes the $m=n/2$ step.  The metric $\gt'_{IJ}$ has the form
\eqref{spform} and its Ricci curvature satisfies
$\Rt_{0I}'=0$, $\Rt_{ij}'=O^{n/2}$, $\Rt_{i\infty}'$, 
$\Rt_{\infty\infty}'=O^{n/2-1}$.  

Now rename what was $\gt'_{IJ}$ to be a new $\gt_{IJ}$.
Proposition~\ref{formRt} shows that this new $\gt_{IJ}$ has $\rho \Rt_{ij}
\in \cA$, so we can write  
\begin{equation}\label{formRij}
(\tfrac{n}{2})!\,\Rt_{ij} = \rho^{n/2}\left( r_{ij}^{(1)}\log|\rho| + 
  r_{ij}^{(0)}\right) +O^{n/2+1} 
\end{equation}
with coefficients $r_{ij}^{(N)}$, $N=0,1$, smooth and independent of $t$.   
Now construct
a new $\gt'_{IJ}$ as in Proposition~\ref{perturblog} with $m=n/2+1$.  This
time we take
\[
\begin{split}
&(\tfrac{n}{2}+1)!\,\,\phi_{00}=\rho^{n/2+1}\left( \ka^{(1)}\log|\rho| +
  \ka^{(0)}\right)\\ 
&(\tfrac{n}{2}+1)!\,\,\phi_{0i}=\rho^{n/2+1}\left( \mu_i^{(1)}\log|\rho| +
  \mu_i^{(0)}\right)\\  
&(\tfrac{n}{2}+1)!\,\,\phi_{ij}=\rho^{n/2+1}\left( \la^{(1)}_{ij}\log|\rho|
  +\la^{(0)}_{ij}\right).   
\end{split}
\]
Referring to \eqref{perturbriccilog}, the requirement $\Rt'_{00}=O^{n/2+1}$ 
is equivalent to $\ka^{(1)}|_{\rho =0}=0$; no condition is imposed on
$\ka^{(0)}$.  However, the requirement $\Rt_{0\infty}'=O^{n/2}$ is
equivalent to $\phi_{00}=O^{n/2+2}$.  Thus the $\ka^{(N)}|_{\rho =0}$ are 
determined and we may as well take $\phi_{00}=0$.   Similarly, the
requirement 
$\Rt'_{0i} = O^{n/2+1}$ is equivalent to  $\mu_i^{(1)}|_{\rho =0}=0$, so we
take $\mu_i^{(1)}=0$.  The requirement  
$\Rt'_{i\infty} = O^{n/2}$ uniquely determines $\mu_i^{(0)}|_{\rho =0}$.    
Note that according to \eqref{iinfty}, we have $\mu_i^{(0)}|_{\rho =0}=0$
(and therefore  $\phi_{0i}=O^{n/2+2}$) if and  
only if $h_{ij},{}^j = D_i$, where  
\begin{equation}\label{Dformula} 
D_i =
2\left(\frac{n-1}{n^2}r_j{}^j,{}_i-r_{i\infty}\right)\Big{|}_{\rho =0}. 
\end{equation}
An easy computation from the third line of \eqref{perturbriccilog} using
\eqref{formRij} shows that the requirement $\Rt_{ij}' = O^{n/2+1}$ uniquely
determines the $\la_{ij}^{(N)}|_{\rho =0}$.  Now the third line of 
\eqref{bcomp} with $m=n/2+1$ shows that $\Rt_{\infty\infty}'=O^{n/2}$.  
This completes the $m=n/2+1$ step.  We have $\gt'_{IJ}\in \cM$
with  
$\Rt_{00}'$, $\Rt_{0i}'$, $\Rt_{ij}'=O^{n/2+1}$, $\Rt_{I\infty}'=O^{n/2}$.
Moreover, it always holds that $\gt_{00}'=2\rho + O^{n/2+2}$, and
$\gt_{0i}'=O^{n/2+2}$ if and only if $h_{ij},{}^j = D_i$.  If $\cO_{ij}=0$, 
then no log terms occur in any of the expansions, and $\gt_{IJ}'$ is
smooth.   

We now prove by induction on $m$ that there is a metric
$\gt\in \cM$, with $\al$, $b_i$, $g_{ij}$ in \eqref{refined}  
uniquely determined mod $O^{m}$, such that $\Rt_{IJ}=O^{m-1}$ for $I$,
$J\neq \infty$ and $\Rt_{I\infty}=O^{m-2}$.  We will also show that the
metric $\gt$ satisfies $\gt_{00}=2\rho + O^{m}$ and $\gt_{0i}=O^{m}$ 
if and only if $h_{ij},{}^j = D_i$.  Moreover, $\gt$ is smooth if and only
if $\cO_{ij}=0$.  We have established this for $m=n/2+2$.     

The argument for the induction step passing from $m$ to $m+1$ differs
depending on whether or not $m=n$.  First assume $m\neq n$, and of course
$m\geq n/2+2$.  
Proposition~\ref{formRt} implies that the Ricci curvature of $\gt_{IJ}$
takes the form
\begin{equation}\label{Riccicoeff}
\begin{split}
(m-1)!\,t^2\Rt_{00}&=\rho^{m-1}\sum_{N=0}^{M_1}r_{00}^{(N)}(\log|\rho|)^N
  +O^{m}\\  
(m-1)!\,t\Rt_{0i}&=\rho^{m-1}\sum_{N=0}^{M_2}r_{0i}^{(N)}(\log|\rho|)^N
  +O^{m}\\ 
(m-1)!\Rt_{ij}&=\rho^{m-1}\sum_{N=0}^{M_3}r_{ij}^{(N)}(\log|\rho|)^N +O^{m} 
\end{split}
\end{equation}
where 
$$
M_1=\left\lfloor \frac{2(m-2)}{n}-1\right\rfloor, \qquad
M_2=\left\lfloor \frac{2(m-1)}{n}-1\right\rfloor, \qquad
M_3=\left\lfloor \frac{2m}{n}\right\rfloor
$$
and the coefficient functions $r_{IJ}^{(N)}$ are smooth and independent
of $t$.   Define $\gt'_{IJ}$ as in Proposition~\ref{perturblog}.  The 
requirement that $\gt'\in \cM$ implies that the perturbation terms 
take the form 
\begin{gather}\label{perturbform}
\begin{gathered}
m!\,\phi_{00}=\rho^m\sum_{N=0}^{M_1}\ka^{(N)}(\log|\rho|)^N \\ 
m!\,\phi_{0i}=\rho^m\sum_{N=0}^{M_2}\mu_i^{(N)}(\log|\rho|)^N \\
m!\,\phi_{ij}=\rho^{m}\sum_{N=0}^{M_3} \la_{ij}^{(N)}(\log|\rho|)^N   
\end{gathered}
\end{gather}
with coefficients smooth and independent of $t$. (It is straightforward to
modify the argument to allow more general perturbations; for example, only
to require that $\gt'_{00}$, $\gt'_{0i} \in \cA$.  One finds that the
only solution is the one constructed here with $\gt' \in \cM$.)
Since $\gt'\in \cM$, the components $\Rt_{00}'$, $\Rt_{0i}'$, 
$\Rt_{ij}'$ take the same form as in \eqref{Riccicoeff} with coefficients 
$r_{IJ}'^{(N)}$ determined by \eqref{perturbriccilog}. 
Upon substituting the first line of \eqref{perturbform} into the first line
of \eqref{perturbriccilog}, one finds by considering the coefficients of 
$(\log|\rho|)^N$ inductively that  
requirement $\Rt_{00}'= O^{m}$ can be satisfied and uniquely determines
the coefficients $\ka^{(N)}$ at $\rho = 0$.  Since $M_1\leq M_2$, the term
$\frac12 \pa^2_{i\rho}\phi_{00}$ in the second line of
\eqref{perturbriccilog} can be written in the same form as that of
$t\Rt_{0i}'$ in 
\eqref{Riccicoeff}.  Substituting the second line of \eqref{perturbform} 
into \eqref{perturbriccilog}, one finds that the requirement 
$\Rt_{0i}'= O^{m}$ can be satisfied and uniquely determines 
the $\mu_i^{(N)}$ at $\rho =0$.  Similarly, using $m\neq n$, one finds that
the 
requirement $\Rt_{ij}'= O^{m}$ can be satisfied and uniquely determines   
the $\la_{ij}^{(N)}$ at $\rho =0$.  Now consider the Bianchi identities
\eqref{bcomp}.  The form of the components $\Rt_{I\infty}'$ is given by 
Proposition~\ref{formRt} and the induction hypothesis.  Substituting 
into \eqref{bcomp} and using the vanishing $\Rt_{IJ}'=O^m$ for $I$, $J \neq
\infty$, one finds successively that 
$\Rt_{0\infty}'=O^{m-1}$, $\Rt_{i\infty}'=O^{m-1}$, and (again using $m\neq 
n$) $\Rt_{\infty\infty}'=O^{m-1}$.  It is evident from Lemma~\ref{0comp} 
that if $\gt$ has the form \eqref{spform}, then $\gt'$ does too, and it is 
also evident that if $\gt$ is smooth, then $\gt'$ is too.  This concludes
the induction step in case $m\neq n$.

Finally consider the case $m=n$.  The argument for the determination of the
components $\phi_{00}$ and $\phi_{0i}$ mod $O^{n+1}$ is unchanged.  The
first two lines of \eqref{bcomp} then show that $\Rt_{0\infty}'$, 
$\Rt_{i\infty}'= O^{n-1}$.  For $m=n$ we have $M_3=2$.  Substituting into
the third line 
of \eqref{perturbriccilog}, one finds that one can uniquely choose the
coefficients $\la_{ij}^{(N)}$ for $N=1,2$ and the trace-free
part $\tf(\la_{ij}^{(0)})$ at $\rho =0$ to make $r_{ij}'^{(0)}$,
$r_{ij}'^{(1)}$, $\tf(r_{ij}'^{(2)})=O(\rho)$.  This leaves us with 
$$
\Rt_{ij}' = c \rho^{n-1} (\log|\rho|)^2g_{ij} + O^{n}
$$
for some smooth $c$, and there remains an indeterminacy of a multiple of 
$\rho^n g_{ij}$ in $\phi_{ij}$.  Proposition~\ref{formRt} and the induction 
hypothesis imply that 
$$
(n-2)!\,\Rt_{\infty\infty}' = \rho^{n-2}\sum_{N=0}^2
r_{\infty\infty}'^{(N)} (\log|\rho|)^N +O^{n-1}.
$$
Substituting this information into the third line of \eqref{bcomp}, one
finds that $c$, $r_{\infty\infty}'^{(2)}$, and $r_{\infty\infty}'^{(1)}$
all vanish at $\rho =0$.  Thus we have $\Rt_{ij}'=O^n$ and 
$(n-2)!\,\Rt_{\infty\infty}'=\rho^{n-2}r_{\infty\infty}'^{(0)} +O^{n-1}$.   
Now an inspection of the last line of \eqref{perturbriccilog} shows that
one can uniquely fix the $\rho^n g_{ij}$ indeterminacy in $\phi_{ij}$ to 
kill this last coefficient in $\Rt_{\infty\infty}'$, completing the
induction step.  In principle, this
argument allows the possibility that a log term might be created in
$\phi_{ij}$ even if $\gt$ is smooth, but the same reasoning as in the proof
for $n$ odd in Theorem~\ref{formal} shows that this potential log term does
not occur.   
\stopthm

We remark that similar arguments using the form of the perturbation
formulae \eqref{perturbriccilog} for the Ricci curvature show that the
metrics constructed in Theorems~\ref{2dim},   
\ref{odd} and \ref{even} are the only formal expansions of metrics
for $\rho>0$ or $\rho <0$ 
involving positive powers of $|\rho|$ and $\log|\rho|$ which are   
homogeneous of degree 2, Ricci-flat to infinite order, and in normal form.      

Convergence of formal series determined by Fuchsian problems such as these
in the case of real-analytic data have been considered by several authors.   
In particular, results of \cite{BaoG} can be applied to establish the 
convergence of the series occuring in Theorems~\ref{2dim} and \ref{odd}
(and also in Theorem~\ref{even} if the obstruction tensor vanishes) 
if $g$ and $h$ are real-analytic.  Convergence 
results including also the case when log terms occur in Theorem~\ref{even}
are contained in \cite{K}.       

Other treatments of various aspects of the construction and properties of
ambient metrics are contained in \cite{CG}, \cite{GP1}, \cite{BrG},
\cite{GP2}.    

\section{Poincar\'e Metrics}\label{poincaresection}

In this chapter we consider the formal theory for Poincar\'e metrics
associated 
to a conformal manifold $(M,[g])$.  We will see that even Poincar\'e
metrics are 
in one-to-one correspondence with straight ambient metrics, if both are in
normal form.  
Thus the formal theory for Poincar\'e metrics is a consequence of the  
results of Chapter~\ref{formalsect}.  The derivation of a Poincar\'e metric
from an ambient metric was described in \cite{FG}, and the inverse 
construction of an ambient metric as the cone metric over a Poincar\'e
metric was given in \S 5 of \cite{GrL}.   

The definition of Poincar\'e metrics is motivated by the example of the
hyperbolic metric   
$4(1-|x|^2)^{-2}g_e$ on the ball, where $g_e$ denotes the Euclidean
metric.  Let $(M,[g])$ be a smooth manifold of dimension $n\geq 2$ with a
conformal 
class of metrics of signature $(p,q)$.   
Let $M_+$ be a manifold with boundary satisfying $\pa M_+=M$.  Let $r$
denote 
a defining function for $\pa M_+$; i.e. $r\in C^\infty(M_+)$ satisfies
$r>0$ in the interior 
$M_+^\circ$, $r=0$ on $M$, and $dr\neq 0$ on $M$.   
A smooth metric $g_+$ on $M_+^\circ$ of signature $(p+1,q)$ 
is said to be conformally compact 
if $r^2g_+$ extends smoothly to $M_+$ and $r^2g_+|_M$ is nondegenerate 
(so $r^2g_+$ has signature $(p+1,q)$ also on $M$).     
A conformally compact metric is said to have conformal infinity $(M,[g])$
if $r^2g_+|_{TM}\in [g]$.  These conditions are independent of the choice
of defining function $r$.    

In the following, we will be concerned only with behavior near $M$.  We
will identify $M_+$ with an open neighborhood of $M\times \{0\}$ in      
$M\times [0,\infty)$, and $r$ will denote the coordinate in the second
factor.  We will use lower case Greek indices to label objects on $M_+$.  
Let $S_{\al\be}$ be a symmetric 2-tensor field in a open neighborhood of  
$M\times \{0\}$ in $M\times [0,\infty)$.  For $m\geq 0$, we will write 
$S= O^+_{\al\be}(r^m)$ if $S= O(r^m)$ and
$\operatorname{tr}_g(i^*(r^{-m}S))=0$ on  
$M$, where $i:M\rightarrow M\times[0,\infty)$ is $i(x)=(x,0)$ and $g$ is a
metric in the conformal class $[g]$.    

\begin{definition}\label{poincare}
A Poincar\'e metric for $(M,[g])$, where $[g]$ is a conformal class of
signature $(p,q)$ on $M$, is a conformally compact metric $g_+$ of
signature $(p+1,q)$ on 
$M_+^\circ$, where $M_+$ is an open neighborhood  
of $M\times \{0\}$ in $M\times [0,\infty)$, such that:
\begin{enumerate}
\item
$g_+$ has conformal infinity $(M,[g])$.   
\item
If $n$ is odd or $n=2$, then $\Ric(g_+)+ng_+$ vanishes to infinite order
along $M$.

\noindent
If $n\geq 4$ is even, then $\Ric(g_+)+ng_+\in O^+_{\al\be}(r^{n-2})$.    
\end{enumerate} 
\end{definition}
\noindent
Alternatively, one can consider metrics $g_-$ on $M_+^\circ$   
of signature $(p,q+1)$ such that $\Ric(g_-)-ng_-$ vanishes to the stated  
order.  This is equivalent to the above upon taking $g_-=-g_+$, with
$g\rightarrow -g$ and $(p,q)\rightarrow (q,p)$. 

If $g_+$ is a conformally compact metric, then 
$|dr/r|_{g_+}=|dr|_{r^2g_+}$ extends smoothly to $M_+$.  
The conformal transformation law for the curvature tensor shows that  
all sectional curvatures of $g_+$ approach $-|dr/r|^2_{g_+}$ at a boundary 
point (see \cite{M}).  We will say that a conformally compact metric $g_+$
is asymptotically hyperbolic if 
$|dr/r|_{g_+}=1$ on $M$.  A Poincar\'e metric is asymptotically hyperbolic.   

There is a normal form for asymptotically hyperbolic metrics analogous to
the normal form for pre-ambient metrics discussed in Chapter~\ref{setup}.  
\begin{definition}
An asymptotically hyperbolic metric $g_+$ 
is said to be in {\it normal form} 
relative to a metric $g$ in the conformal class if 
$g_+=r^{-2}\left(dr^2+g_r\right)$, where $g_r$ is a 1-parameter family
of metrics on $M$ of signature $(p,q)$ such that $g_0=g$.    
\end{definition}

\begin{proposition}\label{poincarenormal}
Let $g_+$ be an asymptotically hyperbolic metric on 
$M_+^\circ$ and let $g$ be a metric in the conformal class.  
Then there exists an open neighborhood $\cU$ of $M\times \{0\}$ in $M\times
[0,\infty)$ 
on which there is a unique diffeomorphism $\phi$ from $\cU$ into $M_+$
such that $\phi|_M$ is the identity map, and such that
$\phi^*g_+$ is in normal form  relative to $g$ on $\cU$.
\end{proposition}

We refer to \S5 of \cite{GrL} for the proof.  The proof in \cite{GrL} is for
the case $M=S^n$ and $g_+$ positive definite, but the same argument applies
in the general case, 
arguing as in the proof of Proposition~\ref{normalform} if $M$ is
noncompact.  

We will say that an asymptotically hyperbolic metric $g_+$ on   
$M_+^\circ$ is {\it even}  
if $r^2g_+$ is the restriction to $M_+$ of a smooth metric $h$ on an open
set $\mathcal V \subset M\times (-\infty,\infty)$ containing  
$M_+$, such that $\cV$ and $h$ are invariant under $r\rightarrow -r$.   
We will say that a diffeomorphism $\psi$ from $M_+$ into $M\times
[0,\infty)$ satisfying $\psi|_{M\times \{0\}}=Id$ is even if $\psi$  
is the restriction of a diffeomorphism of such an open set $\cV$ which 
commutes with $r\rightarrow -r$.  If $\psi$ is an even diffeomorphism and
$g_+$ is an even asymptotically hyperbolic metric, then $\psi^*g_+$ is also  
even.  An examination of the proof in \cite{GrL} 
shows that if $g_+$ in Proposition~\ref{poincarenormal} is even, then
$\phi$ is also even.  

The first main results of this chapter are the following analogues of
Theorems~\ref{main} and \ref{formal}.
\begin{theorem}\label{poincaremain}
Let $(M, [ g ])$ be a smooth manifold of
dimension $n\geq 2$, equipped with a conformal class.  
Then there exists an even Poincar\'e metric for $(M , [g])$.   
Moreover, if $g_+^1$ and $g_+^2$ are two even Poincar\'e metrics for
$(M,[g])$ defined on $(M_+^1)^\circ$, $(M_+^2)^\circ$, resp., then there
are open subsets   
$\cU^1\subset M_+^1$ and $\cU^2\subset M_+^2$ containing $M\times \{0\}$
and an even diffeomorphism $\phi:\cU^1\rightarrow \cU^2$ such that 
$\phi|_{M\times \{0\}}$ is the identity map, and such that:
\begin{enumerate}
\item[(a)] If $n = \dim M$ is odd, then $g_+^1 - \phi^*  
g_+^2$ vanishes to infinite order at every point of $M \times 
\{ 0 \}$.
\item[(b)] If $n = \dim M$ is even, then $g_+^1 - \phi^*  
g_+^2 = O_{\al\be}^+ ( r^{n-2})$. 
\end{enumerate}
\end{theorem}

\begin{theorem}\label{poincareformal}
Let $M$ be a  smooth manifold of dimension $n\geq 2$ and $g$ a
smooth metric on $M$.
\begin{enumerate}
\item[(A)] There exists an even Poincar\'e metric $g_+$ for $(M , [ g ])$  
which is in normal form relative to $g$.  
\item[(B)] Suppose that $g_+^1$ and $g_+^2$ are 
even Poincar\'e metrics for $(M , [ g ])$, both of which are in normal 
form relative to $g$.   
If $n $ is odd, then $g_+^1 - g_+^2$ 
vanishes to infinite order at every point of 
$M \times \{ 0 \}$.  
\noindent
If $n $ is even, then $g_+^1 - g_+^2
= O^+_{\al\be}(r^{n-2})$.
\end{enumerate}
\end{theorem}

Theorem~\ref{poincaremain} follows from Theorem~\ref{poincareformal} and
Proposition~\ref{poincarenormal} just as in the proof of
Theorem~\ref{main}.  
Theorem~\ref{poincareformal} will be proven as a consequence of  
Theorem~\ref{formal} after we establish the equivalence of straight   
ambient metrics and even Poincar\'e metrics in normal form.  

Let $(\cGt,\gt)$ be a straight pre-ambient space for $(M,[g])$.  It 
follows from 
Proposition~\ref{straightequiv} that $\|T\|^2$ vanishes exactly to first
order   
on $\cG\times \{0\}\subset \cGt$.  Therefore (shrinking $\cGt$ if
necessary), the hypersurface    
$\cH = \cGt\cap\{\|T\|^2=-1\}$ lies on one side of $\cG\times \{0\}$.
Since $\|T\|^2$ is homogeneous of degree 2 with respect to the dilations,
it follows (shrinking $\cGt$ again if necessary) 
that each dilation orbit in $\cGt$ on this side intersects $\cH$ 
exactly once.  
We extend the projection $\pi:\cG \rightarrow M$ to 
$\pi:\cGt\subset \cG\times \R \rightarrow M\times \R$ by acting in the
first 
factor.  Define $\chi:M\times \R\rightarrow M\times [0,\infty)$ 
by $\chi(x,\rho) = \left(x,\sqrt{2|\rho|}\right)$.  Then (shrinking $\cGt$
yet again if necessary), there is an open set $M_+$   
in $M\times [0,\infty)$ containing $M\times \{0\}$ so that  
$\chi\circ\pi|_\cH : \cH\rightarrow M_+^\circ$ is a diffeomorphism. 
In the following, we allow ourselves to shrink the domains of definition of 
$\gt$ and $g_+$ without further mention.
\begin{proposition}\label{ambienttopoincare}
If $(\cGt,\gt)$ is a straight pre-ambient space for $(M,[g])$ and $\cH$ and  
$M_+$ are as above, then
\begin{equation}\label{formg+}
g_+:=\left((\chi\circ\pi|_\cH)^{-1}\right)^*\gt
\end{equation}
is an even  
asymptotically hyperbolic metric with conformal infinity $(M,[g])$.   
If $\gt$ is in normal form relative to a metric $g\in [g]$, then  
$g_+$ is also in normal form relative to $g$.
Every even asymptotically hyperbolic metric $g_+$ with conformal infinity 
$(M,[g])$ is of the form \eqref{formg+} for some straight pre-ambient
metric $\gt$ for $(M,[g])$.  If $g_+$ is in normal form 
relative to $g$, then $\gt$ can be taken to be in normal form relative to
$g$, and in this case $\gt$ on $\{\|T\|^2\leq 0\}$ is uniquely determined
by $g_+$.    
\end{proposition}

\noindent
\begin{proof}
Choose a metric $g$ in the conformal class at infinity, with corresponding 
identification $\cG\times \R\cong \R_+\times M\times \R$.  Set
$\gt_{00}=a$, so $a=a(x,\rho)$ is smooth and homogeneous of degree 0 with
respect to 
$t$.  Then $\|T\|^2 = at^2$.  According to condition (2) of
Proposition~\ref{straightequiv}, $\gt$ 
is straight if and only if it has the form
\begin{equation}\label{straightform}  
\gt = adt^2 + tdtda +t^2h,
\end{equation}
where $h=h(x,\rho,dx,d\rho)$ is a smooth 
quadratic form defined in a neighborhood of $M\times \{0\}$ in 
$M\times \R$.  The initial condition $\iota^*\gt = {\bf g}_0$ is equivalent
to 
$a(x,0)=0$ and $h(x,0,dx,0)=g(x,dx)$, and nondegeneracy of $\gt$    
is equivalent to $\pa_\rho a\neq 0$ when $\rho =0$.      

Now $\cH = \{at^2=-1\}$.  On $\{a<0\}$, introduce new variables $u>0$ 
and $s>0$ by $a=-u^2$, $s=ut$.  Elementary computation shows that   
$adt^2+tdtda = s^2u^{-2}du^2 -ds^2$.  Therefore, in terms of these variables 
$\gt$ can be written
$$
\gt=s^2u^{-2}(h+du^2) -ds^2. 
$$  
This is the cone metric over the base $u^{-2}(h+du^2)$.
Thus every straight pre-ambient metric is a cone metric of this 
form.  In the new variables, $\cH$ is defined by the equation $s=1$.  
So the restriction of $\gt$ to $T\cH$ is the metric $u^{-2}(h+du^2)$.    

First suppose that $\gt$ is in normal form relative to $g$.  Comparing 
Lemma~\ref{explicitnormal} with \eqref{straightform}, one sees that this is
equivalent to the conditions that    
$\pa_\rho a=2$ and $h = g_\rho(x,dx)$, where $g_\rho$ is a smooth 
1-parameter family of metrics on $M$ satisfying $g_0=g$.  We obtain 
$a=2\rho$ and so $u =\sqrt{-2\rho}$.  By the definition of $\chi$, we see
that $u=r$ is the coordinate in the second factor of $M\times [0,\infty)$. 
Hence $g_+$ defined by \eqref{formg+} is just the metric  
\begin{equation}\label{g+form}
g_+= r^{-2}\left(dr^2 + g_{-\frac12 r^2}\right)
\end{equation} 
on $M\times [0,\infty)$.  Clearly $g_+$ is an even 
asymptotically hyperbolic metric with conformal infinity $(M,[g])$ in
normal form relative to $g$.  

In the general case, we have $g_+=u^{-2}(h+du^2)$ with  
$u=\sqrt{-a}$ and $h=h(x,\rho,dx,d\rho)$.  Since $r$ is given by   
$r=\sqrt{2|\rho|}$ and $a$ vanishes exactly to first order at $\rho   
=0$, we can write $u=rb(x,r^2)$ for a positive smooth function $b$. 
Using this and writing $\rho = \pm \frac12 r^2$, one sees easily that $g_+$  
is an even asymptotically hyperbolic metric with conformal infinity
$(M,[g])$.  

To see that every even asymptotically hyperbolic metric with conformal
infinity $(M,[g])$ is of this form, take $\gt$ to be given by 
\eqref{straightform} with $a=2\rho$ so that $g_+=r^{-2}(h+dr^2)$ with
$\rho = -\frac12 r^2$.  Writing 
$$
h=h_{ij}(x,\rho)dx^idx^j +2h_{i\infty}(x,\rho)dx^id\rho 
+h_{\infty\infty}(x,\rho)d\rho^2
$$
gives
$$
h+dr^2=h_{ij}(x,-\tfrac12 r^2)dx^idx^j 
-2r h_{i\infty}(x,-\tfrac12 r^2)dx^idr  
+\left[1+r^2h_{\infty\infty}(x,-\tfrac12 r^2)\right]dr^2.  
$$
An asymptotically hyperbolic metric $g_+$ is even if and only if 
the Taylor expansions of the $dx^idx^j$ coefficients of $r^2g_+$ have only
even terms, the Taylor expansions of the $dx^idr$ coefficients of $r^2g_+$
have only odd 
terms, and the $dr^2$ coefficient of $r^2g_+$ equals $1$ at $r=0$ and has
Taylor expansion with only even terms.  Clearly, any such $g_+$ can be
written in the form $r^{-2}(h+dr^2)$ for some $h$ with all
$h_{\al\be}(x,\rho)$ smooth.   

In general, there are many straight pre-ambient metrics $\gt$ such that
\eqref{formg+} gives the same metric $g_+$:  a given such $\gt$ can be
pulled back by a diffeomorphism of $\cGt$ which covers the
identity on $M\times \R$ but which smoothly rescales the $\R_+$-fibers by a
factor depending on the point in the base (and which restricts to the
identity on $\cG\times \{0\}$).  If, however, $g_+$ is in normal form and
$\gt$ is 
required also to be in normal form, then the determination $a=2\rho$ is
forced.  In this case, we saw above that if $h$ in \eqref{straightform} is
written $h = g_\rho(x,dx)$, then $g_+$ is given by \eqref{g+form}. 
Clearly, $g_+$ uniquely determines $h$ on $\{\rho \leq 0\}$.   
\end{proof}

We summarize the relation for metrics in normal form.  A straight
pre-ambient metric in normal form can be written 
\begin{equation}\label{normalgt}
\gt=2\rho dt^2 +2tdtd\rho +t^2 g_\rho
\end{equation}
for a 1-parameter family of metrics $g_\rho$ on $M$.  Under the change of
variables $-2\rho = r^2$, $s=rt$ on $\{\rho \leq 0\}$, $\gt$ takes the form 
\begin{equation}\label{gtg+}
\gt = s^2 g_+ -ds^2
\end{equation}
where $g_+$ is given by \eqref{g+form}. 

We remark that the hypothesis in 
Proposition~\ref{ambienttopoincare} that $\gt$ is straight is important:
if $\gt$ is not assumed to be straight, then the metric      
$g_+$ defined by \eqref{formg+} need not be asymptotically hyperbolic.   
Also, a cone metric $s^2 k -ds^2$ over any base $(N,k)$, where $N$ is a
manifold and $k$ a metric on $N$, is 
straight in the sense that for each $s>0$, $p\in N$, the curve
$\lambda\rightarrow (\lambda s,p)$ is a geodesic.  

The relation between the curvature of a cone metric and that of the base is
well-known.  We include the derivation for completeness.  
\begin{proposition}\label{curvrelation}
Let $g_+$ be a metric on a manifold $M_+$ of dimension $n+1$, and 
define a metric $\gt=s^2g_+ -ds^2$ on $M_+\times \R_+$.  Then 
\[
\begin{split}
\operatorname{Rm}(\gt) 
&= s^2 \left[\operatorname{Rm}(g_+) +g_+\owedge
g_+\right] \\ 
\Ric(\gt) &= \Ric(g_+) + ng_+\\
R(\gt) &= s^{-2} \left[R(g_+) + n(n+1)\right] 
\end{split}
\] 
Here $\operatorname{Rm}$ denotes the Riemann tensor viewed as a
covariant 4-tensor, and for symmetric 2-tensors $u$, $v$ we write 
$$
2(u\owedge v)_{IJKL}
=u_{IK}v_{JL}-u_{IL}v_{JK}+u_{JL}v_{IK}-u_{JK}v_{IL},  
$$
so that $g_+\owedge g_+$ is the curvature tensor of constant sectional
curvature $+1$.  
Also, tensors on $M_+$ are implicitly pulled back to $M_+\times \R_+$.  
\end{proposition}
\begin{proof} 
Set $s=e^y$; then $\gt=e^{2y}\left(g_+ -dy^2\right)$ is a conformal
multiple of a product metric.  Under a conformal change $\gt=e^{2y}h$, the
curvature tensor transforms by 
$$
\operatorname{Rm}(\gt)
=e^{2y}\left[\operatorname{Rm}(h) +2\Lambda \owedge h\right],
$$
where $\Lambda = -\nabla_h^2y+dy^2-\frac12 |dy|_h^2 \,h$.
For $h=g_+-dy^2$ we have 
$\operatorname{Rm}(h)=\operatorname{Rm}(g_+)$, $\nabla_h^2y=0$ and
$|dy|_h^2=-1$ so that $\Lambda = \frac12 (g_+ + dy^2)$.  Thus  
$2\Lambda \owedge 
h= (g_++dy^2)\owedge (g_+-dy^2)=g_+\owedge g_+$, which
gives the first equation.  The second and third equations follow by
contraction.
\end{proof}

Proposition~\ref{curvrelation} shows in particular that $\gt$ is flat if
and only if 
$g_+$ has constant curvature $-1$, $\gt$ is Ricci-flat if and only if  
$\Ric(g_+) = -ng_+$, and $\gt$ is scalar-flat if and only if $g_+$ has 
constant scalar curvature $-n(n+1)$.
Also, it is immediate from the equation for $\operatorname{Rm}(\gt)$
in Proposition~\ref{curvrelation} that 
$\partial_t\into \operatorname{Rm}(\gt)=0$; compare with Lemma~\ref{0comp} 
and the discussion before \eqref{curv}.  

\bigskip\noindent
{\it Proof of Theorem~\ref{poincareformal}}.
Given $M$ and $g$, Theorem~\ref{formal} shows that there exists an ambient 
metric in normal form relative to $g$, which by Propositions~\ref{simpler}
and \ref{tgeodesics} can be taken to be straight.  According to
Proposition~\ref{ambienttopoincare}, the metric $g_+$ defined by
\eqref{formg+} is even and in normal form relative to $g$, and \eqref{gtg+}
holds.  Proposition~\ref{curvrelation} shows   
that the Ricci curvature of $g_+$ vanishes to the correct order, so that 
$g_+$ is a Poincar\'e metric for $(M,[g])$.  This proves (A).
Part (B) follows similarly from 
part (B) of Theorem~\ref{formal}:  every even Poincar\'e metric $g_+$ in
normal form gives rise to a straight ambient metric in normal form by
\eqref{gtg+}.  
\stopthm

Thus for metrics in normal form, the formal asymptotics for Poincar\'e 
metrics is entirely equivalent to that for straight ambient metrics
under the change of variable $\rho = -\frac12 r^2$.  
In particular, via this change of variable one can easily write down the 
analogues for $g_+$ of \eqref{evol}, \eqref{deriv}, \eqref{traces} 
\eqref{prinpart}, and Proposition~\ref{onlyricci}.  
For future reference, we observe from \eqref{initial} that for $n\geq 3$,
the expansion  
for a Poincar\'e metric $g_+=r^{-2}\left(dr^2 +g_r\right)$ in normal form
begins:   
\begin{equation}\label{initialpoincare}
(g_r)_{ij} = g_{ij} - P_{ij}r^2 +\ldots
\end{equation}
(We have returned to our original notation of $g_r$ for the $M$-component
of $r^2g_+$.) 

There is a substantial literature concerning   
the asymptotics of Poincar\'e metrics.  See for example \cite{GrH1} and 
\cite{HSS} for direct analyses of the asymptotics up to order $n$ in $r$
from the Poincar\'e metric point of view, and \cite{R} for a study of 
the asymptotics in the context of general relativity.

As discussed in Chapter~\ref{formalsect}, Theorems~\ref{2dim}, \ref{odd} and
\ref{even} describe all formal expansions involving powers of $\rho$ and
$\log|\rho|$ of generalized ambient metrics in normal form.  
Proposition~\ref{ambienttopoincare} extends
with the same proof to the case of metrics with such expansions;
the relations between the expansions and regularity of $\gt$ and  
$g_+$ are determined by the change of variable $\rho = -\frac12 r^2$.  
This gives the following result describing all formal expansions 
$g_+=r^{-2}\left(dr^2+g_r\right)$ solving $\Ric(g_+)=-ng_+$.    
\begin{theorem}\label{infinitepoincareformal}
Let $M$ be a smooth manifold of dimension $n$, and let $g$ be a smooth
metric of signature $(p,q)$ and $h$ a smooth symmetric 2-tensor on $M$ such
that $g^{ij}h_{ij}=0$.   
\begin{itemize}
\item 
If $n=2$ and if $h$ satisfies $h_{ij,}{}^j=-\frac12 R_{,i}$, then there
exists an even Poincar\'e metric 
in normal form relative to $g$, such that 
$\tf\left(\pa^2_rg_r|_{r=0}\right)=h$.  These conditions uniquely
determine  
$g_r$ to infinite order at $r=0$.  (The solution $g_r$ can be written
explicitly; see Chapter~\ref{flat}.)  
\item 
If $n\geq 3$ is odd and if $h$ satisfies $h_{ij,}{}^j=0$, then  
there exists a Poincar\'e metric in normal form relative to $g$, such that  
$\tf\left(\pa_r^ng_r|_{r=0}\right)=h$.  These conditions uniquely determine      
$g_r$ to infinite order at $r=0$, and the solution satisfies
$\operatorname{tr}_g\left(\pa_r^ng_r|_{r=0}\right)=0$.  
\item
Let $n\geq 4$ be even.  
There is a natural pseudo-Riemannian invariant 
1-form ${\overline D}_i$ so that if $h$ satisfies $h_{ij,}{}^j=   
{\overline D}_i$, then there exists a 
metric $g_+=r^{-2}\left(dr^2 + g_r\right)$ satisfying $g_0=g$ and  
$\Ric(g_+)=-ng_+$ to  
infinite order, such that $g_r$ has an expansion of the form 
$$
g_r \sim \sum_{N=0}^\infty g_r^{(N)}\left(r^n\log r\right)^N,
$$
where each of the $g_r^{(N)}$ is a smooth family of symmetric 2-tensors on
$M$ even in $r$, and  
$\tf\left(\pa_r^ng^{(0)}_r\big{|}_{r=0}\right)=h$.      
These conditions uniquely determine the $g_r^{(N)}$ to 
infinite order at $r=0$.  The solution $g_r$ is smooth (i.e. $g_r^{(N)}$ 
vanishes to infinite order for $N\geq 1$) if and only if the obstruction
tensor $\cO_{ij}$ vanishes on $M$.   
\end{itemize}
\end{theorem}
\noindent
The 1-form ${\overline D}_i$ is given by  
$$
{\overline D}_i=\frac{(-2)^{-n/2}n!}{(n/2)!} D_i,
$$
where $D_i$ is the 1-form appearing in Theorem~\ref{even}.  

As for ambient metrics, if $g$ and $h$ are real-analytic, then the formal
series for $g_r$ converges.    

For $n$ odd, the even solution given by Theorem~\ref{poincareformal}
is the one determined by Theorem~\ref{infinitepoincareformal} upon taking
$h=0$.  

We close this chapter by describing an alternate interpretation of 
Poincar\'e metrics in terms of projective geometry.  This is motivated by  
the Klein model of hyperbolic space, which is the metric 
$$
(1-|x|^2)^{-1}\sum (dx^i)^2 + (1-|x|^2)^{-2}\left(\sum x^idx^i\right)^2 
$$
on the ball, whose geodesics are straight lines.  
Begin by recalling that two torsion free connections  
$\nabla$ and $\nabla'$ on the tangent bundle of a manifold are said 
to be projectively 
equivalent if they have the same geodesics up to parametrization, and that
this condition is equivalent to the requirement that their difference
tensor has the form $v_{(i}\delta_{j)}{}^k$ for some 1-form $v$ on $M$.    
Given a manifold with boundary $N$, consider the class of metrics
on $N^\circ$ which near the boundary have the form 
$h/\rho + d\rho^2/4\rho^2$, 
where $\rho$ is a defining function for the boundary and $h$ is a symmetric 
2-tensor which is smooth up to the boundary and for which $h|_{T\partial
N}$ is nondegenerate of signature $(p,q)$.  It is easily checked that this
class of metrics is independent of the choice of defining function so is
invariantly associated to $N$.  It is also easily checked that the
conformal class of $h|_{T\partial N}$ is also independent of the choice of
$\rho$, so should be called the conformal infinity by analogy with the
conformally compact case.  Elementary calculations show that if $\nabla$
denotes the Levi-Civita connection of such a metric, then there is a
connection $\nabla'$ which is smooth up to the boundary so that the
difference tensor $\nabla -\nabla'$ takes the form $v_{(i}\delta_{j)}{}^k$ 
where $v=-d\rho/\rho$.  (The equivalence class of $v$ modulo smooth 1-forms
is independent of $\rho$.)  Thus it makes sense to call such metrics
projectively compact; their geodesics are the same as those of a smooth
connection up to parametrization.  In the projective formulation, a 
Poincar\'e metric is then defined to be a projectively compact metric with
prescribed conformal infinity and with constant Ricci curvature $-n$.   
The change of variable $\rho =r^2$ transforms the class of projectively 
compact metrics to the class of even asymptotically hyperbolic metrics.  
As the construction of Proposition~\ref{ambienttopoincare} shows, the
projectively compact metrics are more directly related to 
pre-ambient metrics than are the asymptotically hyperbolic metrics.  A
defining function $\rho$ and the smooth structure on the space where the
projectively compact metrics live are induced   
directly from that on the ambient space without the introduction of the
square root.

\section{Self-dual Poincar\'e Metrics}\label{selfdualsect}

In \cite{LeB}, LeBrun showed using twistor methods that if $g$ is a
real-analytic metric on an oriented real-analytic 3-manifold $M$, then 
$[g]$ is the conformal infinity of a real-analytic self-dual Einstein  
metric on a deleted collar neighborhood of $M\times \{0\}$ in $M\times
[0,\infty)$, uniquely determined up to real-analytic diffeomorphism.
As mentioned in \cite{FG}, 
LeBrun's result can be proved as an application of our formal theory of 
Poincar\'e metrics.  In this chapter we show that the corresponding formal
power series 
statement is a consequence of Theorem~\ref{infinitepoincareformal}.  The 
self-duality condition can be viewed as providing a conformally invariant
specification of the formally undetermined term
$\pa_r^3g_r|_{r=0}$. 

Let $M$ be an oriented 3-manifold.  Give $M\times [0,\infty)$ 
the induced orientation determined by the requirement that $(\pa_r,
e_1,e_2,e_3)$ is positively oriented for a positively oriented frame 
$(e_1,e_2,e_3)$ for $T_pM$, $p \in M$.  Here, as above, $r$ denotes the  
coordinate in the $[0,\infty)$ factor.  Throughout this chapter we 
will use lower case Greek indices to label objects on $M\times [0,\infty)$, 
lower case Latin 
indices for objects on $M$, and a '0' index for the $[0,\infty)$ factor.   

Let $g$ be a metric on $M$ of signature $(p,q)$ and let 
$g_+=r^{-2}\left(dr^2 + g_r\right)$ be an asymptotically hyperbolic metric
on $M_+^\circ$ in normal form relative to $g$, where $M_+$ is an open
neighborhood of $M\times \{0\}$ in 
$M\times [0,\infty)$.  Set $\gb = dr^2 +g_r$.  We denote by $\mu$
the volume form of $g_r$ on $M$ (the dependence on $r$ is to be
understood), and by $\mb$ the volume form of $\gb$.  In terms of a
positively oriented coordinate system on $M$, these are given by 
$\mu_{ijk}=\mb_{0ijk}=\sqrt{|\det{g_r}|}\,\ep_{ijk}$, where $\ep_{ijk}$
denotes the sign of the permutation.  We use $\gb$ to raise and lower 
Greek indices and $g_r$ to raise and lower Latin indices.  

The metric $\gb$ and the orientation determine a Hodge $*$ operator on
$\Lambda^2(M_+)$ which is conformally invariant so agrees with the $*$
operator for $g_+$.  This induces an operator $*:\mathcal{W}\rightarrow
\mathcal{W}$, where $\mathcal{W}\subset \otimes^4T^*M_+$ denotes the bundle
of algebraic Weyl tensors; i.e. 4-tensors with curvature tensor 
symmetry which are trace-free with respect to $\gb$ (or equivalently
$g_+$).  If $W_{\al\be\ga\delta}$ is a section of $\mathcal{W}$, then    
\begin{equation}\label{star}
(*W)_{\al\be\ga\delta} = \tfrac12 \,\mb_{\al\be}{}^{\rho\sigma} 
W_{\rho\sigma\ga\delta}.
\end{equation}
One has $*^2=(-1)^q$.  We assume throughout the rest of this chapter that
$q$ is even (i.e. $q=0$ or $q=2$) so that $*^2=1$.    
Then $\mathcal{W}=\mathcal{W}^+\oplus\mathcal{W}^-$, where
$\mathcal{W}^\pm$ denote the $\pm 1$ eigenspaces of $*$.   

We will use an alternate description of the bundles $\mathcal{W}^\pm$.
For $(p,r)\in M_+\subset M\times \R$, define a map
$s:\mathcal{W}_{(p,r)}\rightarrow  \odot^2T^*_pM$ by  
$$
s(W)_{ij}=W_{0i0j}.
$$
Observe that $s(W)$ is trace-free with respect to $g_r$ since $W$ is
trace-free with respect to $\gb$:
$$
0=\gb^{\al\be}W_{0 \al 0 \be}
=(g_r)^{ij}W_{0i0j}+W_{0000}=(g_r)^{ij}s(W)_{ij}.  
$$
We denote by $\odot_0^2T^*M$ the bundle on $M_+$ whose fiber at $(p,r)$ 
consists of the symmetric 2-tensors in the $M$ factor which are 
trace-free with respect to $g_r(p)$.  No confusion should arise with this
notation since this bundle restricts to $M\subset M_+$ to the bundle for
which this notation is usually used, and the context will make clear 
whether the bundle is considered as defined on $M$ or $M_+$.   
\begin{lemma}\label{weyl}
$s|_{\mathcal{W}^\pm}:\mathcal{W}^\pm\rightarrow \odot_0^2T^*M$  
is a bundle isomorphism.  (Here $s$ is restricted to either the bundle 
$\mathcal{W}^+$ or $\mathcal{W}^-$.)
\end{lemma}
\begin{proof}
Since the bundles 
$\mathcal{W}^\pm$ and $\odot_0^2T^*M$ all have rank 5, it suffices to show
that $s|_{\mathcal{W}^\pm}$ is injective.  If $W\in \mathcal{W}$, then    
$$
0=\gb^{\al\be}W_{\al i \be j}
=(g_r)^{kl}W_{kilj}+W_{0i0j}.    
$$
So if $s(W)=0$, then $(g_r)^{kl}W_{kilj}=0$.  But then $W_{kilj}$ is a 
trace-free tensor in 3 dimensions with curvature tensor symmetry, so 
$W_{kilj}=0$.  

If also $W\in \mathcal{W}^\pm$, then \eqref{star} gives
$$
W_{0jkl} = \pm\tfrac12 \mb_{0j}{}^{\rho\sigma}
W_{\rho\sigma kl}.  
$$
The terms all vanish for which $\rho$ or $\sigma$ is 0, and the  
terms vanish in which neither $\rho$ nor $\sigma$ is 0 since $W_{ijkl}=0$.   
Thus all components $W_{0jkl}=0$.  This accounts for all possible nonzero
components of $W$, so $W=0$ and $s|_{\mathcal{W}^\pm}$ is injective.      
\end{proof}

We use Lemma~\ref{weyl} to represent sections of
$\mathcal{W}^\pm$ as one-parameter families of tensors on $M$.   
Let $\Wb$ denote the Weyl tensor of $\gb$ and $\Wb^\pm$ its $\pm$ 
self-dual parts.  Then for either choice of $\pm$, $s(\Wb^\pm)$ can be
regarded as a one-parameter family 
parametrized by $r$ of sections of the fixed bundle $\odot^2T^*M$ on $M$.  
Thus it makes sense to differentiate this family with respect to $r$.  The
following proposition gives an expression for the first derivative 
of $s(\Wb^\pm)$ in terms of derivatives of $g_r$.  We use $'$ to denote  
differentiation with respect to $r$ and we sometimes suppress the subscript
$r$ on $g_r$.  
\begin{proposition}\label{selfdualderiv}
Let $g_r$ be a one-parameter family of metrics on $M$ satisfying $g_r'=0$
at $r=0$, and let $\gb = dr^2 +g_r$.  Then at $r=0$ we have:    
$$
s(\Wb^\pm)_{ij}'=-\tfrac18 \tf (g_{ij}''')    
\pm\tfrac14 \nabla_lg_{k(i}''\mu_{j)}{}^{kl}.
$$
On the right hand side, $\tf$, $\nabla$ and $\mu$ are with 
respect to the initial metric $g_0$.     
\end{proposition}
\noindent
We remark that the star operator on 2-forms on $M$ with respect to $g$ is  
given by $(*\eta)_i=\frac12 \mu_i{}^{kl}\eta_{kl}$.  Thus the last term
above can 
be interpreted as a multiple of $\operatorname{Sym}*dg''$, where $d$ acts 
on $g''$ as a 1-form-valued 1-form, $*$ acts on the 2-form indices   
generated by $d$, and $\operatorname{Sym}$ symmetrizes over
the two indices of $*dg''$.  
\begin{proof}
Specializing the indices in 
\begin{equation}\label{weylbar}
\Wb_{\al\be\ga\delta}=\overline{R}_{\al\be\ga\delta}
-\left(\overline{P}_{\al\ga}\gb_{\be\delta}
-\overline{P}_{\al\delta}\gb_{\be\ga}
-\overline{P}_{\be\ga}\gb_{\al\delta}
+\overline{P}_{\be\delta}\gb_{\al\ga}\right)
\end{equation}
and using the form of $\gb$ gives 
$$
\Wb_{0i0j}=\overline{R}_{0i0j}
-\left(\overline{P}_{00}g_{ij}+\overline{P}_{ij}\right).
$$
But we know that the left hand side is trace-free in $ij$ with respect to
$g_r$, 
so taking the trace-free part with respect to $g_r$ and recalling the
definition of $\overline{P}$ gives 
$$
\Wb_{0i0j}=\tf\left(\overline{R}_{0i0j}
-\tfrac12 \overline{R}_{ij}\right).
$$
Now $\overline{R}_{ij}=\gb^{\al\be}\overline{R}_{\al i\be j}= 
\overline{R}_{0i0j}+g^{kl}\overline{R}_{ikjl}$, giving
\begin{equation}\label{firstweyl}
\Wb_{0i0j}=\tfrac12 \tf\left(\overline{R}_{0i0j}
-g^{kl}\overline{R}_{ikjl}\right).  
\end{equation}
Next, \eqref{star} gives $(*\Wb)_{0i0j}
=\frac12 \mb_{0i}{}^{\rho\sigma}\Wb_{\rho\sigma 0j}
=\frac12 \mb_{0i}{}^{kl}\Wb_{kl 0j}$.  Specializing 
\eqref{weylbar} again, substituting, and simplifying gives easily 
$$
(*\Wb)_{0i0j}=\tfrac12 \mb_{0i}{}^{kl}
\left(\overline{R}_{0jkl}-\overline{R}_{0k}g_{jl}\right).  
$$
However, $\mb_{0i}{}^{kl}\overline{R}_{0k}g_{jl}
=\mb_{0i}{}^k{}_j\overline{R}_{0k}$ is skew in $ij$ and 
$(*\Wb)_{0i0j}$ is symmetric in $ij$, so we can symmetrize to obtain 
\begin{equation}\label{secondweyl}
(*\Wb)_{0i0j}=\tfrac12 \mb_{0(i}{}^{kl}\overline{R}_{j)0lk}.
\end{equation}

Now it is straightforward to calculate the curvature tensor of 
$\gb=dr^2+g_r$; one obtains the curvature tensor of a product metric plus
extra terms involving $r$-derivatives of $g_r$.  The result (cf. 
Gauss-Codazzi equations) is:  
\[
\begin{split}
\overline{R}_{ijkl}&=R_{ijkl}
+\tfrac14 \left(g'_{il}g'_{jk}-g'_{ik}g'_{jl}\right),\\  
\overline{R}_{0jkl}&=
\tfrac12 \left(\nabla_lg'_{jk}-\nabla_kg'_{jl}\right),\\  
\overline{R}_{0i0j}&=-\tfrac12 \,g_{ij}'' +\tfrac14 \,g^{rs}g_{ir}'g_{js}', 
\end{split}
\]
where $R_{ijkl}$ and $\nabla$ are with respect to $g_r$ as a metric on $M$
with $r$ fixed.  Substituting into \eqref{firstweyl} and \eqref{secondweyl}
gives
\[
\begin{split}
\Wb_{0i0j}&=-\tfrac14 \tf\left(g_{ij}''+2R_{ij}
-\tfrac12 g^{kl}g_{kl}'g_{ij}'\right),\\ 
(*\Wb)_{0i0j}&=\tfrac12\, \nabla_lg'_{k(i}\mu_{j)}{}^{kl}.
\end{split}
\]

Now differentiate these equations with respect to $r$ at $r=0$.  Since
$g'=0$ at $r=0$, it follows that raising indices, taking the trace-free
part and applying 
$\nabla$ all commute with the differentiation.
Similarly, the derivatives 
of $R_{ij}$ and $\mu_{jkl}$ vanish at $r=0$.  Thus we have at $r=0$:
\[
\begin{split}
\Wb_{0i0j}'&=-\tfrac14 \tf\left(g_{ij}'''\right),\\  
(*\Wb)_{0i0j}'&=\tfrac12 \nabla_lg''_{k(i}\mu_{j)}{}^{kl}. 
\end{split}
\]
Adding and subtracting yields the stated formula.     
\end{proof}

On an oriented pseudo-Riemannian 3-manifold of (arbitrary) signature
$(p,q)$, the Cotton tensor $C_{jkl}$  
can be reinterpreted as a trace-free symmetric 2-tensor.  Define  
$$
\mathcal{C}_{ij}=\mu_i{}^{kl}C_{jkl}.  
$$
We sometimes write $\mathcal{C}(g)$ to indicate the underlying metric.
The fact that $C_{[jkl]}=0$ implies that $\mathcal{C}_{ij}$ is trace-free: 
$$
g^{ij}\mathcal{C}_{ij}=\mu^{jkl}C_{jkl}=0,
$$
and the fact that $C_{jkl}$ is trace-free implies that $\mathcal{C}_{ij}$
is symmetric:
$$
\mu^{ijm}\mathcal{C}_{ij}=\mu^{ijm}\mu_i{}^{kl}C_{jkl}
=(-1)^q\left(g^{jk}g^{ml}-g^{jl}g^{mk}\right)C_{jkl}=0.  
$$
The fact that $C_{jkl,}{}^j=0$ implies that $\mathcal{C}_{ij,}{}^j=0$.   

We will say that a Poincar\'e metric $g_+$ is a $\pm$ self-dual Poincar\'e
metric if 
$W^\mp(g_+)$ vanishes to infinite order.   
The main result of this chapter is the following theorem.
\begin{theorem}\label{selfdualthm}
Let $g$ be a smooth metric of signature $(p,q)$ on a manifold $M$ of
dimension 3, with $q$ even. 
\begin{enumerate} 
\item
A Poincar\'e metric $g_+=r^{-2}\left(dr^2 +g_r\right)$ in normal form
relative to $g$ is a $\pm$ self-dual Poincar\'e metric if and only if  
\begin{equation}\label{cotton}
\left(\pa_r^3g_r\right)|_{r=0} = \pm 2 \mathcal{C}(g). 
\end{equation}
\item
There exists a $\pm$ self-dual Poincar\'e metric $g_+$ in normal form
relative to 
$g$, and such a $g_+$ is uniquely determined to infinite order. 
\end{enumerate}
\end{theorem}
\noindent
\begin{proof}
If $g_+$ is a Poincar\'e metric in normal form relative to $g$, then 
Theorem~\ref{infinitepoincareformal} shows that $g_r'''|_{r=0}$ is 
trace-free and  
\eqref{initialpoincare} shows that $g_r''|_{r=0}=-2P_{ij}$. 
Thus Proposition~\ref{selfdualderiv} shows that at $r=0$ we have: 
\begin{equation}\label{weylcotton}
\begin{split}
s(\Wb^\mp)_{ij}'&=-\tfrac18 g_{ij}'''    
\pm\tfrac12 \nabla_lP_{k(i}\mu_{j)}{}^{kl}\\
&=-\tfrac18 g_{ij}'''    
\pm\tfrac14 \mu_{(j}{}^{kl}C_{i)kl}\\
&=-\tfrac18 g_{ij}'''    
\pm\tfrac14 \mathcal{C}_{ij}
\end{split}
\end{equation}
This must vanish if $W^\mp$ vanishes to infinite order, which establishes
the necessity of \eqref{cotton}.    

For the converse, we must show that $\Wb^\mp$ vanishes to infinite order if 
\eqref{cotton} holds.  For this we use the fact that $W^\mp(g_+)$ satisfies
a first order equation as a consequence of the Bianchi identity and the
fact that $g_+$ is Einstein.  The contracted second Bianchi identity in
dimension 4 can be written $\nabla^\al W_{\al\be\ga\delta}
=-C_{\be\ga\delta}$.   The Cotton tensor 
vanishes for Einstein metrics, so it follows that
$\nabla^\al W_{\al\be\ga\delta}=O(r^\infty)$,
where $W$, $\nabla$ and index raising are with respect 
to $g_+$.  Now 
rewrite $\nabla^\al W_{\al\be\ga\delta}$ in terms of $\overline{\nabla}$   
and $\Wb$.  The connections of $g_+$ and $\gb$ are related by   
$$
\Gamma^\al_{\be\ga}(g_+)-\Gamma^\al_{\be\ga}(\gb)=
-r^{-1}\left(r_\be \delta^\al{}_\ga+r_\ga \delta^\al{}_\be
-\gb^{\al\lambda}r_\lambda \gb_{\be\ga}\right).
$$
Using this to transform the covariant derivative shows that
$$
\nabla^\al W_{\al\be\ga\delta}=
r^2\,\overline{\nabla}^\al W_{\al\be\ga\delta}
+r r^\al W_{\al\be\ga\delta},
$$
where on the right hand side the index is raised using $\gb$.  
Now substitute $W_{\al\be\ga\delta}=r^{-2}\Wb_{\al\be\ga\delta}$ to
conclude that  
$$
r\,\overline{\nabla}^\al \Wb_{\al\be\ga\delta}
- r^\al \Wb_{\al\be\ga\delta}=O(r^\infty). 
$$
Since $\overline{\nabla}$ commutes with $*$, this
equation also holds for $\Wb^\mp$.   
Using $r^\al=\delta^\al{}_0$, we therefore obtain
\begin{equation}\label{dirac}
r\,\overline{\nabla}^\al \Wb^\mp_{\al\be\ga\delta}
-  \Wb^\mp_{0\be\ga\delta}=O(r^\infty).
\end{equation}

Suppose we know that $\Wb^\mp = O(r^s)$ for some $s$.  
Write $\Wb^\mp = r^sV$ and substitute into \eqref{dirac}.  One obtains 
$(s-1)r^sV_{0\be\ga\delta}+O(r^{s+1})=O(r^\infty)$.  So if
$s\neq 1$, then $V_{0\be\ga\delta}=0$ at $r=0$.  
Taking $\be\ga\delta = i0j$ shows that $s(V)=0$ at $r=0$, which implies 
$V=0$ at $r=0$ by Lemma~\ref{weyl}.  Thus one concludes that
$\Wb^\mp=O(r^{s+1})$ if $s\neq 1$.

We use this observation inductively.  Begin with $s=0$ to conclude that
$\Wb^\mp = O(r)$ for any Poincar\'e metric in normal form.  If also 
\eqref{cotton} holds, then \eqref{weylcotton}
shows that $s(\Wb^\mp)'=0$ at $r=0$, which gives 
$\Wb^\mp=O(r^2)$.  Now the induction proceeds to all higher orders and one
concludes that $\Wb^\mp=O(r^\infty)$ as desired.

Part (2) is an immediate consequence of (1) and 
Theorem~\ref{infinitepoincareformal}.  For any metric $g$ on $M$, the
tensor $\mathcal{C}(g)$ is trace-free and divergence-free.  Take 
$h=\pm 2 \mathcal{C}(g)$ in Theorem~\ref{infinitepoincareformal} to
conclude the existence of a Poincar\'e metric in normal form relative to
$g$ which satisfies \eqref{cotton}.  
Part (1) implies that this Poincar\'e metric is $\pm$ self-dual to infinite 
order.  Uniqueness follows from the necessity of \eqref{cotton} in 
(1) together with the uniqueness statement in 
Theorem~\ref{infinitepoincareformal}.  
\end{proof} 

It follows by putting metrics into normal form that two
self-dual Poincar\'e metrics with the same conformal infinity   
agree up to diffeomorphism and up to terms vanishing to infinite
order.

If $g$ in Theorem~\ref{selfdualthm} is real-analytic, then $\mathcal{C}(g)$
is also real-analytic, so the series for $g_r$ converges as 
mentioned earlier.  Thus one recovers LeBrun's result.   

We remark that for definite signature, one can prove a version of part (1)
of Theorem~\ref{selfdualthm} concluding $\pm$ self-duality in an open set 
near the boundary for metrics which are Einstein in  
an open set and which satisfy \eqref{cotton}.  This is stated as Theorem 3
of \cite{Gr2}.    

\section{Conformal Curvature Tensors}\label{cct}

In this chapter we study conformal curvature tensors of a
pseudo-Riemannian metric $g$.  These are defined in terms of
the covariant derivatives of the curvature tensor of an ambient metric in 
normal form relative to $g$.  Their transformation laws under 
conformal change are given in terms of the action of a subgroup of the
conformal group $O(p+1,q+1)$ on tensors.  We assume throughout this chapter
that $n\geq 3$.  

Let $g$ be a metric on a manifold $M$.  By Theorem~\ref{formal}, there is
an ambient metric in normal form relative to $g$, which 
by Proposition~\ref{straightprop} we may take to be
straight.  Such a metric takes the form \eqref{spform} on a  
neighborhood of $\R_+ \times M\times \{0\}$ in $\R_+ \times M\times \R$.
Equations \eqref{evol} determine the 1-parameter family of metrics
$g_{ij}(x,\rho)$ on $M$ in terms of the initial metric to
infinite order for $n$ odd and modulo $O(\rho^{n/2})$ for $n$ even, except
that also $g^{ij}\pa_{\rho}^{n/2}g_{ij}|_{\rho = 0}$ is determined for $n$
even.  Each of the  
determined Taylor coefficients is a 
natural invariant of the initial metric $g$.   

We consider the curvature tensor and its covariant derivatives for
such an ambient metric $\gt$.  In general, 
we denote the curvature tensor of a pre-ambient metric by $\Rt$, with 
components $\Rt_{IJKL}$.  Its $r$-th covariant derivative will be denoted 
$\Rt^{(r)}$, with components $\Rt^{(r)}_{IJKL,M_1\cdots M_r}$.
Sometimes the superscript ${}^{(r)}$ will be omitted when the list of
indices makes clear the value of $r$.  Using \eqref{cn}, it is  
straightforward to calculate the curvature tensor of a metric of the form 
\eqref{spform}.  One finds that $\Rt_{IJK0}=0$ (another derivation of 
this is given in Proposition~\ref{Tcontract} below) and that the other
components are given by: 
\begin{equation}\label{curv}
\begin{split}
\Rt_{ijkl} = &t^2\left[R_{ijkl} +\frac12(g_{il}g_{jk}'+g_{jk}g_{il}'
-g_{ik}g_{jl}'-g_{jl}g_{ik}')+\frac{\rho}{2}(g_{ik}'g_{jl}'-g_{il}'g_{jk}')\right]\\   
\Rt_{\nf jkl}= &\frac12 t^2\left[ \nabla_lg_{jk}' - \nabla_kg_{jl}'\right]\\ 
\Rt_{\nf jk\nf}=&\frac12 t^2\left[ g_{jk}'' -
  \frac12g^{pq}g_{jp}'g_{kq}'\right]\\ 
\end{split}
\end{equation}
Here $'$ denotes $\pa_\rho$, $R_{ijkl}$ denotes the curvature tensor of the
metric $g_{ij}(x,\rho)$ with $\rho$ fixed, and $\nabla$ its Levi-Civita
connection.  Components of the ambient covariant derivatives of curvature
can then be calculated recursively starting with these formulae and using 
\eqref{cnr}.  

Fix a straight ambient metric $\gt$ in normal form relative to $g$.  
We construct tensors on $M$ from the covariant derivatives of curvature of
$\gt$ as follows.  Choose an order $r\geq 0$ of covariant differentiation.  
Divide the set of symbols $IJKLM_1\cdots M_r$ into three disjoint
subsets labeled $\cS_0$, $\cS_M$ and $\cS_\nf$.  Set the indices in $\cS_0$
equal to $0$, those in $\cS_\nf$ equal to $\nf$, and let those in $\cS_M$
correspond to $M$ in the decomposition $\R_+ \times M\times \R$.  
(In local coordinates, the indices in $\cS_M$ vary between $1$ and $n$.)
Evaluate the resulting component $\Rt_{IJKL,M_1\cdots M_r}$ at $\rho = 0$
and $t=1$.  
This defines a tensor on $M$ which we denote by
$\Rt^{(r)}_{\cS_0,\cS_M,\cS_\nf}$.  
(Recall that the submanifold $\{\rho = 0,\,\, 
t=1\}$ of $\R_+\times M \times \R$ can be invariantly described as the
image of $g$ viewed as a section of $\cG$.)

Consider for example the case $r=0$.  Since
$\Rt_{IJK0}=0$ as noted above, we must choose $\cS_0 = \emptyset$ in order
to get something nonzero.  If we choose also $\cS_\nf = \emptyset$ so that 
$\cS_M = \{IJKL\}$, then the resulting tensor
$\Rt^{(0)}_{\cS_0,\cS_M,\cS_\nf}$ can be determined by setting $\rho = 0$
and $t=1$ 
in the first line of \eqref{curv}.  Recall from \eqref{initial},
\eqref{Ptensor} that $g'_{ij}|_{\rho = 0} = 2P_{ij}$.  This gives 
$\Rt_{ijkl}|_{\rho = 0}=t^2 W_{ijkl}$, so in this case the tensor
$\Rt^{(0)}_{\cS_0,\cS_M,\cS_\nf}$ is the 
Weyl tensor of $g$.  Similarly one finds that
$\Rt_{\nf jkl}|_{\rho = 0,\,t=1} = C_{jkl}$ is the Cotton tensor of $g$.  
The tensors for all other possibilities for the subsets $\cS_0$, 
$\cS_M$, $\cS_\nf$ are 
determined from these and from $\Rt_{\nf jk\nf}$ by the usual symmetries of
the curvature tensor.  Now $\Rt_{\nf jk\nf}$ is given by the last line of
\eqref{curv}.  When $n=4$, the trace-free part of $g''_{jk}|_{\rho =0}$
depends on the specific ambient metric which has been chosen and 
is not determined solely by $g$, so the same holds for $\Rt_{\nf jk\nf}$.
However, when $n\neq 4$, one finds using the first line of \eqref{deriv}
and the last line of \eqref{curv} that  
$\Rt_{\nf jk\nf}|_{\rho =0,\,t=1}=-(n-4)^{-1}B_{jk}$.  Thus the conformal 
curvature tensors which arise for $r=0$ are precisely the Weyl, Cotton, and
Bach tensors of $g$:
\begin{equation}\label{curv0}
\Rt_{ijkl}|_{\rho = 0, t=1}=W_{ijkl},\quad
\Rt_{\nf jkl}|_{\rho = 0,\,t=1} = C_{jkl}\quad
\Rt_{\nf jk\nf}|_{\rho =0,\,t=1}=-\frac{B_{jk}}{n-4}.
\end{equation}

Iterated covariant derivatives  of $\Rt$ can be 
can be calculated recursively using \eqref{cnr}.  For example, one obtains  
$$
\Rt_{ijkl,m}|_{\rho =0,\, t=1} = V_{ijklm},\qquad 
\Rt_{\nf jkl,m}|_{\rho =0,\,t=1} = Y_{jklm},
$$
where
\begin{equation}\label{Vform}
\begin{split}
V_{ijklm}& = W_{ijkl,m} +g_{im}C_{jkl}-g_{jm}C_{ikl}+g_{km}C_{lij}
-g_{lm}C_{kij}\\
Y_{jklm} &= C_{jkl,m} - P_m{}^iW_{ijkl} 
+ (n-4)^{-1}(g_{km}B_{jl} - g_{lm}B_{jk}), 
\end{split}
\end{equation}
and we assume that $n\neq 4$ for $\Rt_{\nf jkl,m}$.  

The usual identities satisfied by covariant derivatives of
curvature imply identities and relations amongst the conformal curvature
tensors.  For example, we conclude that
\[
\begin{aligned}
&V_{ijklm} = V_{[ij][kl]m}\qquad &V_{i[jkl]m}=0\qquad &V_{ij[klm]}=0\\
&Y_{jklm} = Y_{j[kl]m}\qquad  &Y_{[jkl]m} = 0\qquad &Y_{j[klm]} = 0.
\end{aligned}
\]
The differentiated Ricci identity for commuting ambient  
covariant derivatives gives relations involving the conformal
curvature tensors.  
The asymptotic Ricci-flatness of the ambient metric gives relations amongst 
different conformal curvature tensors about which we will be more precise
shortly.   
The covariant derivatives of ambient curvature also satisfy extra
identities involving the infinitesimal dilation $T$ arising from the
homogeneity and straightness conditions, as follows.
\begin{proposition}\label{Tcontract}
The covariant derivatives of the curvature tensor of a straight pre-ambient 
metric $\gt$ satisfy: 
\begin{enumerate}
\item
$T^L \Rt_{IJKL,M_1\cdots M_r}= 
-\sum_{s=1}^r \Rt_{IJKM_s,M_1\cdots \widehat{M_s} \cdots M_r}$
\item
$T^P\Rt_{IJKL,M_1\cdots M_sPM_{s+1}\cdots M_r}= \\
{}\qquad\qquad
-(s+2)\Rt_{IJKL,M_1\cdots M_r} - 
\sum_{t=s+1}^r \Rt_{IJKL,M_1\cdots M_s M_t M_{s+1}\cdots
  \widehat{M_t}\cdots M_r}$.
\end{enumerate}
\end{proposition}
\noindent
Condition (1) in the case $r=0$ is interpreted as the statement
$T^L \Rt_{IJKL}= 0$, or equivalently $\Rt_{IJK0}=0$, mentioned above.  Note
also that the case $s=r$ in (2) 
reduces to \begin{equation}\label{Tend} 
T^P \Rt_{IJKL,M_1\cdots M_rP}
 =-(r+2)\Rt_{IJKL,M_1\cdots M_r}. 
\end{equation}
\begin{proof}
Differentiating the identity $T^K{}_{,I} = \delta^K{}_I$ shows that 
$T^K{}_{,IJ} = 0$.  Commuting the derivatives shows that $T^L\Rt_{IJKL} = 
0$.  Now differentiating this relation successively and using 
again $T^L{}_{,M} = \delta^L{}_M$ proves (1).   

The identity (2) can be proved as a consequence of (1) by commuting the
contracted index to the left using the differentiated Ricci identity and
then applying the second Bianchi identity and (1).  
It can alternately be derived directly as follows.  
Recall the general formula for the Lie derivative of a covariant tensor $U$
in terms of a torsion-free connection:
\begin{equation}\label{covder}
(\mathcal{L}_XU)_{i\cdots j} = X^kU_{i\cdots j,k} 
+  X^k{}_{,i} U_{k\cdots j} + \cdots + X^k{}_{,j} U_{i\cdots k}.   
\end{equation}
{From} the fact that $\mathcal{L}_T\gt = 2\gt$, it follows that 
$\mathcal{L}_T\Rt^{(r)} = 2\Rt^{(r)}$ for all $r$.  
Using this and $T^I{}_{,J}=\delta^I{}_{J}$ in \eqref{covder}, one 
concludes \eqref{Tend}.  Replacing $r$ by $s$ in \eqref{Tend} and 
differentiating $r-s$ more times using again
$T^I{}_{,J}=\delta^I{}_{J}$ gives (2).
\end{proof}

Recall that $T^I=t\delta^I{}_0$.  So a consequence of 
Proposition~\ref{Tcontract} is that $t$ times a component of a 
tensor $\Rt^{(r+1)}$ with a $0$ index can be expressed as a sum 
of components of $\Rt^{(r)}$ with the $0$ index removed and with the
remaining indices permuted.  This implies a corresponding statement
relating tensors $\Rt^{(r+1)}_{\cS_0,\cS_M,\cS_\nf}$ and 
$\Rt^{(r)}_{\cS_0',\cS_M',\cS_\nf'}$ on $M$.

When $n$ is odd, for any choice of $r$ and $\cS_0$, $\cS_M$,
$\cS_\nf$, the tensor
$\Rt^{(r)}_{\cS_0,\cS_M,\cS_\nf}$ depends only on $g$ and each such tensor 
is a natural tensor invariant of $g$.  However, as we have seen in the
examples above, when $n$ is even some of these tensors may depend on the
specific chosen ambient metric and not just on $g$.  The next result gives
precise conditions for this not to happen. 
\begin{proposition}\label{indeterminacy}
Denote by $s_0$, $s_M$, $s_\nf$ the cardinalities of the sets 
$\cS_0$, $\cS_M$, $\cS_\nf$.  Suppose $n$ is even and 
$s_M+2s_\nf\leq n+1$.  Then the tensor $\Rt^{(r)}_{\cS_0,\cS_M,\cS_\nf}$
is independent of the specific ambient metric which has been chosen and is
a natural tensor invariant of $g$. 
\end{proposition}
\begin{proof}
Consider a component $\Rt_{IJKL,M_1\cdots M_r}$ for a 
straight ambient metric $\gt$ in normal form relative to $g$.  Let $s_0$,
$s_M$, $s_\nf$ denote the 
number of $0$'s, indices corresponding to $M$, and $\nf$'s, resp.,
in the list $IJKLM_1\cdots M_r$.  We will prove by induction on $r$ the 
following statement:  all components of $\Rt^{(r)}$ satisfy that 
$\Rt_{IJKL,M_1\cdots M_r} \mod O(\rho^{(n+2-s_M-2s_\nf)/2})$ is independent
of the $O(\rho^{n/2})$ ambiguity of the component $g_{ij}(x,\rho)$ of
$\gt$.  Proposition~\ref{indeterminacy} follows from this upon
restricting to $\rho =0$ and using the fact that the 
Taylor coefficients of $g_{ij}(x,\rho)$ of order $<n/2$ are
natural tensor invariants of $g$. 

We observe first that it suffices to assume that the power of $\rho$  
in our inductive statement satisfies $(n+2-s_M-2s_\nf)/2 \leq n/2 -1$.   
Otherwise $s_M + 2s_\nf \leq 3$, which implies that at most 3 of the
indices $IJKLM_1\cdots M_r$ are not equal to $0$.  Such a component
$\Rt_{IJKL,M_1\cdots M_r}$ vanishes identically by
Proposition~\ref{Tcontract}. 

Now proceed with the induction.  The case $r=0$ follows easily from
\eqref{curv} and $\Rt_{IJK0}=0$.  For the inductive step, 
consider $\Rt_{IJKL,M_1\cdots M_r P}$.  If $P=0$, the fact that the
ambiguity of this component is no worse than that of 
$\Rt_{IJKL,M_1\cdots M_r}$ follows immediately from \eqref{Tend}.  
For $P\neq 0$, write
$$
\Rt_{IJKL,M_1\cdots M_r P}=\pa_P\Rt_{IJKL,M_1\cdots M_r } 
-\Gat_{IP}^A\Rt_{AJKL,M_1\cdots M_r } -\ldots
 -\Gat_{M_rP}^A\Rt_{IJKL,M_1\cdots A } . 
$$
If $P\neq \nf$, then the ambiguity in the first term on the right hand side
vanishes to the same
order as that in $\Rt_{IJKL,M_1\cdots M_r}$, while if $P=\nf$, then it
vanishes to one order less.  Thus the required vanishing 
for the first term follows from the inductive hypothesis.   
Consider next the second term.  The Christoffel symbols are given by
\eqref{cnr}.  They have their own ambiguity of at most $O(\rho^{(n/2-1)})$ 
owing to the $O(\rho^{n/2})$ ambiguity in $g_{ij}(x,\rho)$.  But since by
the observation above we can assume that $(n+2-s_M-2s_\nf)/2 \leq n/2 -1$,   
we can neglect the ambiguity in the Christoffel symbols.  Consider the 
ambiguity in the second term arising from the ambiguity in
$\Rt_{AJKL,M_1\cdots M_r }$.  If $P=\nf$, then we need to show  
that the order of vanishing of this
ambiguity is at most 1 less than that of $\Rt_{IJKL,M_1\cdots M_r }$.  
But this is clear from the inductive hypothesis without even considering
the Christoffel symbol since these two components of
$\Rt^{(r)}$ have at most one different index and the change in ambiguity
from changing one index is at most 1.  If $P=p$ is between $1$ and $n$, we
need to show that the order of vanishing of the ambiguity in
$\Gat_{Ip}^A\Rt_{AJKL,M_1\cdots M_r }$ is at most $\frac12$ less than that
of $\Rt_{IJKL,M_1\cdots M_r }$.  This is clear by the same reasoning as in 
the case $P=\nf$ unless $I=0$ and $A=\nf$.  But \eqref{cnr} shows that 
$\Gat_{0p}^\nf = 0$, so that the inductive statement holds in this case
also.  The same reasoning as for the second term applies to the
remaining terms on the right hand side above. 
\end{proof}

The weighting which appears in Proposition~\ref{indeterminacy} suggests the
following definition.
\begin{definition}\label{strength}
We define the {\it strength} of lists of indices in $\R^{n+2}$ as follows.
Set $\|0\| = 0$,  
$\|i\| = 1$ for $1\leq i \leq n$, and $\|\infty\|=2$.  For a list,  
write $\|I\ldots J\| = \|I\|+\cdots +\|J\|$.
\end{definition}
\noindent
Proposition~\ref{indeterminacy} thus asserts that for $n$ even, the
curvature component $\Rt_{IJKL,M_1\cdots M_r}$ of a straight ambient metric
in normal form relative to $g$ is well-defined at $\rho =
0$ independent of the choice of ambient metric so long as $\|IJKLM_1\cdots  
M_r\|\leq n+1$.      

Next we consider the trace-free condition imposed by the asymptotic
Ricci-flatness of ambient metrics.
\begin{proposition}\label{tracefree}
If $n$ is odd, the covariant derivatives of curvature of a straight
ambient metric in normal form relative to $g$ satisfy at $\rho =0$:
\begin{equation}\label{trace}
\gt^{IK}\Rt_{IJKL,M_1\cdots M_r} = 0.
\end{equation}
If $n$ is even, the same result holds assuming that  
$\|JLM_1\cdots M_r\|\leq n-1$.  
\end{proposition}
\begin{proof}
The result is clear for $n$ odd since the Ricci curvature of $\gt$
vanishes to infinite order.  For $n$ even we prove by induction on $r$ the
statement that 
$$
\gt^{IK}\Rt_{IJKL,M_1\cdots M_r} = O(\rho^{(n-\|JLM_1\cdots M_r\|)/2})
$$
for all components $JLM_1\cdots M_r$.  The desired result then follows upon
setting $\rho = 0$.  

The case $r=0$ of the induction is a consequence of the fact that  
$\Rt_{JL}=O^+_{JL}(\rho^{n/2-1})$.  For the induction, write
\[
\begin{split}
\gt^{IK}\Rt_{IJKL,M_1\cdots M_r} = (&\gt^{IK}\Rt_{IJKL,M_1\cdots
  M_{r-1}}),_{M_r}
= \pa_{M_r}(\gt^{IK}\Rt_{IJKL,M_1\cdots M_{r-1}})\\
&-\Gat^A_{JM_r}\gt^{IK}\Rt_{IAKL,M_1\cdots M_{r-1}}
-\ldots 
-\Gat^A_{M_{r-1}M_r}\gt^{IK}\Rt_{IJKL,M_1\cdots A}.
\end{split}
\]
The bound on the first term on the right hand side follows from the
induction hypothesis and the effect of differentiation on the order of
vanishing.  The bound on the remaining terms is easily seen to be a
consequence of the induction hypothesis and the fact that $\Gat^I_{JK}=0$ 
for $\|I\|>\|J\|+\|K\|$, which follows from \eqref{cnr}.
\end{proof}
\noindent
When the trace in \eqref{trace} is written out in terms of components,
the resulting identity can be interpreted as expressing a trace with
respect to $g$ of 
a conformal curvature tensor $\Rt^{(r)}_{\cS_0,\cS_M,\cS_\nf}$ in terms of
other conformal curvature tensors. 

Suppose now we choose a conformally related metric $\wh{g} = e^{2\Up}g$.
By Proposition~\ref{normalform}, an ambient metric in normal form for $g$
can be put into normal form for $\wh{g}$ by a unique homogeneous
diffeomorphism which restricts to the identity on $\cG$.  By calculating
on $\cG$ the Jacobian of this diffeomorphism and using the fact that the 
ambient curvature tensors are tensors, we will be able to compute the 
transformation laws of the conformal curvature tensors under conformal
change.  We denote the coordinates relative to $\wh{g}$ by   
$(\wh{t}$, $\wh{x}$, $\wh{\rho})$ and the conformal curvature tensors for 
$\wh{g}$ by 
$\wh{\Rt}{}^{(r)}_{\cS_0,\cS_M,\cS_\nf}=\wh{\Rt}_{IJKL,M_1\cdots  
M_r}|_{\wh{\rho} = 0,\,\wh{t}=1}$.
\begin{proposition}\label{ctransform}
Let $g$ and $\wh{g}=e^{2\Up}g$ be conformally related metrics on $M$.  
Let $IJKLM_1\cdots M_r$ be a list of indices, $s_0$ of which are $0$,
$s_M$ of which correspond to $M$, and $s_\nf$ of which are $\nf$.  
If $n$ is even,
assume that $s_M + 2s_\nf \leq n+1$.  Then the conformal curvature tensors
satisfy the conformal transformation law:
\begin{equation}\label{transformula}
\wh{\Rt}_{IJKL,M_1\cdots M_r}|_{\wh{\rho} = 0,\,\wh{t}=1}
= e^{2(1-s_\nf)\Up} \Rt_{ABCD,F_1\cdots F_r}|_{\rho = 0,\,t=1}p^A{}_I
\cdots p^{F_r}{}_{M_r},
\end{equation}
where $p^A{}_I$ is the matrix
\begin{equation}\label{pmatrix}
p^A{}_I = 
\left(
\begin{matrix}
1&\Up_i&-\frac12 \Up_k\Up^k\\
0&\delta^a{}_i&-\Up^a\\
0&0&1
\end{matrix}
\right).
\end{equation}
\end{proposition}

We make several observations and explanations before giving the proof.  
For each division of $IJKLM_1\cdots M_r$ into 
subsets $\cS_0$, $\cS_M$, and $\cS_\nf$, the identity \eqref{transformula}  
is a relation amongst tensors on $M$.  If $\Up$ is constant, then 
$p^A{}_I = \delta^A{}_I$ and \eqref{transformula} just tells how 
$\Rt^{(r)}_{\cS_0,\cS_M,\cS_\nf}$ scales.   Because of the
upper-triangular form of the matrix $p^A{}_I$, in the general case the
other 
terms on the right hand side all involve ``earlier'' conformal curvature
tensors in the sense that each `$i$' can be replaced only by $0$ and each 
$\nf$ only by
an `$i$' or a $0$.  The conformal transformation law of a conformal 
curvature tensor involves only other conformal curvature tensors and first
derivatives of $\Up$.  

Consider the case $r=0$.  If we take $s_0$, $s_\nf =0$, then  
$\wh{\Rt}_{ijkl}|_{\wh{\rho} = 0,\,\wh{t}=1} =\wh{W}_{ijkl}$ is the
Weyl tensor of $\wh{g}$, and since $\Rt_{IJK0}=0$, it follows that
\eqref{transformula} 
reproduces the conformal invariance $\wh{W}_{ijkl} = e^{2\Up}W_{ijkl}$
of the Weyl tensor.  Taking $I =\nf$ reproduces the transformation law
$\wh{C}_{jkl} = C_{jkl} -\Up^iW_{ijkl}$ of the Cotton tensor, and taking 
$I$, $L = \nf$ gives for $n\neq 4$ the transformation law
$$
\wh{B}_{jk} = 
e^{-2\Up}\left( B_{jk}+(n-4)\left[\Up^l(C_{jkl}+C_{kjl})
-\Up^i\Up^l W_{ijkl}\right]\right)
$$
of the Bach tensor.  For $n=4$ we have already noted the
conformal invariance of the Bach tensor in its appearance as the
obstruction tensor.

The conformal transformation law \eqref{transformula} is 
one of the most important features of the conformal curvature tensors
and makes evident the importance of their interpretation as components 
of tensors on the ambient space.  Since the conformal curvature tensors
clearly vanish if $g$ is flat, one consequence of \eqref{transformula}  
is that they also vanish if $g$ is locally conformally flat.  
This is not {\it a priori} obvious since they are not individually
conformally invariant.  We will remove the restriction 
$s_M+2s_\infty\leq n+1$ from this observation in even dimensions in
Chapter~\ref{flat}.     
Collections of tensors  
on $M$ which transform conformally according to rules of the form  
\eqref{transformula} define sections of weighted tensor powers of the
cotractor bundle  
associated to the conformal structure.  See \cite{BEGo}, \cite{Go1} for
development of this point of view, and \cite{CG}, \cite{BrG}, \cite{AL} for
discussion of the relation between tractors and the ambient construction.    
Note that for each $x\in M$, the matrix $p^A{}_I(x)$ 
is in the orthogonal group of the quadratic form
$$
\left(
\begin{matrix}
0&0&1\\
0&g_{ij}(x)&0\\
1&0&0
\end{matrix}
\right).
$$
\noindent
{\it Proof of Proposition~\ref{ctransform}.}
Let $\gt$ be a straight ambient metric in normal form relative to $g$.  
The metric $\gt$ is defined on a neighborhood of $\cG\times \{0\}$ in
$\cG\times \R$.  If we use the trivialization induced by $g$ to identify 
$\cG$ with $\R_+\times M$, then $\gt$ takes the form \eqref{spform}.  
According to Proposition~\ref{normalform}, we can put
$\gt$ in normal form relative to $\wh{g}$ by a unique homogeneous
diffeomorphism $\phi$ defined on a neighborhood of $\cG\times \{0\}$ which
restricts to the identity on $\cG$.  If on the
$\wh{g}$ side we identify $\cG$ with $\R_+\times M$ using the
trivialization induced by $\wh{g}$,  
then $\phi^*\gt$ will take the form \eqref{spform} in coordinates  
$(\wh{t},\wh{x},\wh{\rho})$.  The trivializations are related by
\eqref{triv}.  Therefore we deduce the existence of a homogeneous
diffeomorphism $\psi(\wh{t},\wh{x},\wh{\rho}) = (t,x,\rho)$ 
on a homogeneous neighborhood of $\R_+\times M\times \{0\}$ into 
$\R_+\times M \times \R$ with the properties that $\wh{\gt}:= \psi^*\gt$ 
takes the form \eqref{spform} in $(\wh{t},\wh{x},\wh{\rho})$ and  
$$
\psi(\wh{t},\wh{x},0) = (\wh{t}e^{\Up},\wh{x},0). 
$$

Now differentiation of this relation determines the derivatives of the
components of $\psi$ at 
$\wh{\rho} = 0$ in the $\wh{t}$ and $\wh{x}$ directions.
The derivatives in the $\wh{\rho}$ direction are determined from these and
from the requirement that 
$\psi^*\gt = 2\wh{t}d\wh{t}d\wh{\rho}
+\wh{t}^2\wh{g}_{ij}d\wh{x}^id\wh{x}^j$ at $\wh{\rho}=0$.  This gives 
$$
(\psi')^A{}_I|_{\wh{\rho}=0} = 
\left(
\begin{matrix}
e^\Up&\wh{t}e^\Up \Up_i&-\frac12 \wh{t}e^{-\Up}\Up_k\Up^k\\
0&\delta^a{}_i&-e^{-2\Up}\Up^a\\
0&0&e^{-2\Up}
\end{matrix}
\right).
$$
Observe that this may be factored as 
\begin{equation}\label{psid}
\psi'|_{\wh{\rho}=0}=d_1pd_2,
\end{equation} 
where $p$ is given by \eqref{pmatrix} and 
\begin{equation}\label{diag}
d_1 = 
\left(
\begin{matrix}
\wh{t}e^\Up&0&0\\
0&Id&0\\
0&0&1
\end{matrix}
\right)\qquad\
d_2=
\left(
\begin{matrix}
\wh{t}^{-1}&0&0\\
0&Id&0\\
0&0&e^{-2\Up}
\end{matrix}
\right).
\end{equation}

Each curvature component $\Rt_{IJKL,M_1\cdots M_r}$ is homogeneous
with respect to the dilations.  Since $\delta^*_s \gt = s^2 \gt$, it
follows that for the full tensors we have $\delta^*_s \Rt^{(r)}= s^2
\Rt^{(r)}$.  Since $\pa_t$ is homogeneous of degree $-1$ and $\pa_\rho$ and
$\pa_{x^i}$ are of degree $0$, we deduce that the component 
$\Rt_{IJKL,M_1\cdots M_r}$ is homogeneous of degree $2-s_0$.  Therefore
\begin{equation}\label{homo}
\wh{\Rt}_{IJKL,M_1\cdots M_r}|_{\wh{\rho} = 0,\,\wh{t}=1}=
e^{(2-s_0)\Up}\wh{\Rt}_{IJKL,M_1\cdots M_r}|_{\wh{\rho} = 0,\,
\wh{t}=e^{-\Up}}.  
\end{equation}

Since the covariant derivatives of curvature are tensorial, we have   
$$
\wh{\Rt}_{IJKL,M_1\cdots M_r} = 
\Rt_{ABCD,F_1\cdots F_r}\circ \psi\,\,(\psi')^A{}_I
\cdots (\psi')^{F_r}{}_{M_r}.
$$
Evaluate both sides at $\wh{\rho} = 0$, $\wh{t}=e^{-\Up}$ and use
\eqref{psid} to obtain
$$
\wh{\Rt}_{IJKL,M_1\cdots M_r}|_{\wh{\rho} = 0,\,\wh{t}=e^{-\Up}}
= e^{(s_0-2s_\nf)\Up} \Rt_{ABCD,F_1\cdots F_r}|_{\rho = 0,\,t=1}p^A{}_I 
\cdots p^{F_r}{}_{M_r}.
$$
When combined with \eqref{homo}, this gives \eqref{transformula}.

This argument is valid for all components of ambient curvature whether $n$
is even or odd.  However, when $n$ is even, the components with
$s_M+2s_\nf>n+1$ generally depend on the ambient metric which has been
chosen.  To 
ensure that the component depends only on $g$ as formulated in
Proposition~\ref{ctransform}, we take $s_M+2s_\nf\leq n+1$ 
and use Proposition~\ref{indeterminacy}. 
\stopthm

If $n$ is even and $\|IJKLM_1\cdots M_r\|\geq n+2$, then the restriction
of $\Rt_{IJKL,M_1\cdots M_r}$ to $\rho =0$ in general does depend on the 
choice of ambient metric.  The component $\Rt_{\infty 
  ij\infty,\,\underbrace{\scriptstyle{\infty\ldots\infty}}_{n/2-2}}$ is an 
important example.  The next proposition, which we will use in 
Chapter~\ref{flat}, makes explicit this dependence.
\begin{proposition}\label{ambiguity}
Let $n\geq 4$ be even.  There is a natural trace-free symmetric 2-tensor
$K_{ij}$ depending on a pseudo-Riemannian metric $g$, with the 
following properties.  
Let $\gt_{IJ}$ be an ambient metric in normal form relative to $g$, and
write $\gt_{ij}=t^2 g_{ij}(x,\rho)$.  Then $K_{ij}$ can be    
expressed algebraically in terms of the tensors $g^{ij}|_{\rho =0}$ and 
$\pa_\rho^m g_{ij}|_{\rho =0},$ $0\leq m\leq n/2-1$, and  one has 
\begin{equation}\label{Kformula}
2\Rt_{\infty 
  ij\infty,\,\underbrace{\scriptstyle{\infty\ldots\infty}}_{n/2-2}}|_{\rho
  =0,\,t=1}  
=\tf\left(\pa_\rho^{n/2}g_{ij}|_{\rho=0}\right) + K_{ij}.   
\end{equation}
\end{proposition}
\begin{proof}
First suppose that $\gt$ is straight.  
An easy induction beginning with the last line of \eqref{curv} and using 
\eqref{cnr} shows that for $m\geq 0$, 
\begin{equation}\label{rhoderivs}
2t^{-2}\Rt_{\infty 
  ij\infty,\,\underbrace{\scriptstyle{\infty\ldots\infty}}_{m}}
=\pa_\rho^{m+2}g_{ij} + \cP^{(m)}_{ij}\left(g^{kl},g_{kl},  
\pa_\rho g_{kl},\ldots, \pa_\rho^{m+1}g_{kl}\right),
\end{equation}
where $\cP^{(m)}_{ij}$ is a tensor depending polynomially on the indicated 
arguments.  Define 
$$
K_{ij}=\tf\left(\cP^{(n/2-2)}_{ij}|_{\rho =0}\right).  
$$
Using \eqref{inverse}, we have
$$
t^{-2}g^{ij}\Rt_{\infty ij
  \infty,\,\underbrace{\scriptstyle{\infty\ldots\infty}}_{m}} 
= \gt^{IJ}\Rt_{\infty IJ
  \infty,\,\underbrace{\scriptstyle{\infty\ldots\infty}}_{m}} 
=- \Rt_{\infty   
  \infty,\,\underbrace{\scriptstyle{\infty\ldots\infty}}_{m}}.
$$
Since $\Rt_{IJ}=O(\rho^{n/2-1})$, it follows that 
$\Rt_{\infty ij
  \infty,\,\underbrace{\scriptstyle{\infty\ldots\infty}}_{m}}|_{\rho 
    =0, t=1}$ is trace-free for $m\leq n/2-2$.  
So restricting 
\eqref{rhoderivs} with 
  $m=n/2-2$ to $\rho =0$, $t=1$ and taking the trace-free part gives the
  desired  conclusion.  

It is not hard to check using \eqref{inverse} and \eqref{c} that the last
line of \eqref{curv} and \eqref{rhoderivs} are also valid even if $\gt$ is
not assumed to be straight, so that the same proof is valid in the general 
case.  The $\gt_{00}$ and $\gt_{0i}$ components simply do 
not enter into any of the expressions which occur.
\end{proof}

\noindent
Evaluating the last line of \eqref{curv} using \eqref{initial} shows that 
$K_{ij}=-2\tf\left(P_i{}^kP_{jk}\right)$ for $n=4$.     

The discussion in the rest of this chapter concerns the dimension
dependence of certain 2-tensors.  Consider partial contractions $C$ of  
the form 
$$
\pcontr\left(\nabla^{r_1}R\otimes \cdots
\otimes\nabla^{r_L}R\right)\quad \text{or}\quad
\pcontr\left(\nabla^{r_1}R\otimes \cdots
\otimes\nabla^{r_L}R\otimes g\right)
$$
of the covariant derivatives of the curvature tensor of a metric
$g$ and possibly also $g$ itself, in which 2 indices remain uncontracted.   
For partial contractions of the second type, we require that the
uncontracted indices are the indices on $g$.  
Either type of partial contraction may be viewed formally as  
a choice of non-negative integers $r_1, \cdots, r_L$, a pairing of the 
contracted indices, and an  
ordering of the uncontracted indices.  Consider also formal linear 
combinations $\xi=\sum_{i=1}^M f_i(d) C_i$ of such formal partial
contractions whose coefficients $f_i(d)$ are rational  
functions of a single variable $d$.  Any such formal linear combination
$\xi$ 
defines a natural 2-tensor $\operatorname{Eval}_{d=n}\xi$ in any dimension 
$n\geq 3$ which is a regular point for   
all the coefficients $f_i$ by evaluation on metrics in dimension $n$.  One 
can also define the residue of $\xi$ at any dimension $n\geq 3$ to be the
natural 2-tensor in dimension $n$ given by  
$$
\operatorname{Res}_{d=n}\xi = \sum_{i=1}^M (\operatorname{Res}_{d=n}f_i)
\operatorname{Eval}_{d=n}C_i.  
$$
By considering the product of a metric in a given dimension with flat 
metrics, it can be seen that if $\xi$ satisfies that for some $n_0$, 
$\operatorname{Eval}_{d=n}\xi =0$ for all $n\geq n_0$, then 
$\operatorname{Res}_{d=n}\xi =0$ for all $n$.  Thus the residue at any
dimension is independent of the  
way a given dimension-dependent family of natural tensors is written
formally.  

The construction of the ambient metric in Chapter~\ref{formalsect} shows
that for  
each $m\geq 1$, there is a formal linear combination $\xi_m$ as above such
that if $n$ is odd or if $n$ is even and $m< n/2$, then
for a straight ambient metric in dimension $n$ in normal form relative to
$g$, one has
$\pa_\rho^mg_{ij}|_{\rho=0}=\operatorname{Eval}_{d=n}\xi_m$.
For $m\geq 0$, the curvature component $\Rt_{\infty   
ij\infty,\,\underbrace{\scriptstyle{\infty\ldots\infty}}_m}|_{\rho=0,
t=1}$ can be so written as well for $n$ odd or for $n$ even and 
$m< n/2-2$.       
Fix an even integer $n\geq 4$ and consider 
$\operatorname{Res}_{d=n} \Rt_{\infty 
ij\infty,\,\underbrace{\scriptstyle{\infty\ldots\infty}}_{n/2-2}}|_{\rho=0, 
t=1}$.  As noted above, the residue is independent of the choice of formal
expression, so this notation is justified.     
For $n=4$ we have $\Rt_{\infty ij\infty}|_{\rho =0,   
t=1}=-(d-4)^{-1}B_{ij}$ from \eqref{curv0}.  Thus
$\Rt_{\infty ij\infty}|_{\rho =0, t=1}$ has a  
simple pole at $d=4$ with residue the negative of the Bach tensor.  
The following proposition generalizes this fact to higher $n$.    
\begin{proposition}\label{bachgen}
Let $n\geq 4$ be an even integer.  
$\Rt_{\infty 
ij\infty,\,\underbrace{\scriptstyle{\infty\ldots\infty}}_{n/2-2}}|_{\rho=0, 
t=1}$ has a simple pole at $d=n$ with residue given by
$$
\operatorname{Res}_{d=n}
\Rt_{\infty 
ij\infty,\,\underbrace{\scriptstyle{\infty\ldots\infty}}_{n/2-2}}|_{\rho=0, 
t=1}
=(-1)^{n/2-1}\left[2^{n/2-2}(n/2-2)!\right]^{-1}\cO_{ij}.
$$
\end{proposition}
\begin{proof}
Write $\gt_{ij}=t^2g_{ij}(x,\rho)$ as in Proposition~\ref{ambiguity}.  
The derivatives $\pa_\rho^m g_{ij}|_{\rho=0}$ are all regular in $d$ at 
$d=n$ for $0\leq m\leq n/2-1$, so the argument of
Proposition~\ref{ambiguity} shows that 
\begin{equation}\label{residue}
\operatorname{Res}_{d=n}
\Rt_{\infty 
ij\infty,\,\underbrace{\scriptstyle{\infty\ldots\infty}}_{n/2-2}}|_{\rho=0, 
t=1}=\tfrac12
\operatorname{Res}_{d=n}\tf\left(\pa_\rho^{n/2}g_{ij}|_{\rho=0}\right).  
\end{equation}
For $d>n$, $\tf\left(\pa_\rho^{n/2}g_{ij}|_{\rho=0}\right)$ is determined
by replacing $n$ 
by $d$ in the first line of \eqref{evol}, setting the right hand side equal
to 0, applying $\pa_{\rho}^{n/2 -1}|_{\rho = 0}$, and taking the trace-free 
part.  At $\rho =0$ one has
$$
\pa_\rho^{n/2-1}\left[\rho g_{ij}'' -(d/2-1)g_{ij}'\right]
=\tfrac12 (n-d)\pa_\rho^{n/2}g_{ij}.
$$
For the remaining terms in the right hand side of the first line of
\eqref{evol}, set 
$$
T_{ij}=\tf\left[\pa_\rho^{n/2-1}\left(
 -\rho g^{kl}g_{ik}'g_{jl}' +\tfrac12 \rho
g^{kl}g_{kl}'g_{ij}' 
- \tfrac12 g^{kl}g_{kl}'g_{ij} +R_{ij}\right)
|_{\rho=0}\right].
$$
Then $T_{ij}$ is expressible in terms of the 
$\pa_\rho^{m}g_{ij}|_{\rho =0}$ for $m\leq n/2-1$, so is
regular at $d=n$.  Differentiation of \eqref{evol} as indicated above thus
gives 
$
\tfrac12 (n-d)\tf\left(\pa_\rho^{n/2}g_{ij}|_{\rho =0}\right)+T_{ij}=0,
$
so 
\begin{equation}\label{T}
\tfrac12\operatorname{Res}_{d=n}
\tf\left(\pa_\rho^{n/2}g_{ij}|_{\rho=0}\right)  
=T_{ij}|_{d=n}.
\end{equation}

On the other hand, the obstruction tensor $\cO_{ij}$ in dimension $n$ is
obtained by setting 
$\Rt_{ij} = c_n^{-1}(2\rho)^{n/2-1} \cO_{ij} \mod O(\rho^{n/2})$ in the
first line of \eqref{evol} (without changing $n$ to $d$), applying
$\pa_{\rho}^{n/2 -1}|_{\rho = 0}$, and taking the trace-free part.   This 
gives 
$$
c_n^{-1}2^{n/2-1}(n/2-1)!\,\cO_{ij}=T_{ij}|_{d=n}.  
$$
Combining with \eqref{T} and \eqref{residue} gives the result.
\end{proof}

Thus $(d-n)\Rt_{\infty 
ij\infty,\,\underbrace{\scriptstyle{\infty\ldots\infty}}_{n/2-2}}|_{\rho=0, 
t=1}$ is a natural trace-free symmetric 2-tensor depending on $d$ regular
also at $d=n$ which equals a multiple of 
the obstruction tensor for $d=n$.  Of course, there are many other such
tensors; one can add to  
$\Rt_{\infty 
ij\infty,\,\underbrace{\scriptstyle{\infty\ldots\infty}}_{n/2-2}}|_{\rho=0, 
t=1}$ any tensor regular at $d=n$.  For example, one can replace 
$\Rt_{\infty 
ij\infty,\,\underbrace{\scriptstyle{\infty\ldots\infty}}_{n/2-2}}|_{\rho=0, 
t=1}$ by
$\frac12\tf\left(\pa_\rho^{n/2}g_{ij}|_{\rho=0}\right)$; 
the relation \eqref{Kformula} holds also as dimension-dependent natural
tensors in dimensions $d$ which are odd or which are even and $>n$, with 
$K_{ij}$ regular also at $d=n$.  A main    
advantage of the choice $\Rt_{\infty     
ij\infty,\,\underbrace{\scriptstyle{\infty\ldots\infty}}_{n/2-2}}|_{\rho=0, 
t=1}$ is that its tranformation law under conformal change is known and
relatively simple, being given by \eqref{transformula} for $d$ odd or $d>n$
even.   

One can also consider the continuation of 
$\Rt_{\infty 
ij\infty,\,\underbrace{\scriptstyle{\infty\ldots\infty}}_{n/2-2}}|_{\rho=0, 
t=1}$ to even values of $d<n$.  In general, the pole is not
simple and its order increases with $n-d$.    

\section{Conformally Flat and Conformally Einstein Spaces}\label{flat}

If $n$ is odd, an ambient metric in normal form is uniquely 
determined to infinite order by $(M,g)$.  Theorem~\ref{odd} shows that  
within the family of all formal solutions, this one is 
distinguished by the vanishing of the $\rho^{n/2}$ coefficient,   
or equivalently by the condition that it be smooth.  In the 
Poincar\'e realization, the corresponding condition is that the expansion
have no $r^n$ term, or equivalently that it be
even.  The fact that there is such a normalization giving rise to a unique
diffeomorphism class of ambient metrics associated to the conformal
manifold $(M,[g])$ is crucial for the applications to  
jet isomorphism and invariant theory in Chapters~\ref{jet} and \ref{sinv}. 

When $n$ is even, the situation is different primarily because of the 
existence of the obstruction tensor, which gives rise to the log terms in 
the infinite order expansions determined in Theorem~\ref{even}.  But even 
for conformal classes with 
vanishing obstruction tensor, for which all formal solutions are 
smooth, the coefficient of $\rho^{n/2}$ is undetermined, so there is not a 
unique solution.
One can impose the condition that $\gt$ be straight;
this removes some of the freedom, but there is no apparent analog of the 
the odd-dimensional normalization.

In this chapter we discuss two families of special conformal classes for
which there is a natural conformally invariant normalization in even 
dimensions giving rise to a unique diffeomorphism class of infinite order    
formal expansions for an ambient (or Poincar\'e) metric.  These are the 
locally conformally flat conformal classes and the conformal classes
containing an Einstein metric.  

Consider first the case of a locally conformally flat manifold $(M,[g])$.  
This means that in a neighborhood of any point of $M$, there is a flat
metric, say $\ga$, in the conformal class.  Now \eqref{curv} implies that  
the metric $2\rho dt^2 + 2tdtd\rho +t^2 \ga$, with $\ga$ independent of
$\rho$, is flat.  If $n$ is odd, this implies that every ambient metric is
flat to infinite order, since the ambient metric is unique to infinite
order up to diffeomorphism and the condition that $\gt$ is 
flat is invariant under diffeomorphism.  If $n\geq 4$ is even, the 
condition that $\gt$ be flat to infinite order uniquely determines $\gt$ to
infinite order up to diffeomorphism.  This follows from Theorem~\ref{even}
and Proposition~\ref{ambiguity}:  the obstruction tensor vanishes and
Theorem~\ref{even} shows that the full expansion is determined by  
$\tf\left(\pa_\rho^{n/2}g_{ij}|_{\rho =0}\right)$, which by
Proposition~\ref{ambiguity}  must equal $-K_{ij}$ if $\gt$ is flat.  
Thus we have:
\begin{proposition}\label{confflatprop}
Let $n\geq 3$ and suppose that $(M,[g])$ is locally conformally flat.  Then
there 
exists an ambient metric $\gt$ for $(M,[g])$ which is flat to infinite
order, and such a $\gt$ is unique to infinite order up to diffeomorphism.  
\end{proposition}

\noindent
We remark that it follows from the fact that 
the metric $2\rho dt^2 + 2tdtd\rho +t^2 \ga$ is straight, that $\gt$ in
Proposition~\ref{confflatprop} can be taken to be straight.  
In Proposition~\ref{ss} and Theorem~\ref{flatnbhd},
we will extend Proposition~\ref{confflatprop}  
to establish the existence and uniqueness of a flat ambient metric not just
to infinite order, but in a neighborhood of $\cG\times \{0\}$.   

One consequence of Proposition~\ref{confflatprop} is that the conformal
curvature tensors can be invariantly defined for all orders of
differentiation for locally conformally flat metrics in even dimensions. 
Because of conformal invariance, Proposition~\ref{ctransform} continues to
hold, so all the conformal curvature tensors are defined and vanish for
locally conformally flat metrics.       

When $n=2$ the situation is quite different.  (Recall that when $n=2$,
every metric is locally conformally flat.)  In fact, when $n=2$, every
straight ambient metric is flat to infinite order.  This follows by
consideration of the Poincar\'e metric.  If $\gt$ is a straight ambient
metric, then $\Ric(\gt)$ vanishes to infinite order, so 
Proposition~\ref{curvrelation} shows that the Poincar\'e metric $g_+$ 
defined in Proposition~\ref{ambienttopoincare} satisfies 
that $\Ric(g_+)+ng_+$ vanishes to infinite order.  Since the Ricci
tensor determines the curvature tensor in dimension 3, it follows that 
$\operatorname{Rm}(g_+)+g_+\owedge g_+$ vanishes to infinite order, so by 
Proposition~\ref{curvrelation} again, $\gt$ is flat to infinite order.   
There are many diffeomorphism-inequivalent straight ambient metrics for
$(M,[g])$; they are parametrized in Theorem~\ref{2dim}.    

An explicit formula for flat ambient metrics in normal form relative to   
an arbitrary metric in the conformal class 
was given in \cite{SS} (in the Poincar\'e realization).  
This is closely related to a result of Epstein \cite{PP} describing 
the form of hyperbolic metrics near conformal infinity in terms of data on 
an interior hypersurface.  
In the following, we will denote by $\ga$ an arbitrary metric  
in the conformal class on $M$, and by $g=g_\rho$ a metric obtained by
fixing $\rho$ in the $ij$ component of an ambient metric.  The formula of
\cite{SS} is as follows.
\begin{proposition}\label{ss}
Let $\ga_{ij}$ be a locally conformally flat metric on a manifold $M$ of 
dimension $n$ and let $P_{ij}$ be a symmetric 2-tensor on $M$.  
Set 
\begin{equation}\label{gformula}
(g_\rho)_{ij}=\ga_{ij} +2 P_{ij}\rho + P_{ik}P^k{}_j \rho^2 
\end{equation}
and define      
a metric $\gt$ on $\R_+\times U$ , where $U$ is the subset of $M\times \R$
on which $g_\rho$ is nondegenerate, by   
$$
\gt = 2t dtd\rho +2\rho dt^2 + t^2g_\rho.   
$$
Then $\gt$ is flat if: 
\begin{itemize}
\item
$n\geq 3$ and $P$ is the Schouten tensor of $\ga$.
\item
$n=2$ and $P$ is any symmetric 2-tensor satisfying
$$
P_i{}^i=\tfrac12 R \qquad\text{and}\qquad P_{ij,}{}^j= \tfrac12 R_{,i}. 
$$
\end{itemize}
Here indices are raised and lowered using $\ga$, 
indices preceded by a comma denote covariant derivatives with respect 
to the Levi-Civita connection of $\ga$, and $R$ denotes the scalar
curvature of $\ga$.  
\end{proposition}

Note that the local conformal flatness of $\ga$ implies for $n\geq 3$ that
the Schouten tensor of $\ga$ satisfies  
$P_{ij,k}=P_{(ij,k)}$.  For $n=2$, the hypotheses on $P$ imply that 
$P_{i[j,k]}$ is trace-free.  But in two dimensions, a trace-free tensor
$A_{ijk}$ satisfying $A_{ijk}=A_{i[jk]}$ and $A_{[ijk]}=0$ must vanish
identically (see \cite{W}).  So in all dimensions, a consequence of the
hypotheses of Proposition~\ref{ss} is that $P_{ij,k}=P_{(ij,k)}$.   

As described in Propositions~\ref{ambienttopoincare} and
\ref{curvrelation},  
under the change of variable $\rho =-\frac12 r^2$, the condition that $\gt$
is flat is equivalent to the condition that $g_+$ is hyperbolic (i.e. has 
constant sectional curvature $-1$), where
$$
g_+ = r^{-2}\left(dr^2+g_{-\tfrac12 r^2}\right).
$$
Thus one obtains explicit formulae for hyperbolic metrics.   
Consequences for the coefficients in the expansion of
the volume form of $g_+$ in the case of a 
locally conformally flat infinity are given in \cite{GJ}.  

We will indicate three proofs of Proposition~\ref{ss}.  First we prove it
by direct calculation.   

Since there are two metrics involved, it is important to be careful about  
conventions for raising and lowering indices.  We will use $\ga$ to raise
and lower indices with the one exception that $g^{ij}$ denotes the inverse
of $g_{ij}$.  Hereafter we write $g$ for $g_\rho$; the $\rho$ dependence is
to be understood.  We can write  
$$
g_{ij}=\ga_{il}U^l{}_kU^k{}_j=U_{ik}U^k{}_j 
$$
where
\begin{equation}\label{Udef}
U^i{}_j=\delta^i{}_j +\rho P^i{}_j.
\end{equation}
Note that $g$ is nondegenerate precisely where $U$ is nonsingular.  
Let $V=U^{-1}$ so that $V^i{}_kU^k{}_j = \delta^i{}_j$ and observe that 
$U_{ij}=U_{ji}$, $V_{ij}=V_{ji}$.  It is evident that 
\begin{equation}\label{Vg}
V^k{}_ig_{kj}=U_{ij}
\end{equation}
and
\begin{equation}\label{gderiv}
g'_{ij}=2P_{ik}U^k{}_j.
\end{equation}

First we derive the relation between the Levi-Civita connections   
${}^g\nabla$ and ${}^\ga\nabla$ and the curvature tensors 
${}^gR_{ijkl}$ and ${}^\ga R_{ijkl}$ of $g$ and $\ga$. 
\begin{lemma}\label{diffcurv}
Let $g$ be given by \eqref{gformula}, where $\ga$ is a metric, $P_{ij}$ is
a symmetric 2-tensor satisfying $P_{ij,k}=P_{(ij,k)}$ and $\rho\in \R$.  If
$g$ is nondegenerate, then the Levi-Civita connections are related by 
\begin{equation}\label{difften}
{}^g\nabla_i\eta_j = {}^\ga\nabla_i\eta_j 
-\rho V^k{}_lP^l{}_{i,j}\eta_k 
\end{equation}
and the curvature tensors by
\begin{equation}\label{flatcurvrelation}
{}^gR_{ijkl}={}^\ga R_{abkl}U^a{}_iU^b{}_j.
\end{equation}
\end{lemma}
\begin{proof}
The right hand side of \eqref{difften} defines a torsion-free connection,
so for \eqref{difften} it suffices to  
show that $g$ is parallel.  Differentiating \eqref{gformula} gives
\begin{equation}\label{derivg}
{}^\ga\nabla_kg_{ij}  = 2\rho P_{ij,k} +2\rho^2 P^l{}_{(i}P_{j)l,k}.  
\end{equation}
On the other hand, using \eqref{Vg}, \eqref{Udef} and the symmetry of
$P_{ij,k}$ gives 
$$
V^m{}_lP^l{}_{k,(i}g_{j)m}=P^l{}_{k,(i}U_{j)l}=P_{ij,k}
+\rho P^l{}_{k,(i}P_{j)l}. 
$$
Combining these shows that $g$ is parallel. 

Set $D^i{}_{jk}=\rho V^i{}_lP^l{}_{j,k}$.  The difference of the  
curvature tensors of the connections is given in terms of the difference
tensor by
$$
{}^gR_{mjkl}g^{im}={}^\ga R^i{}_{jkl} + 2D^i{}_{j[l,k]}
+2D^c{}_{j[l}D^i{}_{k]c}.
$$
Now
\[
\begin{split}
D^i{}_{jl,k}
&=\rho\left(V^i{}_{a,k}P^a{}_{j,l}+V^i{}_aP^a{}_{j,lk}\right)\\  
&=-\rho V^i{}_bU^b{}_{c,k}V^c{}_aP^a{}_{j,l} +\rho V^i{}_aP^a{}_{j,lk}\\ 
&=-\rho^2 V^i{}_bP^b{}_{c,k}V^c{}_aP^a{}_{j,l} +\rho V^i{}_aP^a{}_{j,lk}\\ 
&=-D^i{}_{ck}D^c{}_{jl}+\rho V^i{}_aP^a{}_{j,lk} 
\end{split}
\]
so
\[
\begin{split}
{}^gR_{mjkl}g^{im}
&={}^\ga R^i{}_{jkl} + 2\rho V^{ia}P_{aj,[lk]}\\  
&={}^\ga R^i{}_{jkl} + \rho V^{ia}\left(  
{}^\ga R^b{}_{alk}P_{bj}+{}^\ga R^b{}_{jlk}P_{ab}
\right).
\end{split}
\]
Using \eqref{Vg}, we obtain
\[
\begin{split}
{}^gR_{ijkl}
&={}^\ga R^b{}_{jkl}g_{bi} 
+\rho U^a{}_i\left({}^\ga R^b{}_{alk}P_{bj} 
+{}^\ga R^b{}_{jlk}P_{ab}\right)\\
&={}^\ga R^b{}_{jkl}U_{ba}U^a{}_i 
+\rho U^a{}_i\left({}^\ga R^b{}_{alk}P_{bj} 
+{}^\ga R^b{}_{jlk}P_{ab}\right)\\
&=\left[ {}^\ga R^b{}_{jkl}\left(U_{ab}-\rho P_{ab}\right)     
+ {}^\ga R^b{}_{alk}(\rho P_{bj})\right] U^a{}_i \\   
&=\left[ {}^\ga R^b{}_{jkl}\ga_{ab}
+ {}^\ga R_{balk}(\rho P^b{}_j)\right] U^a{}_i \\   
&=\left[ {}^\ga R_{ajkl}
+ {}^\ga R_{abkl}(\rho P^b{}_j)\right]U^a{}_i \\ 
&=\left[ {}^\ga R_{abkl}\left(\delta^b{}_j+\rho P^b{}_j\right)\right] 
U^a{}_i \\ 
&={}^\ga R_{abkl}U^b{}_jU^a{}_i,
\end{split}
\]
which is \eqref{flatcurvrelation}.
\end{proof}

We remark that \eqref{flatcurvrelation} can be viewed as a compatibility 
condition on a solution $P_{ij}$ of the system $P_{i[j,k]}=0$, which is
overdetermined if $n\geq 3$.  

\bigskip
\noindent
{\it Proof of Proposition~\ref{ss}}.
According to \eqref{curv}, $\gt$ is flat if and only if 
\begin{equation}\label{firstone}
{}^gR_{ijkl} +\frac12(g_{il}g_{jk}'+g_{jk}g_{il}'
-g_{ik}g_{jl}'-g_{jl}g_{ik}')+\frac{\rho}{2}(g_{ik}'g_{jl}'-g_{il}'g_{jk}')=0
\end{equation}
\begin{equation}\label{secondone}
{}^g\nabla_kg_{ij}' - {}^g\nabla_jg_{ik}'=0 
\end{equation}
\begin{equation}\label{thirdone}
g_{ij}'' -   \frac12g^{pq}g_{ip}'g_{jq}'=0.
\end{equation}
We will verify \eqref{firstone}--\eqref{thirdone} by direct calculation.       

Begin with \eqref{thirdone}.  Differentiating \eqref{gformula} gives 
$g''_{ij}= 2P_{ik}P^k{}_j$.  Now $g^{ij}=V^{il}V_l{}^j$, so using   
\eqref{gderiv} gives $g^{pq}g_{ip}'g_{jq}'=4P_{ik}P^k{}_j$.  This
proves \eqref{thirdone}.

Next consider \eqref{secondone}.  Differentiating \eqref{derivg} with 
respect to $\rho$ gives
$$
{}^\ga\nabla_kg'_{ij}=
2P_{ij,k} +4P^l{}_{(i}P_{j)l,k}\,\rho.
$$
Using \eqref{gderiv} gives
$$
V^m{}_lP^l{}_{k,(i}g'_{j)m} =2V^m{}_lP^l{}_{k,(i}P_{j)a}U^a{}_m
=2P^l{}_{k,(i}P_{j)l}. 
$$
Applying \eqref{difften} and combining terms gives 
${}^g\nabla_kg'_{ij}=2P_{ij,k}$, so \eqref{secondone} holds. 

To prove \eqref{firstone}, note first that an easy calculation shows that 
$$
g_{ij}-\tfrac12 \rho g'_{ij}=U_{ij}. 
$$
Thus we obtain
\[
\begin{split}
g_{i[l}g'_{k]j}+ g_{j[k}g'_{l]i}-\rho g'_{i[l}g'_{k]j} 
&=\left(g-\tfrac12 \rho g'\right)_{i[l}g'_{k]j}
+\left(g-\tfrac12 \rho g'\right)_{j[k}g'_{l]i}\\   
&=U_{i[l}g'_{k]j}+U_{j[k}g'_{l]i}.
\end{split}
\]
Substituting \eqref{gderiv} gives
$$
g_{i[l}g'_{k]j}+g_{j[k}g'_{l]i}-\rho g'_{i[l}g'_{k]j} 
=2U_i{}^aU_j{}^b\left(\ga_{a[l}P_{k]b}+\ga_{b[k}P_{l]a}\right).
$$
But the hypotheses of Proposition~\ref{ss} imply that 
$$
{}^\ga R_{abkl} = -2\left(  
\ga_{a[l}P_{k]b}+\ga_{b[k}P_{l]a}\right)
$$
in all dimensions:  this is automatically true for $n=3$, is true for
$n\geq 4$ because $\ga$ is locally    
conformally flat, and is true for $n=2$ because $P_i{}^i=\frac12 R$ and the
space 
of algebraic curvature tensors is one-dimensional.  \eqref{firstone}
follows upon combining with \eqref{flatcurvrelation}.   
\stopthm

We next sketch another proof of Proposition~\ref{ss} suggested by the 
discussion in \cite{SS} 
which avoids the calculations of Lemma~\ref{diffcurv}.  One first verifies
\eqref{thirdone} exactly as above, which is very simple.  According to 
\eqref{curv}, this means that $\Rt_{\infty ij \infty}=0$.  This implies
that the
Bianchi identity $\Rt_{\infty i[jk,\infty]}=0$ reduces to 
$\Rt_{\infty ijk,\infty}=0$.  Write the covariant derivative as 
the coordinate derivative plus
Christoffel symbols and use \eqref{cnr}.  One finds that  
the components $\Rt_{\infty ijk}$ satisfy a system of equations of the 
of the form 
$$
\pa_\rho \Rt_{\infty ijk} +A^{pqr}{}_{ijk}\Rt_{\infty pqr}=0,
$$
where $A^{pqr}{}_{ijk}$ is smooth.  For each point of $M$, this is a
linear system of ordinary differential equations in $\rho$, so $\Rt_{\infty
ijk}=0$ so long as this holds at $\rho =0$, which follows easily from the 
hypotheses on $P_{ij}$.  Finally one applies the same argument to the
Bianchi identity $\Rt_{ij[kl,\infty]}=0$ to deduce that $\Rt_{ijkl}=0$. 

A consequence of Proposition~\ref{ss} is the fact that for $(M,[\ga])$
locally conformally flat, an ambient metric can be chosen which is flat in
a 
neighborhood of $\cG\times \{0\}$ (as opposed to just to infinite order),
or equivalently that a Poincar\'e metric can be chosen which is hyperbolic.    

The arguments in \cite{SS} and \cite{PP} used to derive the form of a 
hyperbolic 
Poincar\'e metric also give uniqueness of the solution   
up to diffeomorphism in an open set, not just to infinite order.  This may
be of some interest in hyperbolic geometry, particularly the distinction
between the cases $n=2$ and $n>2$.  We formulate the  
Poincar\'e metric version rather than the ambient metric version of the
result.  
\begin{theorem}\label{flatnbhd}
Let $\ga$ be a smooth metric on a manifold $M$, and let  
$g_+$ be a hyperbolic Poincar\'e metric defined on $M_+^\circ$, where $M_+$
is a neighborhood of $M\times \{0 \}$ in $M\times [0,\infty)$, with 
conformal infinity $(M,[\ga])$.
Then there is a neighborhood $\cU$ of 
$M\times \{0\}$ in $M\times [0,\infty)$ and a diffeomorphism 
$\phi$ mapping $\cU$ into $M_+$ which restricts to the identity on 
$M\times \{0\}$, so that   
\begin{equation}\label{normalagain}
\phi^*g_+ = r^{-2}\left( dr^2 + g_r\right),
\end{equation}
where
$$
(g_r)_{ij}=\ga_{ij} - P_{ij}r^2 +\tfrac14 P_{ik}P^k{}_j r^4   
$$
and:
\begin{itemize}
\item
If $n\geq 3$,  $P$ is the Schouten tensor of $\ga$.       
\item
If  $n=2$, $P$ is some symmetric 2-tensor on $M$ satisfying 
$$
P_i{}^i=\tfrac12 R \qquad\text{and}\qquad P_{ij,}{}^j= \tfrac12 R_{,i}. 
$$
\end{itemize}
\end{theorem} 
\begin{proof}
Take $\cU$, $\phi$ as in Proposition~\ref{poincarenormal}, so that 
$\phi^*g_+$ is in normal form relative to $\ga$.  Since $g_+$ is
hyperbolic, the ambient metric determined by $\phi^*g_+$ by the change of 
variable $\rho = -\frac12 r^2$ is flat and in normal form.  Its curvature
is given by \eqref{curv}.  In particular, we have \eqref{thirdone}.
Differentiating \eqref{thirdone} with respect to $\rho$ and substituting  
\eqref{thirdone} for the second derivatives which occur gives
$g_{ij}'''=0$, where $' = \pa_\rho$.  Thus $g_{ij}$ is a quadratic
polynomial in $\rho$.  

We saw in Chapter~\ref{cct} that for $n\geq 3$ the Weyl and Cotton tensors of
$\ga$ can be recovered by restricting the curvature tensor of $\gt$ to
$\rho=0$.  Since $\gt$ is flat, $\ga$ is locally conformally flat.  
The uniqueness parts of Proposition~\ref{confflatprop} and
Theorem~\ref{2dim} now show that the $\rho$ and 
$\rho^2$ coefficients of $g_{ij}$ must be of the form given 
in Proposition~\ref{ss}.    
\end{proof}

In particular, Theorem~\ref{flatnbhd} implies that if $n\geq 3$, any 
two conformally compact hyperbolic metrics with the same conformal infinity 
are isometric in a neighborhood of the boundary by a diffeomorphism which
restricts to the identity on the boundary.  
The same result holds for $n=2$, 
except that one must require not only that the conformal infinities 
agree but that the $r^2$ coefficients $P_{ij}$ agree as well.  
The proof of Theorem~\ref{flatnbhd} is valid with weaker 
regularity hypotheses than our usual assumption that the compactification  
of $g_+$ is infinitely differentiable.  
There is an equivalent statement of Theorem~\ref{flatnbhd}
in terms of flat ambient metrics, whose formulation we leave to the
reader.  

Finally, we note that the argument of Theorem~\ref{flatnbhd} can be used to   
give yet another proof of Proposition~\ref{ss}.  
If $n\geq 3$ and $\ga$ is locally conformally flat,
Proposition~\ref{confflatprop} shows that there exists a 
straight ambient metric for $(M,[\ga])$ which is flat to infinite order.  
Put this metric into normal form relative to $\ga$.  The argument of 
Theorem~\ref{flatnbhd} shows that modulo terms vanishing to infinite order,
the corresponding $g_\rho$ must be a   
quadratic polynomial in $\rho$ and its coefficients must be given by
\eqref{gformula}.  If one takes $g_\rho$ to equal this quadratic
polynomial, then $\widetilde{g}$ is real analytic in $\rho$.  Since its
curvature tensor vanishes to infinite order at $\rho =0$, it must vanish
identically, so $\widetilde{g}$ is flat as desired.   
The same argument applies for $n=2$ as well, using Theorem~\ref{2dim} and
the discussion after Proposition~\ref{confflatprop} to conclude the
existence of a straight ambient metric having the prescribed
initial asymptotics which is flat to infinite order.  

Consider next the case of conformal classes $[g]$ containing an Einstein
metric.  For this discussion we assume $n\geq 3$.  We already   
noted in   
Chapter~\ref{formalsect} that if $g$ is Einstein, 
with $R_{ij}=2\la(n-1)g_{ij}$
so that $P_{ij}=\la g_{ij}$, 
then 
\begin{equation}\label{confein}
\gt= 2\rho dt^2 + 2tdtd\rho +t^2g_\rho, \qquad \qquad 
g_\rho = (1+\la \rho)^2g, 
\end{equation}
is straight, Ricci-flat, and in normal form relative to $g$.  In the
Poincar\'e realization, this is a reformulation of the familiar 
warped product construction of an Einstein metric in $n+1$ dimensions from
an Einstein metric in $n$ dimensions (see for example \cite{Be}):  for
$\la\neq 0$, under the change of 
variable 
$r =\sqrt{-2\rho}=\sqrt{\frac{2}{|\la|}}e^{-v}$, the associated Poincar\'e  
metric becomes
\begin{equation}\label{poincareeinstein}
g_+ =r^{-2}\left( dr^2 +(1-\tfrac{1}{2}\la r^2)^2g \right) 
=dv^2 +  2|\la|\, \psi(v)\,g, 
\end{equation}
where 
$$
\psi(v)=\left\{ 
\begin{array}{ll} \operatorname{sinh}^2v& \text{ for $\la >0$} 
\\ \operatorname{cosh}^2v& \text{ for $\la <0$}.  
\end{array} \right.
$$
For $\la =0$, the corresponding change of variable is $r=e^{-v}$, giving 
$g_+ = dv^2+ e^{2v}g$.  See also \cite{GrH1}, \cite{Leit}, \cite{Leis},
\cite{Ar}.  The explicit representation \eqref{confein} 
has been generalized to the case of certain products of Einstein metrics in
\cite{GoL}.  See also the remarkable examples of Nurowski \cite{N} of 
explicit ambient metrics.  

When $n$ is odd, then up to diffeomorphism \eqref{confein} is the unique
ambient metric to infinite order.  But when $n$ is even, there are others
satisfying $\Ric(\gt)=O(\rho^\infty)$, corresponding to other choices of   
$\tf\left(\partial_\rho^{n/2}g_{ij}|_{\rho = 0}\right)$ 
in Theorem~\ref{even}.
The following result shows that modulo diffeomorphism, the ambient metric
\eqref{confein} is uniquely determined to infinite order by the conformal
class alone.   
\begin{proposition}\label{ambientunique}
Let $g$ be an Einstein metric and let $\gt$ be the ambient metric defined
by \eqref{confein}.  Suppose that a conformal rescaling  
$\widehat{g} = e^{2\Up}g$ is also Einstein with constant $\widehat{\la}$,
and let $\widehat{\gt}$ denote the ambient metric defined by taking   
$g_\rho = (1+\widehat{\la} \rho)^2 \widehat{g}$ in
\eqref{confein}.  Then $\gt$ and $\widehat{\gt}$ are diffeomorphic to 
infinite order.    
\end{proposition}
\noindent
Of course, the Poincar\'e metrics 
\eqref{poincareeinstein} determined by $g$ and $\widehat{g}$ 
are diffeomorphic to infinite order as well.  As mentioned in the
introduction, Proposition~\ref{ambientunique} was obtained by 
Haantjes-Schouten.  

We begin the proof of Proposition~\ref{ambientunique} by characterizing the
ambient metric \eqref{confein} in terms of its curvature.   

\begin{proposition}\label{lotsinfinities}
Let $g$ be Einstein and let $\gt$ be defined by \eqref{confein}.  The
covariant derivatives of curvature of $\gt$ satisfy 
$
\Rt_{IJKL,M\cdots N}=0
$ 
if at most 3 of the indices $IJKLM\cdots N$ are between 1 and $n$. 
\end{proposition}
\begin{proof}
First, substituting $g_{\rho}$ into \eqref{curv} gives 
$\Rt_{\nf jkl}=0$ and $\Rt_{\nf ij\nf}=0$.  
Now calculating the covariant derivatives inductively using the
observations from \eqref{cnr} that  
$\Gat^A_{\nf\nf}=0$ for all $A$ and $\Gat^A_{i\nf}=0$ unless $1\leq A\leq
n$, one finds that 
$\Rt_{\nf ij\nf,\nf\cdots\nf}=0$ and $\Rt_{\nf ijk,\nf\cdots\nf}=0$.  
Also, recalling $\Rt_{0JKL}=0$ and calculating gives 
$\Rt_{\nf ij\nf,k}=0$.  Inductively differentiating this relation shows
that $\Rt_{\nf ij\nf,k\nf\cdots\nf}=0$.    

Consider now $\Rt_{IJKL,M\cdots N}$ with at most 3 indices between 1 and
$n$.  Recall from Proposition~\ref{Tcontract} that a 0 index can be 
removed at the expense of permuting the remaining indices.  Thus it can
assumed that none of the indices is 0.  
The symmetries of the curvature tensor then show that  
$\Rt_{IJKL,M\cdots N}$ vanishes unless at least two of $IJKL$ are between 1
and $n$.  Thus at most one of $M\cdots N$ can be between 1 and $n$.  We
have shown above that all such components vanish, except for components 
$\Rt_{\nf ij\nf,\nf \cdots \nf k\nf\cdots\nf}$ with at least one $\nf$ to 
the right of the comma before the $k$.  This can be shown to vanish by
commuting the $k$ to the left:  by the differentiated Ricci identity, one
has that
$$
\Rt_{\nf ij\nf,\,\underbrace{\scriptstyle{\nf\cdots\nf}}_{a} 
\nf k \underbrace{\scriptstyle{\nf\cdots\nf}}_{b}}- 
\Rt_{\nf ij\nf,\,\underbrace{\scriptstyle{\nf\cdots\nf}}_{a} k
  \nf\underbrace{\scriptstyle{\nf\cdots\nf}}_{b}} 
$$
is a sum of terms of the form 
$$
\Rt^P{}_{M\nf k,\nf\cdots\nf} 
\Rt....,......
$$
where $M$ is one of $i$, $j$, $\nf$ and $P$ contracts against one of the
indices of the second factor.  But it follows from what we have already
shown that $\Rt_{QM\nf k,\nf\cdots\nf} =0$ for all choices of $M$ and $Q$.   
\end{proof}

In particular, for $n$ even one has  
\begin{equation}\label{fixamb}
\Rt_{\nf ij\nf,\, \underbrace{\scriptstyle{\nf\cdots\nf}}_{n/2-2}}=0 
\quad\text{at } \rho = 0.
\end{equation}
Recall from Proposition~\ref{ambiguity} that this condition determines  
$\tf\left(\partial_\rho^{n/2}g_{ij}|_{\rho = 0}\right)$ (which clearly
vanishes in this case, so that $K_{ij}$ in Proposition~\ref{ambiguity}
vanishes for an Einstein metric), and by  
Theorem~\ref{even} this determines the solution to infinite order.  Thus we 
have:  
\begin{proposition}
Let $n\geq 4$ be even.  If $g$ is an Einstein metric, then $\gt$ given by    
\eqref{confein} is to infinite order the unique ambient metric  in normal 
form relative to $g$ 
satisfying $\Ric(\gt)=O({\rho^\nf})$ and \eqref{fixamb}.   
\end{proposition}

The proof of Proposition~\ref{ambientunique} uses relations between the 
curvature of an Einstein metric $g$ 
and the conformal factor relating $g$ to another Einstein metric.  Recall
from Chapter~\ref{cct} that if $r\geq 0$ and if we divide the indices
$IJKLM_1\cdots M_r$ into 
three disjoint subsets $\mathcal{S}_0$, $\mathcal{S}_M$, $\mathcal{S}_\nf$ 
of cardinalities $s_0$, $s_M$, $s_\nf$, resp., 
then we can define a covariant tensor
$\Rt_{\mathcal{S}_0,\mathcal{S}_M,\mathcal{S}_\nf}$ 
on $M$ of rank $s_M$ as follows.  Take the 
derivative of ambient curvature $\Rt_{IJKL,M_1\cdots M_r}|_{\rho = 0,t=1}$, 
set the 
indices in $\mathcal{S}_0$ to be 0's, the indices in $\mathcal{S}_\nf$ to
be $\nf$'s, and let those in $\mathcal{S}_M$ correspond to $M$ in the 
identification $\cG\times\R \cong \R_+ \times M \times \R$ induced by  
$g$.  

\begin{proposition}\label{relations}
Let $g$ and $\widehat{g}=e^{2\Up}g$ be conformally related Einstein
metrics.  Divide the indices $IJKLM_1\cdots M_r$ into
three disjoint subsets $\mathcal{S}_0$, $\mathcal{S}_M$, $\mathcal{S}_\nf$ 
and construct the tensor $\Rt_{\mathcal{S}_0,\mathcal{S}_M,\mathcal{S}_\nf}$ 
on $M$ as described above.  Now divide $\mathcal{S}_M$ into two disjoint
subsets: $\mathcal{S}_M = \mathcal{T}_1 \cup \mathcal{T}_2$, and set  
$\operatorname{card}(\mathcal{T}_i)=t_i$, $i=1,\,2$.  
Construct a covariant tensor $Y$ on $M$ of rank
$t_2$ by 
contracting the vector field $\operatorname{grad}\Up$ into each of the
indices of $\mathcal{T}_1$ in
$\Rt_{\mathcal{S}_0,\mathcal{S}_M,\mathcal{S}_\nf}$.  
If $t_2\leq 3$, then $Y=0$.    
\end{proposition}
\noindent
For example, Proposition~\ref{relations} asserts that 
$$
Y_{ijk} = \Up^a\Up^b\Up^c\Up^d \Rt_{i\nf 0a,\nf bc j\nf d\nf k}|_{\rho = 
  0,t=1}=0. 
$$
Observe that the case $\mathcal{T}_1 = \emptyset$ follows from  
Proposition~\ref{lotsinfinities}.  In this case there are no contractions 
with $\operatorname{grad}\Up$ and one need not restrict to $\{\rho = 0\}$.      

\begin{proof}
The proof is by induction on $r$. 
By applying Proposition~\ref{Tcontract}, we may as well assume   
$\mathcal{S}_0=\emptyset$.  By Proposition~\ref{lotsinfinities} we may
assume $t_1\geq 1$.  Also it suffices to consider the case $t_2 = 3$, since 
any $Y$ with $t_2<3$ may be viewed as a contraction against 
$\operatorname{grad}\Up$'s of a $Y$ with $t_2=3$
(Proposition~\ref{lotsinfinities} shows  
that $Y$ vanishes unless $s_M=t_1+t_2\geq 4$, so there must exist enough
indices contracted against $\operatorname{grad}\Up$'s).      

For $r=0$, it need only be shown that $\Up^l \Rt_{lijk}|_{\rho =
0,t=1}=0$. 
We have seen from \eqref{curv} that $\Rt_{lijk}|_{\rho = 0,t=1}=   
W_{lijk}$.  The desired conclusion follows from the conformal 
transformation law $\widehat{C}_{ijk} = C_{ijk} - \Up^l W_{lijk}$ for the
Cotton tensor and the fact that the Cotton tensor vanishes for Einstein
metrics.  

Supppose the result is true up to $r-1$.  Construct a tensor 
$\Rt_{\emptyset,\mathcal{S}_M,\mathcal{S}_\nf}$ from  
$\Rt_{IJKL,M_1\cdots M_r}|_{\rho = 0,t=1}$ as described above, where 
$s_M = t_1 +3$ and $t_1\geq 1$.  We need to show that the contraction of    
$\operatorname{grad}\Up$ into $t_1$ indices in $\mathcal{S}_M$ in 
the tensor 
$\Rt_{\emptyset,\mathcal{S}_M,\mathcal{S}_\nf}$ vanishes.  
Applying the Bianchi identity 
if necessary, we may assume that at least one of the indices of 
$\mathcal{S}_M$ is after the comma.  By commuting derivatives, we want to
arrange that $M_r\in \mathcal{S}_M$.  Take the rightmost element of
$\mathcal{S}_M$ in the list $M_1\cdots M_r$ and consider the effect of
commuting it successively past each of the $\nf$'s to its right.  According
to the differentiated Ricci identity, the expression
$$
\Rt_{IJKL,M_1\cdots M_s a\nf \underbrace{\scriptstyle{\nf\cdots\nf}}_{k}} - 
\Rt_{IJKL,M_1\cdots M_s \nf a \underbrace{\scriptstyle{\nf\cdots\nf}}_{k}}
$$
is a sum of terms of the form 
$$
\Rt^P{}_{Na\nf,\,\underbrace{\scriptstyle{\nf\cdots\nf}}_{l}}
\Rt....,......
$$
where $l\leq k$, $N\in \{IJKLM_1\cdots M_s\}$, and the second $\Rt$ has
some list of indices including $P$ but with $N$ removed.
Proposition~\ref{lotsinfinities} implies that 
$\Rt^P{}_{Na\nf,\nf\cdots\nf}=0$. 
Therefore the commutations introduce no new terms and we may assume that
$M_r\in \mathcal{S}_M$.  We write $M_r = a$.  

Now write
\begin{equation}\label{writeoutderiv}
\begin{split}
\Rt_{IJKL,M_1\cdots M_{r-1}a} = \partial_a &\Rt_{IJKL,M_1\cdots M_{r-1}}\\
&-\Gat_{aI}^P \Rt_{PJKL,M_1\cdots M_{r-1}} -\ldots 
-\Gat_{aM_{r-1}}^P \Rt_{IJKL,M_1\cdots P}.
\end{split}
\end{equation}
Consider the term $\Gat_{aI}^P \Rt_{PJKL,M_1\cdots M_{r-1}}$.  If $I=\nf$,
\eqref{cnr} shows that only terms with $1\leq P\leq n$ can 
contribute.  Write $p=P$.  The number of indices 
between $1$ and $n$ in the list $pJKLM_1\cdots M_{r-1}$ 
is the same as in 
$IJKLM_1\cdots M_{r-1}a$, so the induction hypothesis implies that this
term vanishes when
restricted to $\rho = 0$ and contracted against $t_1$ factors of 
$\operatorname{grad}\Up$.  Similarly for the other terms involving
Christoffel symbols.  Therefore we need consider only the terms involving
Christoffel symbols of the form $\Gat^P_{ab}$ with $b\in \mathcal{S}_M\cap
\{IJKLM_1\cdots M_{r-1}\}$.  If, for example, $I=b\in  
\mathcal{S}_M$, such a term is of the form  
$$
\Gat_{ab}^P \Rt_{PJKL,M_1\cdots M_{r-1}}=\Gat_{ab}^0 \Rt_{0JKL,M_1\cdots
  M_{r-1}} + \Gat_{ab}^p \Rt_{pJKL,M_1\cdots M_{r-1}}
+\Gat_{ab}^\nf \Rt_{\nf JKL,M_1\cdots M_{r-1}}. 
$$  
The number of elements of $\mathcal{T}_2$ in the list 
$JKLM_1\cdots M_{r-1}$ is of course at most 3.  
Therefore the induction hypothesis 
implies that upon setting $\rho = 0$ and contracting against the $t_1$
factors
of $\operatorname{grad}\Up$, the terms $\Gat_{ab}^0 \Rt_{0JKL,M_1\cdots
M_{r-1}}$ and $\Gat_{ab}^\nf \Rt_{\nf JKL,M_1\cdots   M_{r-1}}$ both
vanish.  Similarly for the other Christoffel symbol terms corresponding to
indices in $\mathcal{S}_M  \cap \{IJKLM_1\cdots M_{r-1}\}$.

We have shown that we need only include terms on the right
hand side of \eqref{writeoutderiv} corresponding to indices in 
$\mathcal{S}_M\cap \{IJKLM_1\cdots M_{r-1}\}$ and for such terms need only
sum over $p\in \{1, \ldots, n\}$.  Since these ambient
Christoffel symbols agree with the ones for the metric $g$ on $M$, 
this means that the right hand side can be interpreted as 
$\nabla_a \Rt_{IJKL,M_1\cdots M_{r-1}}$, where $\nabla_a$ is the covariant
derivative on $M$ and the tensor $\Rt_{IJKL,M_1\cdots M_{r-1}}$ is to be
interpreted as the rank $s_M - 1$ tensor on $M$ determined by fixing the
indices in $\mathcal{S}_\nf$ to be $\nf$ and letting the indices in
$\mathcal{S}_M \setminus \{a\}$ vary between $1$ and $n$.  Since the
derivative involved is tangential to $\cG$, we can restrict to $\rho =0$.  
We must show that we get 0 if we contract $t_1$ factors of 
$\operatorname{grad}\Up$ into $\nabla_a \Rt_{IJKL,M_1\cdots M_{r-1}}$. 
Use the Leibnitz rule to factor the $\nabla_a$ outside the product of 
all the $\operatorname{grad}\Up$ terms with $\Rt_{IJKL,M_1\cdots M_{r-1}}$, 
at the expense of the sum of the
terms where the derivative hits each one of the $\operatorname{grad}\Up$
factors in turn.  By the induction hypothesis, the contraction of all the  
$\operatorname{grad}\Up$ factors (except $\Up^a$ if it occurs) into
$\Rt_{IJKL,M_1\cdots M_{r-1}}$ 
vanishes.  To handle the terms where a derivative hits a
$\operatorname{grad}\Up$, recall the transformation law for the Schouten
tensor:  
$$
\widehat{P}_{ij} = P_{ij} -\Up_{ij} + \Up_i\Up_j - \frac12 \Up^k\Up_k
g_{ij}.
$$
We have $\widehat{P}_{ij} = \widehat{\la}\widehat{g}_{ij}
= \widehat{\la}e^{2\Up}g_{ij}$ and $P_{ij}=\la g_{ij}$.  
Therefore it follows that 
$$
\Up_{ij} = \Up_i\Up_j + fg_{ij}
$$ 
for some function $f$ on $M$.  
Substituting this relation for the second derivatives of $\Up$ which arise,
one sees easily that the induction hypothesis implies that all the terms 
vanish, completing the induction.
\end{proof}

\noindent
{\it Proof of Proposition~\ref{ambientunique}.} 
We have already observed that for $n$ odd, this follows from the uniqueness
of the ambient metric up to diffeomorphism.  So we can assume that $n\geq
4$ is even.

Recall that any pre-ambient metric for a given conformal class can be put
into 
normal form corresponding to any choice of representative metric.  So there 
is a diffeomorphism $\psi$ such that to infinite order, $\psi^*\gt$ is of  
the form \eqref{confein} but with $g_\rho$ replaced by some
$\widehat{g}_\rho$ satisfying $\widehat{g}_0= \widehat{g}=e^{2\Up}g$.
Moreover, the expansion of $\widehat{g}_\rho$ 
to order $n/2$ is uniquely determined by the Einstein condition 
$\Ric(\psi^*\gt) =0$ and the expansion to infinite order is determined once
$\tf\left(\partial_\rho^{n/2}\widehat{g}_\rho|_{\rho = 0}\right)$ is fixed.
It follows that  
$\widehat{g}_\rho = (1+\widehat{\la}\rho)^2\widehat{g} + O(\rho^{n/2})$ 
and $\widehat{g}_\rho = (1+\widehat{\la}\rho)^2\widehat{g} +
O(\rho^{\nf})$ if and only if the curvature tensor of $\psi^*\gt$ satisfies
\eqref{fixamb}.  We 
will show that the curvature tensor of $\psi^*\gt$ satisfies
\eqref{fixamb}, which will therefore prove the theorem.  

The curvature tensor of $\psi^*\gt$ is obtained from that of $\gt$ by
tranforming it tensorially.  The calculation of the Jacobian of $\psi$ and
this tensorial transformation are carried out in the proof of 
Proposition~\ref{ctransform}.  The result is that all components of all  
covariant derivatives of curvature of $\psi^*\gt$ are given by the right
hand side of \eqref{transformula}.  Thus \eqref{fixamb} is equivalent to: 
$$
\Rt_{ABCD,M_1\cdots M_{n/2-2}}|_{\rho = 0}
p^A{}_\nf p^B{}_ip^C{}_jp^D{}_\nf p^{M_1}{}_\nf\cdots p^{M_{n/2-2}}{}_\nf  
= 0,
$$
where $p^A{}_I$ is given by \eqref{pmatrix}.  
Expanding this out, one obtains a linear combination with smooth
coefficients of contractions of some number (possibly 0) of factors of
$\operatorname{grad}\Up$   
with tensors $\Rt_{\mathcal{S}_0,\mathcal{S}_M,\mathcal{S}_\nf}$ as 
considered in Proposition~\ref{relations}.  Each contraction which occurs
has $t_2\leq 2$, so by Proposition~\ref{relations}, it vanishes. 
\stopthm

The existence of a unique normalized infinite order ambient metric up to 
diffeomorphism in the  locally conformally flat and conformally Einstein
cases has as a consequence that other invariant objects 
exist for such structures which do not exist for general conformal 
structures in even dimensions.  An example is the family of ``conformally
invariant powers of the Laplacian'' constructed in \cite{GJMS}.  In
\cite{GJMS} it was shown that if $n\geq 3$ is odd and $k\in \mathbb{N}$,
then there is a natural scalar differential operator with principal part
$(-\Delta)^k$ depending on  
a metric $g$, which defines a conformally invariant operator 
$P_{2k}:\cE(-n/2+k)\rightarrow \cE(-n/2-k)$, where $\cE(w)$ denotes the
space of conformal densities of weight $w$.  
(Recall that our convention is $\Delta = \nabla^i\nabla_i$.)  When $n\geq
4$ is even, the   
same result holds with the restriction $k\leq n/2$, and it was shown in
\cite{GoH}, following the special case $k=3$, $n=4$ established in
\cite{Gr1}, that no such natural operator exists for $k>n/2$ for general 
conformal structures.  The conformal invariance of the operator produced by
the algorithm in \cite{GJMS} depends only on the existence and invariance
of the ambient metric.  So it follows that when $n$ is 
even and $(M,[g])$ is locally conformally flat or
conformally Einstein, there exists for all $k\geq 1$ an 
operator $P_{2k}:\cE(-n/2+k)\rightarrow \cE(-n/2-k)$ with principal part
$(-\Delta)^k$ determined solely by the conformal structure.  

One can write the operator $P_{2k}$ explicitly for an Einstein 
metric $g$.  In \cite{GJMS}, the operator $P_{2k}$ arose as the obstruction
to ``harmonic extension'' of the density.  We briefly review the
construction.  A density of weight $w$ can be
viewed as a homogeneous function of degree $w$ on $\cG$, and one asks for
an extension to a homogeneous function of degree $w$ on $\cG\times \R$
which solves $\Dt u=0$ to high order along $\cG\times \{0\}$, where $\Dt$ 
denotes the Laplacian with respect to an ambient metric $\gt$.  In the
identification $\cG\times \R \cong \R_+\times M\times \R$ induced by a
choice of 
metric $g$ in the conformal class, we can write $u = t^wf$, where 
$f=f(x,\rho)$ is a function whose restriction to $\rho =0$ is the initial
density on $M$.  The equation $\Dt u=0$ is  
equivalent to the following equation on $f$:
\begin{equation}\label{harmonic}
-2\rho f''+(2w+n-2-\rho g^{ij}g_{ij}')f' +\Delta_\rho f +\tfrac12 w 
g^{ij}g_{ij}'f=0,
\end{equation}
where ${}'$ denotes $\pa_\rho$, $g_{ij}(x,\rho)$ is as in
\eqref{spform}, and $\Delta_\rho$ denotes the Laplacian with respect to 
$g_{ij}(x,\rho)$ with $\rho$ fixed.  
For $w=-n/2 +k$, the derivatives $\pa_\rho^mf|_{\rho =0}$ for   
$1\leq m\leq  k-1$  
are determined successively by differentiating \eqref{harmonic} with
respect to $\rho$ at $\rho =0$, but the expression obtained by
differentiating the left hand side of \eqref{harmonic} $k-1$ times 
depends only the previously determined derivatives so defines an
obstruction to harmonic extension which can be expressed modulo a
normalizing constant as the invariant operator $P_{2k}$ applied to
$f(x,0)$.  For general metrics one must require $k\leq n/2$ if $n$ is even
because $g_{ij}$ is then only determined to order $n/2$, but this is not 
necessary if the initial metric is Einstein.  

Suppose now that the initial metric $g$ is Einstein with
$\Ric(g)=2\la(n-1)g$ as above.  Then
$g_{ij}(x,\rho)=(1+\la\rho)^2g_{ij}(x,0)$, so for $w=-n/2+k$, 
\eqref{harmonic} becomes 
\begin{equation}\label{specialharmonic}
-2\rho f''+\left(2k-2-\frac{2n\la\rho}{1+\la\rho}\right)f' 
+\frac{1}{(1+\la\rho)^2}\Delta f  
+ \frac{\la n(-n/2+k)}{1+\la\rho} f=0,
\end{equation}
where $\Delta=\Delta_0$ denotes the Laplacian with respect 
to the initial metric $g$.  It is clear from the form of this equation and
from the inductive procedure above that in this case 
the operator $P_{2k}$ takes the 
form $P_{2k}=p_{2k}(\Delta,\la)$, where $p_{2k}$ is a polynomial in two
real variables depending on $n$ and $k$ as parameters.  Moreover, rescaling
$\rho$ 
in \eqref{specialharmonic} and using the fact that $P_{2k}$ has leading
coefficient $(-\Delta)^k$ shows that  
$p_{2k}(\Delta,\la)=\la^kp_{2k}(\Delta/\la,1)$. 
Thus $P_{2k}$ for Einstein metrics is determined by the polynomial
$p_{2k}(\Delta,1)$.  In \cite{Br}, Branson showed that on $S^n$,  
any differential 
operator $P_{2k}$ which defines an equivariant map $:\cE(-n/2-k)\rightarrow 
\cE(-n/2+k)$ on densities viewed as sections of homogeneous vector bundles
for the conformal group is necessarily a multiple of  
$$
\prod_{j=1}^k\left(-\Delta +c_j\right),   
\qquad c_j=\left(\tfrac{n}{2}+j-1\right)\left(\tfrac{n}{2}-j\right) 
$$
when written with respect to the metric of constant sectional curvature 1.  
The operator $P_{2k}$ has this invariance property as a consequence of its
naturality and invariance under conformal rescaling, so this determines the
polynomial $p_{2k}$.  Thus we deduce 
\begin{proposition}\label{gjmseinstein}
If $n\geq 3$ and if $g$ is Einstein with $\Ric(g)=2\la (n-1)g$, then  
the operators $P_{2k}$ are given for $k\geq 1$ by  
$$
P_{2k}=\prod_{j=1}^k\left(-\Delta +2\la c_j\right),  
\qquad c_j=\left(\tfrac{n}{2}+j-1\right)\left(\tfrac{n}{2}-j\right).
$$
\end{proposition}
The above argument showing that the operators produced by the GJMS
algorithm for the ambient metric \eqref{confein} of an Einstein metric  
are given by the same formula as on the sphere was presented by the second
author at the 2003 
AIM Workshop on Conformal Structure in Geometry, Analysis, and Physics in
response to a question of Alice Chang (see \cite{BEGW}).  One can also
argue via Proposition~\ref{onlyricci} if $n$ is odd or if $n$ is even and
$k\leq n/2$.  In these cases, Proposition~\ref{onlyricci} implies that 
the coefficients  
of the $P_{2k}$ can be written solely in terms of the Ricci curvature and
its covariant derivatives.  For an Einstein metric, the covariant
derivatives of Ricci vanish, so the formula necessarily reduces to a
universal formula involving only $\Delta$ and $\lambda$.  Another    
treatment of these results has been given by Gover in   
\cite{Go2} using tractors.  

Proposition~\ref{gjmseinstein} can also be derived 
without reference to Branson's formula on the sphere directly from the
recursion relation obtained by successively  
differentiating \eqref{specialharmonic}.  
Upon setting $\la =1$, multiplying \eqref{specialharmonic} by $(1+\rho)^2$,
and differentiating  $m$ times, one obtains at $\rho =0$: 
\[
\begin{split}
2(k-m-1)f^{(m+1)}  & +\left[ \Delta +4m(k-m-n/2)+n(k-n/2)\right]f^{(m)}\\ 
& +m\left[2(m-1)(k-m+1-n)+n(k-n/2)\right]f^{(m-1)}=0.
\end{split}
\]
Set $(a)_0=1$ and $(a)_m = a(a+1)\ldots (a+m-1)$ for $a\in \C$, $m\in\N$.  
It follows that 
$$
(k-m)_mf^{(m)}(x,0)= q_m(y)f(x,0)\qquad  0\leq m\leq k-1
$$
and that 
\begin{equation}\label{Pform}
p_{2k}(\Delta,1)=2^kq_k(y),   
\end{equation}
where   
$y=-\frac12\left[ \Delta - n(n/2-1)\right]$ and the $q_m(y)$ are the 
polynomials of one variable determined by $q_{-1}=0$, $q_0=1$ and the
recursion relation 
\begin{equation}\label{qrecur}
\begin{split}
q_{m+1}=\left[ y\right.&-2m(k-m-n/2)\left.- n(k-1)/2\right]q_m\\ 
&-m(m-k)(m-1+n/2)(m-1+n/2-k)q_{m-1}
\end{split}
\end{equation}
for $m\geq 0$.  

Up to a normalizing factor, the polynomials $q_m$ are 
a case of the dual Hahn polynomials, a family of discrete orthogonal
polynomials  which may be expressed explicitly in terms of
${}_3F_2$ hypergeometric functions (see \cite{KM}, \cite{NSU} and the
appendix of \cite{AW}).  This leads to the following formula 
for the $q_m$:  
$$
q_m(y)=\sum_{l=0}^m(-1)^{m-l}(n/2+l)_{m-l}(k-m)_{m-l}
\begin{pmatrix}m\\l\end{pmatrix}\prod_{j=1}^l\left[y-j(j-1)\right],  
$$
where for $l=0$ the empty product is interpreted as 1.  Indeed, it is not
difficult to verify directly that $q_m$ defined by this formula satisfies  
\eqref{qrecur}.  Taking $m=k$ gives  
$$
q_k(y) = \prod_{j=1}^k\left[y-j(j-1)\right],  
$$
which via \eqref{Pform} gives Proposition~\ref{gjmseinstein}.   

We remark that Proposition~\ref{gjmseinstein} gives a formula
for Branson's $Q$-curvature for an Einstein metric.  From Branson's
original definition in \cite{Br}, one obtains immediately from  
Proposition~\ref{gjmseinstein} that if $n$ is even and $\Ric(g)=2\la
(n-1)g$, then 
$$
Q(g)=(2\la)^{n/2}(n-1)\prod_{j=1}^{n/2-1}(\tfrac{n}{2}+j-1)(\tfrac{n}{2}-j).  
$$


\section{Jet Isomorphism}\label{jet}

A fundamental result in Riemannian geometry is the jet isomorphism theorem
which asserts that at the origin in geodesic normal coordinates, the full
Taylor expansion of the metric may be recovered from the iterated 
covariant derivatives of curvature.  As a
consequence, one deduces that any local invariant of Riemannian metrics
has a universal expression in terms of the curvature tensor and its
covariant derivatives.    
Geodesic normal coordinates are determined up to the orthogonal group, 
so problems involving local invariants are reduced to purely algebraic
questions concerning invariants of the orthogonal group on tensors.  

Our goal in this chapter is to prove an analogous jet isomorphism theorem
for conformal geometry.  By making conformal changes, the Taylor
expansion of a metric in geodesic normal coordinates can be further
simplified, resulting in a ``conformal normal form'' for metrics about a 
point.  
The jet isomorphism theorem states that the map from the Taylor
coefficients of metrics in conformal normal form to the space of
all conformal curvature tensors, realized in terms of covariant derivatives 
of ambient curvature, is an isomorphism.  If $n$ is even, the theorem holds
only up to a finite order.  In the conformal case, the role
of the orthogonal group in Riemannian geometry is played by a parabolic
subgroup of the conformal group.  We assume throughout this chapter and the
next that $n\geq 3$.   

We begin by reviewing the Riemannian theorem in the form we will 
use it.  Fix a reference quadratic form 
$h_{ij}$ of signature $(p,q)$ on $\R^n$.  In the positive definite case 
one typically chooses $h_{ij}=\delta_{ij}$.  The background coordinates are
geodesic normal coordinates for a metric $g_{ij}$ defined near the origin
in $\R^n$ if and only if $g_{ij}(0)=h_{ij}$ and the radial vector field 
$\pa_r$ satisfies $\nabla_{\pa_r}\pa_r = 0$, which is equivalent to
$\Ga_{jk}^ix^jx^k=0$ or 
\begin{equation}\label{geonormal}
2\pa_kg_{ij}x^jx^k=\pa_ig_{jk}x^jx^k.
\end{equation}
It is easily seen that the coordinates are normal if and only if 
\begin{equation}\label{equivcond}
g_{ij}x^j=h_{ij}x^j.  
\end{equation}
In fact, set $F_i = g_{ij}x^j-h_{ij}x^j$ and observe that 
\begin{equation}\label{F}
\pa_kg_{ij}x^jx^k = x^j\pa_jF_i -F_i, \qquad 
\pa_ig_{jk}x^jx^k = \pa_i(x^jF_j)-2F_i.
\end{equation}
Clearly $F_i=0$ implies that \eqref{geonormal} holds so that the
coordinates are normal.  
Conversely, if the coordinates are normal, then the fact that $g$ and 
$\pa_r$ are parallel along radial lines implies $x^jF_j=0$.  
Substituting this and \eqref{F} into \eqref{geonormal} gives
$x^j\pa_j F_i =0$, which together with the fact that $F_i$ is smooth and
vanishes at the origin implies $F_i=0$. 

We are interested in the space of infinite order jets of metrics at the
origin.  Taylor expanding shows that \eqref{equivcond} holds to infinite   
order if and only if $g_{ij}(0) =h_{ij}$ and the derivatives of $g_{ij}$
for all orders satisfy  
$\pa_{(k_1}\ldots \pa_{k_r} g_{i)j}(0) = 0$.  Since a 3-tensor symmetric in
2 indices and skew in 2 indices must vanish, this implies in particular that  
all first derivatives of $g$ vanish at the origin.  
So we define the space $\cN$ of jets of metrics in geodesic normal
coordinates as follows.   
\begin{definition}
The space $\cN$ is the set of lists $(g^{(2)},g^{(3)},\ldots )$,
where for  
each $r$, $g^{(r)} \in \bigodot^2\R^n{}^*\otimes \bigodot^r\R^n{}^*$
satisfies  
$g^{(r)}_{i(j,k_1\cdots k_r)} = 0$.  Here the comma just serves to separate  
the first two indices.  For $N\geq 2$, $\cN^N$ will denote the set of  
truncated lists $(g^{(2)},g^{(3)},\ldots, g^{(N)})$ with the same
conditions on the $g^{(r)}$.
\end{definition}
\noindent
To a metric in geodesic normal coordinates near the
origin, we associate the element of $\cN$ given by 
$g^{(r)}_{ij,k_1\cdots k_r} =  \pa_{k_1}\ldots \pa_{k_r} g_{ij}(0)$
for $r\geq 2$.  Conversely, to an element of $\cN$ we associate the metric 
determined to infinite order by these relations together with 
$g_{ij}(0) = h_{ij}$ and $\pa_kg_{ij}(0) = 0$.  In the following, we
typically identify the element of $\cN$ and the jet of the metric $g$.

\begin{definition}\label{Rdef}
The space $\cR$ is the set of lists 
$(R^{(0)},R^{(1)},\ldots )$, where for each $r$, 
$R^{(r)} \in \bigwedge^2\R^n{}^*\otimes \bigwedge^2\R^n{}^*\otimes 
\bigotimes^r\R^n{}^*$, and the usual identities satisfied by covariant 
derivatives of curvature hold:
\begin{enumerate}
\item
$R_{i[jkl],m_1\cdots m_r}=0$ 
\item
$R_{ij[kl,m_1]m_2\cdots m_r}=0$
\item
$R_{ijkl,m_1\cdots [m_{s-1}m_s]\cdots m_r} =  
Q^{(s)}_{ijklm_1\cdots m_r}$\\
Here $Q^{(s)}_{ijklm_1\cdots m_r}$ denotes the quadratic polynomial in   
the $R^{(r')}$ with $r'\leq r-2$ which one 
obtains by covariantly differentiating the usual Ricci identity for
commuting covariant  
derivatives, expanding the differentiations using the Leibnitz rule, and
then setting equal
to $h$ the metric which contracts the two factors in each term.     
\end{enumerate}
For $N\geq 0$, $\cR^N$ will denote the set of truncated lists  
$(R^{(0)},R^{(1)},\ldots, R^{(N)} )$ with the same conditions on the
$R^{(r)}$. 
\end{definition}

There is a natural map $\cN\rightarrow \cR$ induced by evaluation of the
covariant derivatives of curvature of a metric, which is polynomial in the
sense that the corresponding truncated maps $\cN^{N+2}\rightarrow \cR^N$
are polynomial.  We will say that a map on a subset of a finite-dimensional 
vector space is 
polynomial if it is the restriction of a polynomial map defined on the
whole space, and a map between subsets of finite-dimensional vector spaces
is a polynomial equivalence if it is a bijective polynomial map whose
inverse is also polynomial.  
The jet    
isomorphism theorem for pseudo-Riemannian geometry is the following.
\begin{theorem}\label{riemjet}
The map $\cN\rightarrow \cR$ is bijective and the truncated  
maps $\cN^{N+2}\rightarrow \cR^N$ are polynomial equivalences.   
\end{theorem}

\noindent
There are two parts to the proof.  One is a linearization   
argument showing that it suffices to show that the linearized map 
is 
an isomorphism.  The other is the observation that the linearized map  
is the direct sum over $r$ of isomorphisms 
between two realizations according to different Young projectors of
specific irreducible representations of $GL(n,\R)$.  The reduction to the
linearization can be carried out in different ways; one is to argue as we
do below in the conformal case.  The analysis of the linearized map     
is contained in \cite{E}.  See also \cite{ABP} for a different  
construction of a left inverse of the map $\cN\rightarrow \cR$.    

We now consider the further freedom allowed by the possibility of making
conformal changes to $g$.  
\begin{proposition}\label{cnf}
Let $g$ be a metric on a manifold $M$ and let $p\in M$.  Given $\Omega_0\in 
C^{\infty}(M)$ with $\Omega_0(p)>0$, there is $\Omega \in C^{\infty}(M)$,
uniquely determined to infinite order at $p$, such that 
$\Omega -\Omega_0$ vanishes to second order at $p$ and 
$\operatorname{Sym}(\nabla^r\Ric)(\Omega^2g)(p)=0$ for all $r\geq 0$.   
Here $\operatorname{Sym}(\nabla^r\Ric)(\Omega^2g)$ denotes the full
symmetrization of the rank $r+2$ tensor $(\nabla^r\Ric)(\Omega^2g)$.  
\end{proposition}
\begin{proof}
Write $\Omega = e^{\Up}$ and set $\wh{g} = e^{2\Up}g$.  We are given 
$\Up(p)$ and $d\Up(p)$.  
Recall the transformation law for conformal change of Ricci curvature:  
$$
\wh{R}_{ij} = R_{ij} -(n-2)\Up_{ij} -\Up_k{}^kg_{ij} 
+(n-2)(\Up_i\Up_j -\Up_k\Up^kg_{ij}).
$$
Differentiating and conformally transforming the covariant derivative
results in
$$
\wh{(\nabla^r\Ric)}_{ij,m_1\cdots m_r} = (\nabla^r\Ric)_{ij,m_1\cdots m_r} 
 -(n-2)\Up_{ijm_1\cdots m_r} -\Up_k{}^k{}_{m_1\cdots m_r}g_{ij} +
 \mbox{lots}, 
$$
where lots denotes terms involving at most $r+1$ derivatives of $\Up$.  
We may replace covariant derivatives of $\Up$ by coordinate derivatives on
the right hand side of this formula.  Then symmetrizing gives 
$$
\operatorname{Sym}\wh{(\nabla^r\Ric)} = 
\operatorname{Sym}(\nabla^r\Ric)
 -(n-2)\pa^{r+2}\Up
-\operatorname{Sym}(\operatorname{tr}(\pa^{r+2}\Up)g) +  
 \mbox{lots}.
$$
If $c>0$, the map $s\rightarrow
s+c\operatorname{Sym}(\operatorname{tr}(s)g)$ is invertible on 
symmetric $(r+2)$-tensors.  Therefore, for any $r\geq 0$, one can uniquely 
determine $\pa^{r+2}\Up(p)$ to make
$\operatorname{Sym}\wh{(\nabla^r\Ric)}(p) =0$.  The result follows upon
iterating and applying Borel's Lemma.  
\end{proof}
\noindent
We remark that there are other natural choices for normalizations of the
conformal factor.  For example, one such is that the 
symmetrized covariant derivatives of the tensor $P_{ij}$ vanish at $p$,
where $P_{ij}$ is given by \eqref{Ptensor}.  Another is that in normal
coordinates, $\det{g_{ij}}-1$ vanishes to infinite order at $p$ (see
\cite{LP}).   

\begin{definition}\label{cnfdef}
The space $\cN_c\subset \cN$ of jets of metrics in conformal normal
form is the subset of $\cN$ consisting of jets of metrics in geodesic
normal coordinates for which   
$\operatorname{Sym}(\nabla^r\Ric)(g)(0)=0$ for all $r\geq 0$. 
For $N\geq 2$, $\cN_c^N$ will denote the subset of $\cN^N$ obtained by
requiring that this relation hold for $0\leq r\leq N-2$.   
\end{definition}

Later we will need the following consequence of the proof of
Proposition~\ref{cnf}. 
\begin{lemma}\label{extend}
For each $N\geq 3$, there is a polynomial map 
$\eta^N:\cN_c^{N-1}\rightarrow \bigodot^2\R^n{}^*\otimes
\bigodot^N\R^n{}^*$ with zero constant  
and linear terms, such that the map $g\rightarrow (g,\eta^N(g))$
maps $\cN_c^{N-1}\rightarrow \cN_c^N$.  Here  
$g= (g^{(2)},g^{(3)},\ldots, g^{(N-1)})\in \cN_c^{N-1}$.  
\end{lemma}
\begin{proof}
View $g$ as the metric on a neighborhood of $0\in \R^n$ given by the
prescribed finite Taylor expansion of order $N-1$.  Then the components of  
$\operatorname{Sym}(\nabla^{N-2}\Ric)(g)(0)$ are polynomials in the  
$g^{(r)}$ with no constant or linear terms.  According to the proof of
Proposition~\ref{cnf}, there is a function $\Omega$ of the form
$\Omega = 1+p_{N}$ with $p_{N}$ a homogeneous polynomial of degree
$N$ whose coefficients are polynomials in the $g^{(r)}$ with no constant
or linear terms, such that
$\operatorname{Sym}(\nabla^{N-2}\Ric)(\Omega^2g)(0)=0$.  Now $\Omega^2g$
is in geodesic normal coordinates to order $N-1$ but not necessarily $N$.   
However by the construction of geodesic normal coordinates, there is a 
diffeomorphism $\psi = I + q_{N+1}$, where $q_{N+1}$ is a vector-valued 
homogeneous polynomial of order $N+1$ whose coefficients are linear in the
order $N$ Taylor coefficient of $\Omega^2g$, such that 
$\psi^*(\Omega^2g)$ is in geodesic normal coordinates to order $N$.   
Since the condition $\operatorname{Sym}(\nabla^{r}\Ric)(\Omega^2g)(0)=0$  
is invariant under diffeomorphisms, $\eta^N$ 
defined to be the order $N$ Taylor coefficient of $\psi^*(\Omega^2g)$
has the required properties. 
\end{proof}

If $g_1$ and $g_2$ are jets of metrics of signature $(p,q)$ at the origin
in $\R^n$, we say that $g_1$ and $g_2$ are equivalent if
there is a local diffeomorphism $\psi$ defined near $0$ satisfying 
$\psi(0)=0$, and a positive smooth function $\Omega$ defined near $0$, so
that $g_2 = \psi^*(\Omega^2 g_1) $ to infinite order.  It is  
clear from Proposition~\ref{cnf} and the existence of geodesic normal
coordinates that any jet of a metric is equivalent to one in $\cN_c$.  
In choosing $\Omega$ we have the freedom of $\R_+ \times \R^n$, and in
choosing $\psi$ a freedom of $O(p,q)$.  We next describe how these
freedoms can be realized as an action on $\cN_c$ of a subgroup of the
conformal group.    

Recall that $h_{ij}$ is our fixed background quadratic form of signature
$(p,q)$ on $\R^n$.  Define a quadratic form $\htt_{IJ}$ on $\R^{n+2}$ by 
$$
\htt_{IJ} = 
\left(
\begin{matrix}
0&0&1\\
0&h_{ij}&0\\
1&0&0
\end{matrix}
\right)
$$
and the quadric 
$\mathcal{Q} = \{[x^I]:\htt_{IJ}x^Ix^J = 0\}\subset \mathbb{P}^{n+1}$.  
The metric $\htt_{IJ}dx^Idx^J$ on $\R^{n+2}$ induces a conformal structure
of signature $(p,q)$ on $\mathcal{Q}$.  
The standard action of the orthogonal group $O(\htt)$ on $\R^{n+2}$ induces
an action on $\mathcal{Q}$ by conformal tranformations and the adjoint
group $O(\htt)/\{\pm I\}$ can be identified with the conformal group of all
conformal transformations of $\mathcal{Q}$.  Let 
$e_0 = 
\left(
\begin{matrix}
1\\
0\\
0
\end{matrix}
\right)\in \R^{n+2}$ and let $P$ be the image in $O(\htt)/\{\pm I\}$ of the
isotropy group $\{A\in O(\htt): Ae_0 = a e_0, a \in \R\}$ of
$[e_0]$.  It is clear that each element of $P$ is represented by exactly
one $A$ for which $a >0$, so we make the identification
$P=\{p\in O(\htt): pe_0 = a e_0, a>0\}$.  The first column of 
$p\in P$ is 
$\left(
\begin{matrix}
a\\
0\\
0
\end{matrix}
\right)$; combining this with the fact that $p\in O(\htt)$, one finds that 
$$
P=\left\{\left(
\begin{matrix}
a&b_j&c\\
0&m^i{}_j&d^i\\
0&0&a^{-1}
\end{matrix}
\right):
a>0,\, m^i{}_j\in O(h),\, c=-\frac{1}{2a}b_jb^j,\,d^i =
-\frac{1}{a}m^{ij}b_j
\right\},
$$
where $h_{ij}$ is used to raise and lower lower case indices.  It is
evident that $P=\R_+\cdot \R^n \cdot O(h)$, where the subgroups 
$\R_+$, $\R^n$, $O(h)$ arise by varying $a$, $b_j$, $m^i{}_j$, resp.  

The intersection of $\mathcal{Q}$ with the cell $\{[x^I]: x^0\neq 0\}$ 
can be identified with $\R^n$ via 
$\R^n\ni x^i\rightarrow 
\left[
\begin{matrix}
1\\
x^i\\
-\frac12 |x|^2
\end{matrix}
\right]\in \mathcal{Q}$, where $|x|^2 = h_{ij}x^ix^j$.   
In this identification, the conformal structure is represented by the 
metric $h_{ij} dx^idx^j$ on $\R^n$.  If $p\in P$ is as above, the conformal 
transformation determined by $p$ is realized as 
$$
(\varphi_p(x))^i= \frac{m^i{}_jx^j -\frac12 |x|^2 d^i}{a+b_jx^j-\frac12 c
|x|^2}
$$
and one has $\varphi_p^*h = \Omega_p^2 h$
for 
$$
\Omega_p = (a+b_jx^j-\tfrac12 c|x|^2)^{-1}.
$$

This motivates the following definition of an action of $P$ on $\cN_c$.  
Given $p\in P$ as above and $g\in \cN_c$, by Proposition~\ref{cnf} there is
a positive smooth function $\Omega$ uniquely determined to infinite order
at $0$ so that $\Omega$ agrees with $\Omega_p$ to second order and such
that 
$\operatorname{Sym}(\nabla^r\Ric)(\Omega^2g)(0)=0$ for all $r\geq 0$.   
Now $(\Omega^2g)(0) = a^{-2}h$, so by the construction of geodesic normal
coordinates, there is a diffeomorphism $\varphi$, uniquely determined to
infinite order at $0$, so that $\varphi(0) = 0$, $\varphi'(0) = a^{-1}m$, 
and such that $(\varphi^{-1})^*(\Omega^2g)$ is in geodesic normal
coordinates to 
infinite order.  We define $p.g = (\varphi^{-1})^*(\Omega^2g)$.  It is
clear by  
construction of $\varphi$ that $p.g\in \cN$, and since the condition of 
vanishing of the 
symmetrized covariant derivatives of Ricci curvature is
diffeomorphism-invariant, it follows by construction of $\Omega$ that  
$p.g\in \cN_c$.  It is straightforward to check that this defines a left
action of $P$ on $\cN_c$.  Note that if $g=h$, then $\Omega = \Omega_p$ 
and $\varphi = \varphi_p$, so that $h$ is a fixed point of the action.   
A moment's thought shows that $(p.g)^{(r)}$ depends only on $g^{(s)}$ for 
$s\leq r$.  Therefore for each $N\geq 2$, there is an induced action on 
$\cN_c^N$.  

It is clear from the construction of the action that $p.g$ is equivalent to
$g$ for all $g\in \cN_c$ and $p\in P$.  In fact, the $P$-orbits are exactly 
the equivalence classes.
\begin{proposition}\label{orbits}
The orbits of the $P$-action on $\cN_c$ are precisely the equivalence
classes of jets of metrics in $\cN_c$ under diffeomorphism and conformal
change. 
\end{proposition}
\begin{proof}
It remains to show that equivalent jets of metrics in $\cN_c$ are in the
same $P$-orbit.  Suppose that $g_1$, $g_2 \in \cN_c$ are equivalent.  
Then we can write $g_2= (\varphi^{-1})^*(\Omega^2g_1)$ to infinite order
for a 
diffeomorphism $\varphi$ with $\varphi(0)=0$ and a positive smooth
$\Omega$.  We can uniquely choose the parameters $a$ and $b$ 
of $p\in P$ so that $\Omega -\Omega_p$ vanishes to second order.  
Since $g_1$ and $g_2$ both equal $h$ at $0$, it follows that 
$a\varphi'(0) \in O(h)$, so that we can write  
$\varphi'(0) = a^{-1}m$ with $m\in O(h)$.  Together with the already
determined parameters $a$ and $b$, this choice of $m$ uniquely determines 
a $p\in P$.  Since $g_2\in \cN_c$, all 
symmetrized covariant derivatives of Ricci curvature of $\Omega^2g_1$
vanish at the origin, so $\Omega$ must be the conformal factor determined
when constructing the action of $p$ on $g_1$.  And since $g_2$ is in
geodesic normal coordinates, $\varphi$ must be the correct diffeomorphism, 
so that $g_2=p.g_1$.  
\end{proof}

It is straightforward to calculate from the definition the action on
$\cN_c$ of the $\R_+$ and $O(h)$ subgroups of $P$.  If we denote by $p_a$ 
the element of $P$ obtained by taking $b=0$ and $m=I$ and by $p_m$ the
element given by $a=1$ and $b=0$, then one finds that $p_a$ acts by
multiplying $g^{(r)}$ by $a^r$, and 
$p_m$ acts by transforming each $g^{(r)}$ as an element of  
$\bigotimes^{r+2}\R^n{}^*$, where $\R^n$ denotes the standard defining  
representation of $O(h)$.  

The problem of understanding local invariants of metrics under 
diffeomorphism and conformal change reduces to understanding this
action of 
$P$ on $\cN_c$.  However, it is very difficult to analyze or even
concretely exhibit the action of the $\R^n$ part of $P$ directly from the
definition.  The ambient curvature tensors enable the reformulation of the
action in terms of standard tensor representations of $P$.  

Propositions~\ref{Tcontract} and \ref{tracefree} show that 
the covariant derivatives of curvature of an ambient metric $\gt$ satisfy
relations arising from the homogeneity and Ricci-flatness of the metric.
These conditions suggest the following definition.  
\begin{definition}\label{Rt}
The space $\cRt$ is the set of lists 
$(\Rt^{(0)},\Rt^{(1)},\ldots )$, such that 
$\Rt^{(r)} \in \bigwedge^2\R^{n+2}{}^*\otimes
\bigwedge^2\R^{n+2}{}^*\otimes \bigotimes^r\R^{n+2}{}^*$, and such that the
following relations hold: 
\begin{enumerate}
\item
$\Rt_{I[JKL],M_1\cdots M_r}=0$
\item
$\Rt_{IJ[KL,M_1]M_2\cdots M_r}=0$
\item
$\htt^{IK}\Rt_{IJKL,M_1\cdots M_r}=0$ 
\item
$\Rt_{IJKL,M_1\cdots [M_{s-1}M_s]\cdots M_r}=
\tilde{Q}^{(s)}_{IJKLM_1\cdots M_r}$
\item
$\Rt_{IJK0,M_1\cdots M_r}= 
-\sum_{s=1}^r \Rt_{IJKM_s,M_1\cdots \widehat{M_s} \cdots M_r}$
\item
$\Rt_{IJKL,M_1\cdots M_s0M_{s+1}\cdots M_r}= \\
{}\qquad\qquad
-(s+2)\Rt_{IJKL,M_1\cdots M_r} - 
\sum_{t=s+1}^r \Rt_{IJKL,M_1\cdots M_s M_t M_{s+1}\cdots
  \widehat{M_t}\cdots M_r}$.
\end{enumerate}
Here, as in Definition~\ref{Rdef}, $\tilde{Q}^{(s)}_{IJKLM_1\cdots M_r}$  
denotes the quadratic polynomial in the components of  
$\Rt^{(r')}$ for $r'\leq r-2$ which one
obtains by differentiating the Ricci identity for commuting covariant 
derivatives, expanding the differentiations using the Leibnitz rule, 
and then setting equal
to $\htt$ the metric which contracts the two factors in each term.       
Condition (5) in case $r=0$ is interpreted as $\Rt_{IJK0}=0$.   
\end{definition}
\noindent
We remark, as we did in the proof of Proposition~\ref{Tcontract}, that 
condition (6) is superfluous:  it is
a consequence of (2), (4), (5).  But we will not use this fact.  It will be
convenient also to introduce the vector space $\cTt$ consisting of the set
of lists $(\Rt^{(0)},\Rt^{(1)},\ldots )$ with  
$\Rt^{(r)} \in \bigwedge^2\R^{n+2}{}^*\otimes
\bigwedge^2\R^{n+2}{}^*\otimes \bigotimes^r\R^{n+2}{}^*$ such that  
(1)--(3), (5), (6) hold.  

We prepare to define truncated spaces $\cRt^N$ for $\cRt$.  
Recall the notion of strength from Definition~\ref{strength}.
Note that it is clear that for each of the relations (1)--(6) except (4) in 
Definition~\ref{Rt}, all components $\Rt_{IJKL,M_1\cdots M_r}$ which 
occur in the relation have the same value for the strength of the index 
list of the component.  For $N\geq 0$, define the following vector spaces
of lists of components of tensors.  We denote by $\cM=M_1\cdots M_r$ a list
of indices of length $|\cM|=r$.  
\[
\begin{aligned}
\cTt^N & = \{(\Rt_{IJKL,\cM})_{|\cM|\geq 0,\,\,\|IJKL\cM\|\leq N+4}:  
\mbox{ (1)--(3), (5), (6) of Definition~\ref{Rt} hold}\} \\
\ttt^N & = \{(\Rt_{IJKL,\cM})_{|\cM|\geq 0,\,\,\|IJKL\cM\|= N+4}:
\mbox{ (1)--(3), (5), (6) of Definition~\ref{Rt} hold}\}
\end{aligned}
\]
If $\Rt_{IJKL,M_1\cdots M_r}$ is a component appearing in an element of
$\cTt^N$ and $r>N$, then at least one of the indices $IJKLM_1\cdots M_r$
must be $0$.  Therefore by (5) and (6), $\Rt_{IJKL,M_1\cdots M_r}$ can be
written as a linear combination of components with $r$ replaced by $r-1$.  
It follows that $\cTt^N$ and $\ttt^N$ are finite-dimensional.
Since (1)--(3), (5) and (6) imply that $\Rt_{IJKL,M_1\cdots M_r}=0$ if
$\|IJKLM_1\cdots M_r\|\leq 3$, we have $\cTt^N = \bigoplus_{M=0}^N\ttt^M$  
and $\cTt = \prod_{M=0}^\infty\ttt^M$.

As for (4), a typical term in
$\tilde{Q}^{(s)}_{IJKLM_1\cdots M_r}$ is 
$\htt^{AB}\Rt_{AIM_{s-1}M_s,\,\mathcal M'}
\Rt_{BJKL,M_1\cdots M_{s-2}\mathcal M''}$,
where $\mathcal M'$ and $\mathcal M''$ are lists of indices such that 
$\mathcal M' \mathcal M''$ is a rearrangement of $M_{s+1}\cdots M_r$.  
In order that $\htt^{AB}\neq 0$, it must be that $\|AB\|=2$.  
Therefore $\|AIM_{s-1}M_s\mathcal M'\| + 
\|BJKLM_1\cdots M_{s-2}\mathcal M''\|  = \|IJKLM_1\cdots M_r\| +2$.
This implies that if either
$\|AIM_{s-1}M_s\mathcal M'\|$ or $\|BJKLM_1\cdots M_{s-2}\mathcal M''\|$ is 
greater than or equal to $\|IJKLM_1\cdots M_r\|$, then 
the other is less than or equal to 2.  
The same reasoning applies to all terms in   
$\tilde{Q}^{(s)}_{IJKLM_1\cdots M_r}$.
Since (1)--(3), (5) and (6) imply that a component of an
$\Rt^{(t)}$ vanishes if its strength is at most 3, it follows that any
component of an $\Rt^{(t)}$ which occurs in  
$\tilde{Q}^{(s)}_{IJKLM_1\cdots  M_r}$ with a nonzero coefficient must have
strength strictly less than $\|IJKLM_1\cdots M_r\|$.  
Hence we will regard (4) as a relation 
involving components with indices of strength at most  
$\|IJKLM_1\cdots M_r\|$, and the quadratic terms in (4) only involve
components of strength less than that of the linear terms.  With this
understanding, we can now define  
\[
\cRt^N = \{(\Rt_{IJKL,M_1\cdots M_r})\in \cTt^N: \mbox{ (4) of
  Definition~\ref{Rt} holds}\}.
\]

We will also need the corresponding linearized spaces, in which the
quadratic term $\tilde{Q}^{(s)}_{IJKLM_1\cdots M_r}$ in (4) is replaced
by $0$.  Define vector spaces
\[
\begin{aligned}
T\cRt &= \{(\Rt^{(0)},\Rt^{(1)},\ldots )\in \cTt:
\Rt_{IJKL,M_1\cdots M_r} = \Rt_{IJKL,(M_1\cdots M_r)}\}\\
T\cRt^N &= \{(\Rt_{IJKL,M_1\cdots M_r})\in \cTt^N: 
\Rt_{IJKL,M_1\cdots M_r} = \Rt_{IJKL,(M_1\cdots M_r)}\}\\
\st^N &=\{(\Rt_{IJKL,M_1\cdots M_r})\in \ttt^N:
\Rt_{IJKL,M_1\cdots M_r} = \Rt_{IJKL,(M_1\cdots M_r)}\}
\end{aligned}
\]
Then $T\cRt^N=\bigoplus_{M=0}^N\st^M$ and $T\cRt=\prod_{M=0}^\infty\st^M$.   
Note that in the presence of the condition
$\Rt_{IJKL,M_1\cdots M_r} = \Rt_{IJKL,(M_1\cdots M_r)}$,  
condition (6) of Definition~\ref{Rt} becomes 
\begin{equation}\label{remove0}
\Rt_{IJKL,M_1\cdots M_r0} = -(r+2)\Rt_{IJKL,M_1\cdots M_r}
\end{equation}
and  is a consequence of (2) and (5).  Note also that if 
$(\Rt_{IJKL,M_1\cdots M_r}) \in \st^N$, then  
$\htt^{LM_s}\Rt_{IJKL,M_1\cdots M_r} = 0$ and  
$\htt^{M_sM_t}\Rt_{IJKL,M_1\cdots M_r} = 0$.  

For $w\in \C$, let $\sigma_w:P\rightarrow \C$ denote the character
$\sigma_w(p) = a^{-w}$.  Since $P\subset O(\htt)$ and $P$ preserves $e_0$
up to scale, it is easily seen that $\cTt$, $\cRt$ and $T\cRt$ are
invariant subsets of the $P$-space 
$\prod_{r=0}^{\infty}\left(\bigotimes^{4+r}\R^{n+2}{}^*
\otimes \sigma_{-2-r}\right)$,
where $\R^{n+2}$ denotes the standard representation of $P\subset
GL(n+2,\R)$. 
These inclusions therefore define actions of $P$ on these spaces.  These
actions of 
$P$ do not preserve strength, but because $P$ consists of block
upper-triangular matrices, a component of $p.(\Rt)$
depends only on components of $\Rt$ of no greater 
strength.  So for $N\geq 0$ there are also actions of $P$ on $\cTt^N$,
$\cRt^N$ and $T\cRt^N$.  
An easy computation shows that the element $p_a\in P$ acts by multiplying  
a component of strength $S$ by $a^{S-2}$.  

We next define our main object of interest.  If $g$ is a metric defined in
a neighborhood of $0\in \R^n$, we construct a straight ambient metric in
normal form for $g$ as in Chapter~\ref{formalsect}.  We then evaluate the
covariant derivatives of curvature of the ambient metric at
$\rho =0$, $t=1$ as described in Chapter~\ref{cct}.  
If $n$ is odd, the values of 
all components of these covariant derivatives at the origin depend only the
derivatives of $g$ at the origin, while if $n$ is even, this is true for
components of strength at most $n+1$ by Proposition~\ref{indeterminacy}.   
If $g_{ij}(0) = h_{ij}$, then $\gt_{IJ} = \htt_{IJ}$ at $\rho = 0$, $t=1$, 
$x=0$.  In this case, 
Propositions~\ref{Tcontract} and  \ref{tracefree} show that the resulting
lists of components satisfy the relations of Definition~\ref{Rt}.  This
procedure therefore defines a map $c:\cN_c\rightarrow \cRt$ for $n$
odd, and $c:\cN_c\rightarrow \cRt^{n-3}$ for $n$ even.  Since 
the conformal curvature tensors are natural polynomial invariants of the
metric $g$, $c$ is a polynomial map.   
\begin{proposition}\label{equivariance}
The map $c: \cN_c\rightarrow \cRt$ (or $\cRt^{n-3}$ if $n$ is even) is
equivariant with respect to the $P$-actions. 
\end{proposition}
\begin{proof}
Recall that the action of $P$ on $\cN_c$ is given by 
$p.g = (\varphi^{-1})^*(\Omega^2g)$, where $\varphi$ and $\Omega$ are 
determined 
to map $g$ back to conformal normal form given the initial normalizations
defined by $p$; see the discussion preceding Proposition~\ref{orbits}
above.  By naturality of the conformal curvature tensors, 
$c((\varphi^{-1})^*(\Omega^2g))=(\varphi^{-1})^*(c(\Omega^2g))$,
where $(\varphi^{-1})^*$ on the right hand side is interpreted as the
pullback in the   
indices between $1$ and $n$ of each of the tensors in the list, leaving the
$0$ and $\infty$ indices alone.  
And $c(\Omega^2g)$ is given by 
Proposition~\ref{ctransform}.  We use these observations to check 
for each of the generating subgroups $O(h)$, $\R_+$ and $\R^n$ of $P$  
that $c(p.g) = p.(c(g))$, where the $P$-action on the right hand
side is that on
$\prod_{r=0}^{\infty}\left(\bigotimes^{4+r}\R^{n+2}{}^*
\otimes \sigma_{-2-r}\right)$.     

For $p=p_m$, we have $\Omega=1$ and $\varphi(x) = mx$, so  
$c(p_m.g)=(m^{-1})^*c(g)$ is obtained from $c(g)$ by transforming
covariantly under $O(h)$ the  
indices between $1$ and $n$.  But this is precisely how $p_m$ acts on  
$\prod_{r=0}^{\infty}\bigotimes^{4+r}\R^{n+2}{}^*$.

For $p=p_a$, we have $\varphi(x) = a^{-1}x$ and $\Omega =a^{-2}$.  By 
Proposition~\ref{ctransform}, the component $\Rt_{IJKL,M_1\cdots M_r}$ 
for $a^{-2} g$ is that for $g$ multiplied by $a^{2s_{\infty}-2}$.  Since
$(\varphi^{-1})^*$ acts by multiplying this component  
by $a^{s_M}$, it follows that the components of $c(p_a.g)$ are those
of $c(g)$ multiplied by $a^{s_M+2s_{\infty}-2} = a^{S-2}$, where 
$S=\|IJKLM_1\cdots M_r\|$.  But we noted above that 
this is precisely how $p_a$ acts on 
$\prod_{r=0}^{\infty}\left(\bigotimes^{4+r}\R^{n+2}{}^*
\otimes \sigma_{-2-r}\right)$.     

Finally, for $p=p_b$, we have $\varphi'(0)=I$ and $\Omega = 1-b_ix^i 
+O(|x|^2)$, so that the components of $c((\varphi^{-1})^*(\Omega^2g))$
are given 
by Proposition~\ref{ctransform} with $\Up_i=-b_i$.  But this is precisely
how $p_b$ acts on $\prod_{r=0}^{\infty}\bigotimes^{4+r}\R^{n+2}{}^*$. 
\end{proof}

Let us examine more carefully the equivariance of $c$ with respect to
the subgroup $\R_+\subset P$.  The
component $\Rt_{IJKL,M_r\cdots M_r}$ of $c(g)$ is a polynomial in the
components of the $g^{(s)}$, and this equivariance says that when the
$g^{(s)}$ 
are replaced by $a^sg^{(s)}$, then $\Rt_{IJKL,M_1\cdots M_r}$ is multiplied
by $a^{S-2}$ with $S=\|IJKLM_1\cdots M_r\|$.  In particular, 
$\Rt_{IJKL,M_1\cdots M_r}$ can only involve $g^{(s)}$ for $s\leq S-2$.   
This implies that for each $N\geq 0$ (satisying also $N\leq n-3$ for $n$
even), $c$ induces a $P$-equivariant polynomial map
$c^N:\cN_c^{N+2}\rightarrow \cRt^N$. 
Clearly these induced maps
satisfy the compatibility conditions
$c^{N-1}\pi_{N+2} = \widetilde{\pi}_Nc^N$, where   
$\pi_N:\cN_c^N\rightarrow \cN_c^{N-1}$ and 
$\widetilde{\pi}_N:\cRt^N\rightarrow \cRt^{N-1}$ are the natural
projections.    

The main result of this chapter is the following jet isomorphism theorem.
\begin{theorem}\label{jetiso}
Let $N\geq 0$ and assume that $N\leq n-3$ if $n$ is even.  Then 
$\cRt^N$ is a smooth submanifold of $\cTt^N$ whose
tangent space at $0$ is $T\cRt^N$, and the map $c^N:\cN_c^{N+2}\rightarrow   
\cRt^N$ is a $P$-equivariant polynomial equivalence. 
\end{theorem}
\noindent
For $n$ odd, it follows that $c:\cN_c\rightarrow \cRt$ is a $P$-equivariant
bijection since $c$ is the projective limit of the $c^N$.   

It will be convenient in the proof of Theorem~\ref{jetiso} to use
Theorem~\ref{riemjet} to realize $\cN_c$ in 
terms of curvature tensors on $\R^n$ rather than Taylor coefficients of  
metrics.  So we make the following definition.
Recall the space $\cR$ introduced in Definition~\ref{Rdef}.
\begin{definition}\label{Rcdef}
Define the space $\cR_c\subset \cR$ to be the subset consisting of lists of
tensors $(R^{(0)},R^{(1)},\ldots )$ satisfying in addition to the
conditions in Definition~\ref{Rdef} the following: for each $r\geq 0$,
\begin{equation}\label{Riccicond}
\operatorname{Sym}(h^{ik}R^{(r)}_{ijkl,m_1\cdots m_r}) =0.
\end{equation}
Here $\operatorname{Sym}$ refers to the symmetrization over the free 
indices $jlm_1\cdots m_r$.  For $N\geq 0$, by $\cR_c^N$ we will denote the
corresponding set of truncated lists $(R^{(0)},R^{(1)},\ldots, R^{(N)})$. 
\end{definition}
\noindent
For $r, N \geq 0$, we define also the following finite-dimensional vector 
spaces:

\bigskip
\noindent
$
\tau^r = \{R^{(r)} \in \bigwedge^2\R^n{}^*\otimes  
\bigwedge^2\R^n{}^*\otimes \bigotimes^r\R^n{}^*: \eqref{Riccicond} 
\mbox{ and (1), (2) of Definition~\ref{Rdef} hold}\}
$

\bigskip
\noindent
$
\sigma^r =\{R^{(r)}\in \tau^r: 
R_{ijkl,m_1\cdots m_r} = R_{ijkl,(m_1\cdots m_r)}\}
$

\bigskip
\noindent
$
\cT^N = \bigoplus_{r=0}^N\tau^r
$

\bigskip
\noindent
$
T\cR_c^N = \bigoplus_{r=0}^N\sigma^r\subset \cT^N
$

\bigskip

The bijection $\cN\rightarrow \cR$ asserted by Theorem~\ref{riemjet} 
clearly restricts to a bijection    
$\cN_c\rightarrow \cR_c$ whose truncated maps  
$\cN_c^{N+2}\rightarrow \cR_c^N$ are polynomial equivalences.     
By composition,   
we can regard $c$ and the $c^N$ as defined on the corresponding  
$\cR_c$ and $\cR_c^N$.  In the following, we will not have occasion to
refer to $\cN_c$ and $\cN_c^{N+2}$, so no confusion should arise from  
henceforth using the 
same symbols $c$ and $c^N$ for the maps defined on
$\cR_c$ and $\cR_c^N$.  We can transfer
the action of $P$ on $\cN_c$ to $\cR_c$.
The element $p_a$ acts on $\cR_c$ by multiplying $R^{(r)}$ by $a^{r+2}$;
this same prescription gives an $\R_+$ action on $\cT^N$ for $N\geq 0$.
The
$\R_+$-equivariance  of $c^N$ implies that a component
$\Rt_{IJKL,M_1\cdots M_r}$ of $c^N(R)$ with  
$\|IJKLM_1\cdots M_r\| = N+4$ can be written as a linear combination of
components 
of $R^{(N)}$ plus quadratic and higher terms in the components of the
$R^{(r)}$ with $r\leq N-2$.   

Our starting point for the proof of Theorem~\ref{jetiso} is the following
lemma.  
\begin{lemma}\label{Rnlinear}
For each $N\geq 0$, the subset $\cR_c^N\subset
\cT^N$ is a smooth submanifold whose tangent 
space at $0$ is $T\cR_c^N$. 
\end{lemma}
\begin{proof}
We will show that for each $N\geq 1$, there is a polynomial equivalence 
$\Phi^N: \cT^N \rightarrow \cT^N$ satisfying $d\Phi^N(0) = I$ 
and $\Phi^N(\cR_c^{N-1}\times \sigma^N)= \cR_c^N$.
Upon iterating this statement and using $\sigma^0 = \cR_c^0 = 
T\cR_c^0$, we conclude the existence of a polynomial
equivalence $:\cT^N \rightarrow \cT^N$ 
whose differential at $0$ is the identity, and which maps 
$T\cR_c^N \rightarrow \cR_c^N$.  The desired conclusion follows
immediately. 

When reformulated in terms of the spaces $\cR_c^N$, 
Lemma~\ref{extend} asserts 
the existence for each $N\geq 1$ of a polynomial map 
$\eta^N:\cT^{N-1}\rightarrow  \bigotimes^{N+4}\R^n{}^*$ with zero  
constant and linear term, such that the map $\Lambda^N:R\rightarrow  
(R,\eta^N(R))$ sends $\cR_c^{N-1}\rightarrow \cR_c^N$.  
Here $R$ denotes the list constituting an element of $\cT^{N-1}$.    
There is no loss
of generality in assuming that $\eta^N(\cT^{N-1})\subset \tau^N$ so that 
$\Lambda^N: \cT^{N-1}\rightarrow \cT^N$.  
Define $\Phi^N:\cT^N\rightarrow \cT^N$ by 
$\Phi^N(R,R^{(N)})=(R,R^{(N)} + \eta^N(R))$.  It is evident   
from the form of the relations defining $\cR_c^N$ that   
$\Phi^N(\cR_c^{N-1}\times \sigma^N)= \cR_c^N$, and clearly  
$(\Phi^N)^{-1}(R,R^{(N)})=(R,R^{(N)} - \eta^N(R))$.     
\end{proof}
\noindent
We remark that the same proof could have been carried out in terms of the
spaces of normal form coefficients, and shows that the subset 
$\cN_c^N\subset \cN^N$ is a smooth submanifold whose tangent space at $0$ 
is obtained by linearizing the equation obtained by writing
\eqref{Riccicond} in terms of the normal form coefficients 
$g_{ij,k_1\cdots k_r}$.

At this point we do not know  
that $\cRt^N$ is a submanifold of $\cTt^N$ with tangent space $T\cRt^N$, 
but it is clear that the tangent vector at $0$ to a smooth curve in
$\cRt^N$ must lie in $T\cRt^N$.  So we conclude for the differential of
$c^N$ at the origin that  $dc^N:T\cR_c^N\rightarrow T\cRt^N$.   
The differentiation of the action of $P$ on 
$\cR_c^N$ gives a linear action of $P$ on $T\cR_c^N$, and  
$dc^N:T\cR_c^N\rightarrow T\cRt^N$ is $P$-equivariant.  
By $\R_+$-equivariance and linearity of 
$dc^N$, it follows that $dc^N$ decomposes as a direct sum of maps   
$\sigma^M\rightarrow \st^M$ for $0\leq M \leq N$.  By the compatibility of 
the $c^N$ as $N$ varies, the map $\sigma^M\rightarrow \st^M$  
is independent of the choice of $N\geq M$, so we may as well denote it as 
$dc^N:\sigma^N\rightarrow \st^N$.  
The main algebraic fact on which 
rests the proof of Theorem~\ref{jetiso} is the following. 
\begin{proposition}\label{lineariso}
For $N\geq 0$ (and $N\leq n-3$ if $n$ is even),
$dc^N:\sigma^N\rightarrow \st^N$ is an isomorphism. 
\end{proposition}
\noindent
{\it Proof of Theorem~\ref{jetiso} using Proposition~\ref{lineariso}.}
We prove by induction on $N$ that there is a polynomial equivalence 
$\Pt^N:\cTt^N\rightarrow \cTt^N$ satisfying $d\Pt^N(0)=I$ and
$\Pt^N(\cRt^{N-1}\times \st^N)= \cRt^N$, and that 
$c^N:\cR_c^N\rightarrow \cRt^N$ is a polynomial equivalence.  
Just as in the proof of Lemma~\ref{Rnlinear}, iterating the 
first statement provides a polynomial equivalence 
$:\cTt^N\rightarrow \cTt^N$ whose differential at $0$ is the identity
and which maps $T\cRt^N\rightarrow \cRt^N$, from which follows the first
statement of Theorem~\ref{jetiso}.  

For $N=0$, we can take $\Pt^N$ to be the identity.  Since
$\cR_c^0 = T\cR_c^0 = \sigma^0$, 
$\cRt^0 = T\cRt^0 = \st^0$, and $c^0$ is linear and can be identified
with $dc^0$, the second 
statement is immediate from Proposition~\ref{lineariso}.   

Suppose for some $N\geq 1$ that we have the polynomial equivalence
$\Pt^{N-1}$ and we know that 
$c^{N-1}$ is a polynomial equivalence.  Recall the polynomial maps
$\eta^N:\cT^{N-1}\rightarrow \tau^N$, $\Lambda^N:\cT^{N-1}\rightarrow
\cT^N$ and $\Phi^N:\cT^N\rightarrow \cT^N$ 
constructed in
Lemma~\ref{extend} and Lemma~\ref{Rnlinear}.  
By the induction hypothesis that 
$c^{N-1}$ is a polynomial equivalence, we conclude that there is a
polynomial map  
$\Lt^N:\cTt^{N-1}\rightarrow \cTt^N$ such that 
$\Lt^N(\cRt^{N-1})\subset \cRt^N$ and such that the diagram
$$
\begin{CD}
\cR_c^{N-1}  @>\Lambda^N>>  \cR_c^N\\
@VV{c^{N-1}}V          @VV{c^N}V\\
\cRt^{N-1}   @>\Lt^N>>    \cRt^N
\end{CD}
$$
\noindent
commutes.  Using the compatibility of $c^{N-1}$ and $c^N$ and  
the form of the map $\Lambda^N$, one sees that $\Lt^N$ can be taken to 
have the form $\Lt^N(\Rt)=(\Rt,\ett^N(\Rt))$ where 
$\ett^N:\cTt^{N-1}\rightarrow \ttt^N$  has no
constant or linear terms.  Now define the map
$\Pt^{N}:\cTt^N\rightarrow \cTt^N$ by
$\Pt^N(\Rt,\Rt^{(N)})=(\Rt,\Rt^{(N)}+\ett^N(\Rt))$.  Clearly $d\Pt^N(0)=I$, 
and $\Pt^N(\cRt^{N-1}\times \st^N)= \cRt^N$ 
by the form of the relations defining $\cRt^{N}$.
It is a straightforward matter to check that the diagram 
$$
\begin{CD}
\cR_c^{N-1}\times \sigma^N  @>\Phi^N>>  \cR_c^N\\
@VV{c^{N-1}\times dc^N}V          @VV{c^N}V\\
\cRt^{N-1}\times \st^N   @>\Pt^N>>    \cRt^N
\end{CD}
$$

\bigskip
\noindent
commutes.  By the induction hypothesis and Proposition~\ref{lineariso}, the
vertical map on the left is a polynomial equivalence.  We conclude that 
$c^N$ is also a polynomial equivalence, completing the induction
step. \stopthm

\noindent
{\it Proof of Proposition~\ref{lineariso}.}  
The proof has two parts.  We will first construct an injective map 
$L:\st^N\rightarrow \sigma^N$.  We will then show that 
$dc^N:\sigma^N\rightarrow \st^N$ is injective.  These statements 
together imply that $\dim(\sigma^N) = \dim(\st^N)$, from which it then
follows that $dc^N$ is an isomorphsim.

Let $(\Rt_{IJKL,M_1\cdots M_r})\in \st^N$.  We can consider the components
$\Rt_{ijkl,m_1\cdots m_N}$ in which all the indices lie between $1$ and 
$n$.  This defines a map $L:\st^N\rightarrow \bigwedge^2\R^n{}^*\otimes  
\bigwedge^2\R^n{}^*\otimes \bigodot^N\R^n{}^*$, and clearly everything in
the range 
of $L$ satisfies conditions (1) and (2) of Definition~\ref{Rdef}.  We claim
that everything in the range of $L$ also satisfies \eqref{Riccicond}, so
that $L:\st^N\rightarrow \sigma^N$.  Condition (3) of Definition~\ref{Rt} 
implies $\htt^{IK}\Rt_{IjKl,m_1\cdots m_N}=0$.  This can be written as
$$
\Rt_{j\infty l0,m_1\cdots m_N}
+\Rt_{l\infty j0,m_1\cdots m_N}+h^{ik}\Rt_{ijkl,m_1\cdots m_N}=0.
$$
If we apply condition (5) of Definition~\ref{Rt} to
$\Rt_{j\infty l0,m_1\cdots m_N}$ and then symmetrize over 
$jlm_1\cdots m_N$, the result is $0$ by the skew symmetry of $\Rt^{(N)}$ in
the second pair of indices.  Similarly for  
$\Rt_{l\infty j0,m_1\cdots m_N}$.  It follows that the symmetrization of 
$h^{ik}\Rt_{ijkl,m_1\cdots m_N}$ vanishes.  This proves that
\eqref{Riccicond} holds. 

Next we show that $L$ is injective.  We claim that for
$(\Rt_{IJKL,M_1\cdots M_r})\in \st^N$, any component
$\Rt_{IJKL,M_1\cdots M_r}$ can be written as a linear combination of
components in which none of the indices $IJKLM_1\cdots M_r$ is $\infty$.  
We first show that any component can be written as a linear combination of
components in which none of $IJKL$ is $\infty$. To see this, note that 
\eqref{remove0} and $\htt^{AB}\Rt_{IJKA,BM_1\cdots M_r}=0$ imply
$$
-(r+2)\Rt_{IJK\infty,M_1\cdots M_r}=\Rt_{IJK\infty,0M_1\cdots M_r}= 
-\Rt_{IJK0,\infty M_1\cdots M_r}
-h^{ab}\Rt_{IJKa,bM_1\cdots M_r}.  
$$
Thus a component $\Rt_{IJKL,M_1\cdots M_r}$ in which $L=\infty$ can be
rewritten as a linear combination of components in which $L\neq \infty$ and
$IJK$ remain unchanged. 
Repeating this procedure allows the removal of any $\infty$'s in $IJKL$.
The same method allows the removal of $\infty$'s in $M_1\cdots M_N$: 
from \eqref{remove0} and $\htt^{AB}\Rt_{IJKL,ABM_2\cdots M_r}=0$, one has  
$$
-2(r+2)\Rt_{IJKL,\infty M_2\cdots M_r}=2\Rt_{IJKL,0\infty M_2\cdots M_r}=   
-h^{ab}\Rt_{IJKL,abM_2\cdots M_r}.
$$
Thus all $\infty$'s can be removed as indices.  Now \eqref{remove0} and 
(5) of Definition~\ref{Rt} can be used to remove any $0$'s as indices at
the expense of permuting the remaining indices between $1$ and $n$.  It
follows that any component $\Rt_{IJKL,M_1\cdots M_r}$ can be written as a
linear combination of components in which all indices are between $1$ and
$n$.  Thus $L$ is injective.

It remains to prove that $dc^N$ is injective.  If 
$R_{ijkl,m_1\cdots m_N}\in \sigma^N$, we set
$R_{jl,m_1\cdots m_N} = h^{ik}R_{ijkl,m_1\cdots m_N}$,
$R_{m_1\cdots m_N} = h^{jl}R_{jl,m_1\cdots m_N}$ and 
$$
(n-2)P_{jk,m_1\cdots m_N} = R_{jk,m_1\cdots m_N}
-R_{m_1\cdots m_N}h_{jk}/2(n-1). 
$$
We also denote by 
$W_{ijkl,m_1\cdots m_N}\in \bigwedge^2\R^n{}^*\otimes
\bigwedge^2\R^n{}^*\otimes \bigodot^N\R^n{}^*$ the tensor obtained by
taking the trace-free part in $ijkl$ while ignoring $m_1\cdots m_N$:    
\begin{equation}\label{decomp}
W_{ijkl,m_1\cdots m_N} = R_{ijkl,m_1\cdots m_N}
-2h_{l[j}P_{i]k,m_1\cdots m_N}h_{jl} +2h_{k[j}.P_{i]l,m_1\cdots m_N}.
\end{equation}
Then $W_{ijkl,m_1\cdots m_N}$ satisfies (1) of Definition~\ref{Rdef} but not
necessarily (2).  We also define
$$
C_{jkl,m_1\cdots m_{N-1}}= 2P_{j[k,l]m_1\cdots m_{N-1}}.
$$
Contracting the second Bianchi identity (2) of Definition~\ref{Rdef} in the
usual way shows that $W^i{}_{jkl,im_1\cdots m_{N-1}}=(3-n)C_{jkl,m_1\cdots 
m_{N-1}}$  and $P^i{}_{j,im_1\cdots m_{N-1}}= P^i{}_{i,jm_1\cdots
m_{N-1}}$.  
\begin{lemma}\label{Weyl}
Let $R_{ijkl,m_1\cdots m_N} \in \sigma^N$.  If $n\geq 4$ and
$W_{ijkl,m_1\cdots m_N}=0$, then $R_{ijkl,m_1\cdots m_N}=0$. 
If $n=3$ and $C_{jkl,m_1\cdots m_{N-1}}= 0$, then 
$R_{ijkl,m_1\cdots m_N}=0$. 
\end{lemma}
\begin{proof}
If $N=0$, this follows from the decomposition of the curvature tensor into 
its Weyl piece and its Ricci piece.  Suppose $N\geq 1$.  If $n\geq 4$, the
contracted Bianchi identity above shows that $C_{jkl,m_1\cdots m_{N-1}}=
0$, which is our hypothesis if $n=3$.  
(The hypothesis $W_{ijkl,m_1\cdots m_N}=0$ for $n\geq 4$ is of course
automatic for $n=3$.)  Thus we conclude for any 
$n$ that $P_{j[k,l]m_1\cdots m_{N-1}}=0$, so 
$P_{ij,m_1\cdots m_N}= P_{(ij,m_1\cdots m_N)}$.  Since 
$R_{ijkl,m_1\cdots m_N} \in \sigma^N$, we also have $R_{(ij,m_1\cdots
  m_N)}=0$.  Therefore
$$
(n-2)P_{ij,m_1\cdots m_N} = (n-2)P_{(ij,m_1\cdots m_N)} =  
-R_{(m_1\cdots m_N}h_{ij)}/2(n-1).
$$
Now $h^{ij}P_{(ij,m_1\cdots m_N)} =  h^{ij}P_{ij,m_1\cdots m_N} =  
R_{m_1\cdots m_N}/2(n-1)$.  Hence the symmetric tensor 
$P=P_{(ij,m_1\cdots m_N)}$ is in the kernel of the operator
$P\rightarrow (n-2)P + \operatorname{Sym}({\operatorname{tr}}(P)h)$.  
This operator is injective on symmetric tensors, so we conclude that 
$P_{ij,m_1\cdots m_N}=0$.  The conclusion now follows from \eqref{decomp}.    
\end{proof}

We will prove that $dc^N$ is injective by showing that if
$R_{ijkl,m_1\cdots m_N}\in \ker(dc^N)$, then 
the hypotheses of Lemma~\ref{Weyl} hold. 
To get the flavor of the argument,
consider first the cases $N=0,1$.  Now $\sigma^0$ is the space of trace-free
curvature tensors.  In Chapter~\ref{cct}, 
we found that $\Rt_{ijkl} = W_{ijkl}$
is the Weyl piece of such a curvature tensor.  So $c^0$ is linear and
is obviously injective.  In Chapter~\ref{cct}, 
we also calculated the curvature 
components $\Rt_{\infty jkl} = C_{jkl}$ and $\Rt_{ijkl,m}$ (see
\eqref{Vform}).  We see that $c^1$ is also linear, so can be identified
with $dc^1$.  If $R_{ijkl,m}\in \ker(dc^1)$, we first conclude
by considering $\Rt_{\infty jkl}$ that $C_{jkl}=0$, and then 
by considering $\Rt_{ijkl,m}$ that $W_{ijkl,m}=0$ (for $n\geq 4$), as
desired.  

For the general case we need to understand the relation between
covariant derivatives 
with respect to the ambient metric and covariant derivatives with respect
to a representative $g$ on $M$.  Recall that the conformal curvature
tensors 
are tensors on $M$ defined by evaluating components of $\Rt_{IJKL,M_1\cdots
M_r}$ at $\rho=0$ and $t=1$.  We can take further covariant derivatives   
of such a tensor with respect to $g$.  We will denote by 
$\Rt_{IJKL,M_1\cdots M_r|p_1\cdots p_s}$  the tensor on $M$ obtained by
such further covariant differentiations.  For example, 
$\Rt_{ijkl,|m}= W_{ijkl,m}$, whereas $\Rt_{ijkl,m}$ is the tensor 
$V_{ijklm}$ given by \eqref{Vform}.  An inspection of \eqref{cnr} shows the
relation between $\Rt_{IJKL,M_1\cdots M_rp}$ and
$\Rt_{IJKL,M_1\cdots M_r|p}$.  Recalling that $g_{ij}' =2P_{ij}$ at $\rho
=0$, one sees that $\Rt_{IJKL,M_1\cdots M_rp}-\Rt_{IJKL,M_1\cdots M_r|p}$ 
is a linear combination of components of $\tilde{\nabla}^r\Rt$, possibly
multiplied 
by components of $g$, plus quadratic terms in curvature.  Iterating, it
follows that
$\Rt_{IJKL,M_1\cdots M_rp_1\cdots p_s}-\Rt_{IJKL,M_1\cdots M_r|p_1\cdots  
p_s}$ is a linear combination of terms of the form 
$\Rt_{ABCD,F_1\cdots F_u|q_1\cdots q_t}$ with $t<s$, possibly multiplied by 
components of $g$, plus nonlinear terms in curvature.    

We return now to consideration of $dc^N:\sigma^N\rightarrow \st^N$. 
The symbol $\Rt_{IJKL,M_1\cdots M_r}$ 
with $\|IJKLM_1\cdots M_r\|=N+4$ is now to be interpreted  
as the linear function of the $R_{ijkl,m_1\cdots m_N}$ 
obtained by applying $dc^N$.  Similarly, we now interpret the  
symbol $\Rt_{IJKL,M_1\cdots M_r|p_1\cdots p_s}$ for 
$\|IJKLM_1\cdots M_r\|+s=N+4$ as a linear function of the 
$R_{ijkl,m_1\cdots m_N}$.  Suppose that 
$R_{ijkl,m_1\cdots m_N}\in \ker(dc^N)$.
We claim that $\Rt_{IJKL,M_1\cdots M_r|p_1\cdots p_s}=0$ for 
$\|IJKLM_1\cdots M_r\|+s=N+4$.  The proof is by induction on 
$s$.  For $s=0$, this is just the hypothesis that  
$R_{ijkl,m_1\cdots m_N}\in \ker(dc^N)$.  Since we are considering the
linearization $dc^N$, the quadratic terms may be ignored in the
relation derived in 
the previous paragraph between ambient covariant derivatives and covariant
derivatives on $M$.  A moment's thought shows that this relation 
provides the induction step to increase $s$ by $1$.  

Taking $r=0$, we conclude 
that $\Rt_{IJKL,|m_1\cdots m_s}=0$ for $\|IJKL\|+s=N+4$.  For 
$IJKL = ijkl$ we obtain $W_{ijkl,m_1\cdots m_N}=0$ and 
for $IJKL=\infty jkl$ we obtain $C_{jkl,m_1\cdots m_{N-1}}=0$.    
Lemma~\ref{Weyl} then shows that $R_{ijkl,m_1\cdots m_N}=0$ as desired.   
\stopthm

\section{Scalar Invariants}\label{sinv}
The jet isomorphism theorem~\ref{jetiso} reduces the study of local 
invariants of conformal structures to the study of $P$-invariants of 
$\cRt$ (we must of course impose the usual finite-order truncation for $n$  
even).  
An invariant theory for scalar $P$-invariants of $T\cRt$ 
was developed in \cite{BEGr}.  In this chapter we 
show how to derive a characterization of scalar invariants of conformal
structures by reduction to the relevant results of \cite{BEGr}.

Recall that a scalar 
invariant $I(g)$ of metrics of signature $(p,q)$ is a polynomial in the
variables $(\partial^{\alpha}g_{ij})_{|\alpha| \geq 0}$
and $|\det{g_{ij}}|^{-1/2}$, which is coordinate-free in the sense that its  
value is independent of orientation-preserving changes of the coordinates
used to express and differentiate $g$.  Such a scalar invariant of metrics 
is said to be {\it even} if it is also unchanged under
orientation-reversing changes of coordinates, and {\it odd} if it changes
sign under orientation-reversing coordinate changes. 

It follows from the jet isomorphism theorem for
pseudo-Riemannian geometry and Weyl's classical invariant theory 
that every scalar invariant of metrics is a linear combination of 
complete contractions 
\begin{gather}\label{comcontr}
\begin{gathered}
\contr(\nabla^{r_1}R\otimes \cdots \otimes\nabla^{r_L}R)\\
\contr(\mu\otimes\nabla^{r_1}R\otimes \cdots \otimes\nabla^{r_L}R),
\end{gathered}
\end{gather}
where $\mu_{i_1\cdots i_n} = |\det{g}|^{1/2}\varepsilon_{i_1\cdots i_n}$ is 
the volume form with respect to a chosen orientation and the contractions  
are with respect to $g$.  Here $\varepsilon_{i_1\cdots i_n}$ denotes the 
sign of the permutation.  Complete contractions of the first type are even
and the second type are odd.  Such an invariant of metrics is said to be     
conformally invariant 
of weight $w$ if $I(\Omega^2 g) = \Omega^{w}I(g)$ for smooth positive
functions $\Omega$.  Under a constant rescaling $g\rightarrow a^2g$, 
the complete contractions above are multiplied by $a^{-2L-\sum r_i}$.   
Since the total number of contracted indices must be even, it follows that 
the weight of a nonzero even scalar conformal invariant must be a negative 
even integer, and that of a nonzero odd scalar conformal invariant must be
a negative integer which agrees with $n \mod 2$.    

Scalar conformal invariants can be constructed quite simply via the ambient
metric.  Consider complete contractions
\begin{gather}\label{ambcontr}
\begin{gathered}
\contr(\nt^{r_1}\Rt\otimes \cdots \otimes \nt^{r_L}\Rt)\\
\contr(\mut\otimes \nt^{r_1}\Rt\otimes \cdots \otimes \nt^{r_L}\Rt)\\ 
\contr(\mut_0\otimes \nt^{r_1}\Rt\otimes \cdots \otimes \nt^{r_L}\Rt), 
\end{gathered}
\end{gather}
where now $\Rt$ is the curvature tensor of a straight ambient metric $\gt$
for $[g]$, $\mut$ denotes the volume form of $\gt$ with respect to an
orientation on $\cGt$ induced from a choice of orientation on $M$,
$\mut_0 = T\into \mut$, 
and the contractions are taken with respect to $\gt$.     
When evaluated for a specific ambient metric, these  
contractions define functions on $\cGt$.  The homogeneity of $\gt$ and $T$
imply that the 
functions defined by contractions of the first two types 
are homogeneous of degree $-2L-\sum r_i$ with respect to the dilations 
$\delta_s$, while those of the third type are homogeneous of degree
$1-2L-\sum r_i$.  Their restrictions to $\cG\subset \cGt$  
are independent of the diffeomorphism ambiguity of $\gt$ (we restrict 
to orientation-preserving diffeomorphisms for the second and third types) 
since the diffeomorphism restricts to the identity on $\cG$.
If $n$ is odd, then the restriction of any such contraction to $\cG$ is 
also clearly independent of the infinite-order ambiguity in the ambient
metric, so depends only on $[g]$.  A representative metric $g$ defines a
section of $\cG$, so composing with this section gives a function $I(g)$ on
$M$.  The homogeneity of the contraction as a function on $\cG$ implies
that $I(\Omega^2 g) = \Omega^{w}I(g)$, where $-w = 2L+\sum r_i$ for the
first two and $-w = 2L+\sum r_i-1$ for the third.  If we take   
$\gt$ to be in normal form relative to $g$ and recall the discussion of 
conformal curvature tensors in Chapter~\ref{cct}, it follows that in local
coordinates $I(g)$ does indeed have the required polynomial dependence on
the Taylor 
coefficients of $g$.  Thus $I(g)$ is a scalar conformal invariant of weight 
$w$.  This proves the following proposition in case $n$ is odd.
\begin{proposition}\label{ambcontrprop}
If $n$ is odd, the complete contractions \eqref{ambcontr} define scalar
conformal invariants.  The first is even and has 
$-w=2L+\sum r_i$; 
the second and third are odd and have $-w=2L+\sum r_i$  
and $-w=2L+ \sum r_i -1$.  If $n$ is even, the same
statements are true with the restrictions 
$L\geq 2$ and $-w\leq n+2$ for the first contraction, and 
$-w\leq 2n-2$ for the last two.   
\end{proposition}

By a Weyl conformal invariant of metrics, we will mean a linear combination
of  
complete contractions \eqref{ambcontr}, all of which have the same weight,
and which satisfy the restrictions of Proposition~\ref{ambcontrprop} if $n$ 
is even.   
Every Weyl invariant can be written as a sum of an even Weyl invariant and
an odd Weyl invariant.

Before giving the proof of Proposition~\ref{ambcontrprop} for $n$ even, we
make some observations concerning 
odd invariants.  The Bianchi and Ricci identities imply that the
skew-symmetrization over any three indices of $\nabla^r R$ can be written
as a quadratic expression in the $\nabla^{r'}R$ with $r'\leq r-2$.  An odd
contraction in \eqref{comcontr} with $L<n/2$ necessarily has at least three 
of the indices of at least one of the $\nabla^{r_i} R$ contracted against
indices in $\mu$.  So by induction, it follows 
that an odd contraction in \eqref{comcontr}  
can be written as a linear combination of contractions of
the same form but with $L\geq n/2$.  Similarly, contractions of the second
and third types in \eqref{ambcontr} can be written as linear combinations
of contractions of the same type and the same weight, but with 
$L\geq (n+2)/2$ for the second type and $L\geq (n+1)/2$ for the third type.   
This leads to the following theorem.
\begin{theorem}\label{vanishing}
There are no nonzero odd scalar conformal
invariants with $-w<n$.  There are no nonzero odd Weyl
invariants with $-w<n+1$.      
\end{theorem}
\begin{proof}
The first statement is clear from the above observations since $-w=2L+\sum 
r_i$.  The same reasoning 
shows that there are no nonzero contractions of the second type in 
\eqref{ambcontr} with $-w<n+2$.  For contractions of the third
type, one has $-w = 2L+\sum r_i-1$, which is $<n+1$ only if 
$L=(n+1)/2$ (so $n$ must be odd) and each $r_i=0$.  However, any such  
contraction must vanish.  To see this, it suffices to take $\gt$ to be
in normal form.  Since one of the indices of $\mut_0$ is  
'$\infty$' and is contracted against one of the $\Rt$ factors, the result
follows from the facts that $\Rt_{0JKL}=0$ and $\gt^{I\infty}=0$ at $\rho
=0$ unless $I=0$.    
\end{proof}

\medskip
\noindent
{\it Proof of Proposition~\ref{ambcontrprop} for $n$ even.}   
The discussion in the $n$ odd case is valid also for $n$ even so long as we
show that the 
hypotheses guarantee that the restriction of \eqref{ambcontr} to $\cG$   
is independent of the $O^+_{IJ}(\rho^{n/2})$ ambiguity in $\gt$.   
In showing this, we may assume by Proposition~\ref{normalform} that $\gt$
is in normal form relative to a representative metric $g$.  

Consider the expansion of the first contraction of \eqref{ambcontr}
into components at $\rho =0$.  By the normal form assumption, we have that 
$\gt^{IJ}=0$ at $\rho=0$ unless $\|IJ\|=2$.  Therefore each pair of
contracted indices must have strength 2 to contribute.  
The total number of indices being contracted is  
$4L+\sum r_i$.  It follows that if $S_1,\ldots, S_L$  
are the strengths of the factors occurring in a contributing monomial, then 
$\sum S_i = 4L+\sum r_i = 2L -w \leq 2L + n+2$.  As we have observed 
previously,  
Proposition~\ref{Tcontract} implies that for any $r$, a component of
$\nt^r\Rt$ 
vanishes unless its strength is at least 4.  So we may assume that 
$S_i\geq 4$ for each $i$.  Thus for any $i_0$ we have
$$
S_{i_0}+4(L-1)\leq S_{i_0} + \sum_{i\neq i_0}S_i= \sum S_i \leq 2L+n+2.
$$
Therefore $S_{i_0}\leq n+6-2L$ with strict inequality if $S_i>4$ for some
$i\neq i_0$.  This implies that $S_{i_0}\leq n$ if $L\geq 3$, and that
$S_{i_0}\leq n+1$ if $L=2$ and $S_i>4$ for some $i \neq i_0$.  
Proposition~\ref{indeterminacy}  
guarantees in these cases that the first contraction of \eqref{ambcontr}  
is independent of the ambiguity in $\gt$ and defines a scalar conformal
invariant.  

It remains to consider the case where $L=2$ and $S_i = 4$ for some $i$.  We
may relabel the factors in \eqref{ambcontr} if necessary so that $S_1=4$.   
As noted above, Proposition~\ref{indeterminacy} implies the result if
$S_2\leq n+1$, so we are reduced to consideration of the case $S_1=4$, 
$S_2 = n+2$.  The only nonvanishing component with $S_1 = 4$ is 
$\Rt_{ijkl}$, which equals $W_{ijkl}$ when evaluated at $\rho =0$, $t=1$. 
This must be contracted with respect to $g$ with a component of strength  
$n+2$.  This component of strength $n+2$ must therefore have exactly four 
indices between $1$ and $n$; all other indices must be $0$ or $\infty$.  
Using \eqref{Tcontract}, we may rewrite such a component as a linear
combination of components in which $0$ does not occur as an index.  Such a 
component has exactly $4$ indices between $1$ and $n$ and exactly  
$(n-2)/2$ indices which are $\infty$.  We investigate the dependence of
such components on the $O(\rho^{n/2})$ ambiguity in  
the expansion of the $g_{ij}(x,\rho)$ coefficient in the ambient metric 
\eqref{spform}.  

Consider the formula for a covariant derivative of curvature
$\Rt_{IJKL,M_1\cdots M_r}$ of a general pseudo-Riemannian metric $\gt_{IJ}$
in terms of coordinate derivatives of $\gt_{IJ}$.  Beginning with 
$$
\Rt_{IJKL} =  \frac12 
(\pa^2_{IL}\gt_{JK}+\pa^2_{JK}\gt_{IL}-\pa^2_{JL}\gt_{IK}-\pa^2_{IK}\gt_{JL}) 
 +(\Gat_{IL}^Q\Gat_{JKQ}-\Gat_{IK}^Q\Gat_{JLQ})
$$
and successively covariantly differentiating, one sees that 
$\Rt_{IJKL,M_1\cdots M_r} = \text{I} + \text{II} + \text{III}$, where 
$\text{I}$ is a linear combination
of terms of the form $\pa^{r+2}_{P_1\cdots P_{r+2}}\gt_{AB}$ in which  
$ABP_1\cdots P_{r+2}$ is a permutation of $IJKLM_1\cdots M_r$, 
$\text{II}$ is a linear combination
of terms of the form $\Gat_{CD}^Q\pa^{r+1}_{QP_1\cdots P_r}\gt_{AB}$ 
and $\Gat_{CD}^Q\pa^{r+1}_{BP_1\cdots P_r}\gt_{AQ}$ in which 
$ABCDP_1\cdots P_r$ is a permutation of $IJKLM_1\cdots M_r$, and 
$\text{III}$ involves only $\pa^s\gt$ with $s\leq r$.     

Apply this observation to a curvature 
component $\Rt_{IJKL,M_1\cdots M_{(n-2)/2}}$ of the ambient metric,  
where $IJKLM_1\cdots M_{(n-2)/2}$ is a permutation of     
$ijkl\underbrace{\infty \cdots \infty}_{(n-2)/2}$.  When restricted to 
$\rho =0$, the terms in I and III
involve only $\pa^r_{\rho}g(x,0)$ for $r\leq n/2-1$.  The terms in 
II may have a factor of $\pa^{n/2}_\rho g(x,0)$, but any such factor is
multiplied by $\Gat_{ab}^\infty$, where $ab$ are two of $ijkl$.  By
\eqref{cnr}, we have $\Gat_{ab}^\infty|_{\rho = 0} = -g_{ab}$.  
Consequently, such terms in II drop out when contracted against $W_{ijkl}$.
It follows that the contraction \eqref{ambcontr} is independent of
the $O(\rho^{n/2})$ ambiguity as desired.     

The argument for the odd contractions is similar to that of the second
paragraph above.  We give the details for
contractions of the second type.  As noted previously, we may assume that 
$L\geq (n+2)/2$.  The total number of indices being contracted is now
$4L + \sum r_i +n+2$.  If $S_i$ are the strengths of the factors in a
contributing monomial in the components of the $\nt^r\Rt$, then since
the strengths of the indices of $\mut$ sum to $n+2$, we again have   
$\sum S_i =4L + \sum r_i = 2L-w$, and now this is $\leq 2L+2n-2$.  As
in the argument above, this implies that 
$S_{i_0}\leq 2n+2-2L \leq n$, so Proposition~\ref{indeterminacy} implies
that the contraction is independent of the 
ambiguity in $\gt$.
\stopthm

We remark that one can also consider the dependence on the ambiguity  
in $\gt$ of the first contraction in \eqref{ambcontr} in case
$L=1$.  Modulo contractions with $L\geq  
2$, the only possibilities are $\tilde{\Delta}^{r/2} \tilde{S}$, where 
$\tilde{S}$ denotes 
the scalar curvature and $\tilde{\Delta}$ the $\gt$-Laplacian.  
It is not hard to see that these are
independent of the $O_{IJ}^+(\rho^{n/2})$ ambiguity in $\gt_{IJ}$ as long
as $r\leq n-2$, which corresponds to $-w\leq n$.  However, the resulting 
conformal invariants all vanish 
since $\Rt_{IJ}= O_{IJ}^+(\rho^{n/2-1})$.   

For an example, consider a complete contraction \eqref{ambcontr} in which 
$r_i=0$ for all 
$i$.  Since $\Rt_{IJK0}=0$, it follows that the resulting conformal
invariant is $\contr(W\otimes \cdots \otimes W)$ with the same pairing of
the indices, where $W$ denotes the
Weyl tensor and the contractions are now with respect to $g$.  Of course,
this is invariant with no restrictions on $L$ when $n$ is even.   

A more interesting example is $\|\nt\Rt\|^2$.
Proposition~\ref{ambcontrprop} implies that this defines a conformal
invariant with $-w=6$ in all dimensions $n\geq 3$.  Expanding the 
contraction and evaluating the components shows that the conformal 
invariant is 
\begin{equation}\label{normsquared}
\|V\|^2 + 16(W,U) + 16\|C\|^2,
\end{equation}
where $V_{ijklm}$ is given by \eqref{Vform} and 
$$
U_{mjkl}= C_{jkl,m}-P_m{}^iW_{ijkl}.
$$
When $n=3$ this reduces to a nonzero multiple of $\|C\|^2$.  
(For $n=4$, the expression involving $B_{ij}/(n-4)$ in the formula for
$\Rt_{\nf jkl,m}|_{\rho =0,\,t=1} = Y_{jklm}$ is replaced by an expression
involving $g_{ij}''$, but this 
drops out as guaranteed by Proposition~\ref{ambcontrprop}.)

We are interested in the question of the extent to which all scalar
conformal invariants are Weyl invariants.  The results presented here 
resolve this question for all invariants when $n$ is odd and for invariants  
with $-w\leq n$ when $n$ is even.  A scalar conformal invariant
is said to be {\it exceptional} if it is not a Weyl invariant. 
Theorem~\ref{vanishing} shows that any nonzero odd scalar conformal  
invariant with $-w=n$ is exceptional.  These are classified as follows.  

\begin{theorem}\label{oddlowweight}
If $n\not\equiv 0 \mod 4$, there are no nonzero odd scalar conformal  
invariants with $-w=n$.   

If $n\equiv 0 \mod 4$, the space of odd scalar conformal invariants 
with $-w=n$ has dimension $p(n/4)$, the number of partitions of $n/4$.   
Every nonzero such invariant is exceptional.    
A basis for this space may be taken to be the 
Pontrjagin invariants whose integrals give the Pontrjagin numbers of a
compact oriented $n$-dimensional manifold.  
\end{theorem}
\begin{proof} 
The proof does not use conformal invariance; these are facts about the
space spanned by odd contractions \eqref{comcontr} with $2L+\sum r_i =n$.  
As noted above, we may restrict consideration to contractions with 
$L\geq n/2$.  We must therefore have 
$L=n/2$ and $r_i = 0$ for each $i$.  In particular, $n$ must be even for
such a contraction to be nonzero.  Upon decomposing $R$ into its Weyl and
Ricci pieces, the fact that at most two indices of $\mu$ can be contracted
against any Weyl factor implies that the Ricci
curvature cannot appear in any nonzero contraction.
Therefore, an odd invariant with $-w=n$ must be a linear 
combination of contractions of the form 
$\contr(\mu\otimes \underbrace{W \otimes \cdots \otimes W}_{n/2})$, 
where $W$ is the Weyl tensor.  Any such contraction is clearly conformally
invariant.  It is easily seen using the symmetries of
$W$ that all such contractions vanish if $n\equiv 2 \mod 4$ and that there
are at most $p(n/4)$ linearly independent such contractions if $n\equiv 0
\mod 4$.  The Pontrjagin invariants are $p(n/4)$ linearly independent
invariants which are linear combinations of contractions of this type.  
It follows that the dimension of the space of such invariants is 
$p(n/4)$ and that the Pontrjagin invariants form a basis.  
\end{proof}

The main result of this chapter is the following.
\begin{theorem}\label{invariants}
If $n$ is odd, every scalar conformal invariant is a Weyl invariant. 
If $n$ is even, every even scalar conformal invariant with 
$-w\leq n$ is a Weyl invariant.   
\end{theorem}
\noindent
For $n$ even, an extension of Theorems~\ref{oddlowweight} and
\ref{invariants} covering all 
invariants is described in \cite{GrH2}.  A theory of even invariants of a   
conformal structure coupled with a conformal density is given in
\cite{Al1}, \cite{Al2}
which recovers Theorem~\ref{invariants} when applied to even invariants 
depending only on the conformal structure.  

As an illustration of the use of Theorem~\ref{invariants}, consider 
even invariants with $-w=4$ or $6$.  For $-w=4$, the only possible   
even contraction \eqref{ambcontr} is $\|\Rt\|^2$, so all even scalar
conformal invariants with $-w=4$ are multiples of $\|W\|^2$.  For even
invariants with $-w=6$, it is  
easily seen that modulo contractions with $L=3$, the only possibility with
$L=2$ is $\|\nt\Rt\|^2$.  It follows that all even scalar conformal
invariants with $-w=6$ (in dimension $n\neq 4$) are linear combinations 
of \eqref{normsquared} and cubic terms in the Weyl tensor.  In general, 
for a given weight, work is required to determine which contractions are
independent, but Theorem~\ref{invariants} gives a finite list of spanning
conformal invariants.    

The first step in the proof of Theorem~\ref{invariants} is to reformulate 
scalar conformal invariants in terms of invariants of the $P$-action on the
space of conformal normal forms.  In general, by a scalar invariant of a
group acting on a space is meant a scalar function on the space
transforming by a  
character of the group.  Accordingly, if $\cS$ is a $P$-space, a real
valued function $Q$ on $\cS$ is said to be an even (resp. odd)
$P$-invariant of weight $w$ if $Q(p.s) = \sigma_w(p)Q(s)$
(resp. $\det(p)\sigma_w(p)Q(s)$) for $s\in \cS$, $p\in P$.   
Equivalently,
$Q$ defines a $P$-equivariant map $Q:\cS\rightarrow \sigma_w$ (resp. 
$\det\otimes \sigma_w$).    
A function is said to be a $P$-invariant of $\cS$ (or a $P$-invariant
function on $\cS$) if it is a sum of  
an even and an odd invariant of weight $w$ for some $w$.     
\begin{proposition}\label{sci}
There is a one-to-one correspondence between scalar conformal invariants of
weight $w$ and $P$-invariant polynomials $Q$ on $\cN_c$ of weight $w$.
\end{proposition}
\begin{proof}
A scalar conformal invariant $I(g)$ of weight $w$ defines a polynomial $Q$
on $\cN_c$ by evaluation at $0\in \R^n$.  We can write $I = I_+ + I_-$,
where $I_+$ (resp. $I_-$) are even (resp. odd) scalar conformal invariants
of weight $w$.  This gives $Q=Q_++Q_-$ for the polynomials on
$\cN_c$.  The definition 
of the action of $P$ on $\cN_c$ and the invariance of $I_+$, $I_-$ show
immediately that $Q_+$ (resp. $Q_-$) is an even (resp. odd) invariant of
$\cN_c$ of weight $w$.     

For the other direction, let $Q$ be a polynomial invariant of $\cN_c$ of 
weight $w$.  If $g$ is any metric defined near the  
origin in $\R^n$, by Proposition~\ref{cnf} and the existence of geodesic
normal coordinates, we can find $\Omega$ with $\Omega(0)=1$ and 
$\varphi$ with $\det\varphi'(0)>0$ such that  
$(\varphi^{-1})^*(\Omega^2g)\in \cN_c$.  Define $I(g) = 
Q((\varphi^{-1})^*(\Omega^2g))$.  The element
$(\varphi^{-1})^*(\Omega^2g)$ of $\cN_c$ is not uniquely determined, 
but Proposition~\ref{orbits} and the normalizations $\Omega(0)=1$, 
$\det\varphi'(0)>0$ imply that
it is determined up to the action of an element $p\in P$
such that $a=1$ and $\det(p)=1$.  The assumed invariance of $Q$  
therefore shows that $I(g)$ is well-defined.  This well-definedness 
makes it clear that $I((\varphi^{-1})^*(\Omega^2g))=I(g)$ for
general $\varphi$, $\Omega$ such that $\det\varphi'(0)>0$, $\Omega(0)=1$.     
The transformation law $Q(p_{a^{-1}}.g) = a^{w}Q(g)$ for $a>0$ and $g\in
\cN_c$ implies that $I(a^2g) = a^wI(g)$ for $a>0$ and general $g$. 
Therefore $I$ satisfies the conformal transformation law
$I(\Omega^2g)=\Omega(0)^wI(g)$.   

For $g\in \cN$, the $\Omega$ and $\varphi$ are uniquely determined to
infinite order if we
require also that $d\Omega(0)=0$ and $\varphi'(0) = Id$, and the
construction of the conformal normal form shows that the Taylor
coefficients  
of such $\Omega$ and $\varphi$ depend polynomially on $g\in \cN$.
Therefore $I(g)$ defines an $SO(h)$-invariant polynomial on $\cN$.
By the pseudo-Riemannian jet isomorphism theorem and   
Weyl's invariant theory, an $SO(h)$-invariant polynomial on $\cN$ is a 
linear combination 
of complete contractions \eqref{comcontr}, so that $I(g)$ is a scalar
invariant of metrics.  Combining this 
with the conformal transformation law established above, we deduce that $I$
is a scalar conformal invariant as desired.

The maps $I\rightarrow Q$ and 
$Q\rightarrow I$ are easily seen to be inverses of one another.
\end{proof} 

We will use Theorem~\ref{jetiso} to transfer $P$-invariant polynomials from 
$\cN_c$ to $\cRt$.  The next Lemma will assure when $n$ is even that we are
in the range in which Theorem~\ref{jetiso} applies.   
\begin{lemma}\label{weightrestrict}
Let $I(g)$ be an even scalar conformal invariant of weight $w$.
The associated polynomial $Q(g)$ on $\cN_c$ determined by
Proposition~\ref{sci}  
can be written as the restriction to $\cN_c$ of a polynomial in derivatives
of $g$ of order $\leq -w-2$.
\end{lemma}
\begin{proof}
The polynomial $Q(g)$ determined by Proposition~\ref{sci} agrees with the
restriction to $\cN_c$ of a linear combination of even complete
contractions \eqref{comcontr} with 
$2L+\sum r_i=-w$.  Thus for each $i_0$ we have 
$r_{i_0}\leq \sum r_i = -w-2L$.  If
$L\geq 2$, we obtain $r_{i_0}\leq -w-4$, so the contraction involves at
most $-w-2$ derivatives of $g$ as claimed.  
Modulo contractions with $L\geq 2$, the only possibility with $L=1$ is
$\Delta^{(-w-2)/2}S$, where $S$ denotes the scalar curvature.  But
contracting $\operatorname{Sym}(\nabla^r\Ric)(g)(0)=0$ and using the
Bianchi and Ricci identities shows that 
$\Delta^{(-w-2)/2}S$ agrees on $\cN_c$ with a linear combination of
complete contractions with $L\geq 2$.   
\end{proof}

Let $\et$ denote the usual volume form on $\R^{n+2}$ and  
$\et_0 = e_0\into \et$.  Complete contractions with respect to $\htt$  
of the form    
\begin{gather}\label{ambcontrR}
\begin{gathered}
\contr(\Rt^{(r_1)}\otimes \cdots \otimes \Rt^{(r_L)})\\ 
\contr(\et\otimes \Rt^{(r_1)}\otimes \cdots \otimes \Rt^{(r_L)})\\  
\contr(\et_0\otimes \Rt^{(r_1)}\otimes \cdots \otimes \Rt^{(r_L)})  
\end{gathered}
\end{gather}
define $P$-invariant polynomials on $\cTt$ with $-w=2L+\sum r_i$ for
the first two and $-w=2L+\sum r_i-1$ for the third. 
By a Weyl invariant of weight $w$ of $\cRt$ (resp. $T\cRt$), we shall mean
the restriction to $\cRt$ (resp. $T\cRt$) of a linear combination of
contractions \eqref{ambcontrR} of weight $w$.   
\begin{lemma}\label{counting}
When restricted to $\cRt$, the first contraction in 
\eqref{ambcontrR} defines a $P$-invariant polynomial on 
$\cRt^{-w-4}$, i.e. it factors through the
projection $\cRt\rightarrow \cRt^{-w-4}$.  The same statement also holds
with $\cRt$ replaced everywhere by $T\cRt$.   
\end{lemma}
\begin{proof}
The argument of the 
second paragraph of the proof of Proposition~\ref{ambcontrprop} shows that 
when expanded, the first contraction in \eqref{ambcontrR} involves only
components    
$\Rt_{IJKL,\cM}$ with $\|IJKL\cM\|\leq 4-2L-w$.  The defining relations of  
$\cRt$ and $T\cRt$ show that when restricted to either of these spaces,
any such contraction with $L=1$   
can be rewritten as a linear combination of contractions with $L\geq 2$.
This gives the bound $\|IJKL\cM\|\leq -w$, which is the desired statement. 
\end{proof}

The main result we will use characterizing invariants of $\cRt$ is the
following.   
\begin{theorem}\label{Rinvariants}
If $n$ is odd, every $P$-invariant polynomial on $\cRt$ is a 
Weyl invariant.  If $n$ is even, every even $P$-invariant polynomial 
on $\cRt$ is a Weyl invariant.    
\end{theorem}
\noindent
{\it Proof of Theorem~\ref{invariants} using Theorem~\ref{Rinvariants}.} 
For $n$ odd, Proposition~\ref{sci} and Theorem~\ref{jetiso} immediately
reduce Theorem~\ref{invariants} to Theorem~\ref{Rinvariants}.  For $n$
even, an even scalar conformal invariant $I(g)$ with $-w\leq n$ 
determines by   
Proposition~\ref{sci} and Lemma~\ref{weightrestrict} an even
$P$-invariant polynomial of weight $w$ on $\cN_c^{n-2}$.
Theorem~\ref{jetiso} shows that this defines an even $P$-invariant
polynomial of weight $w$ on $\cRt^{n-4}$.  Theorem~\ref{Rinvariants}
asserts that as a function on $\cRt$, this
polynomial agrees with a linear combination of complete contractions
\eqref{ambcontrR} of weight $w$, and Lemma~\ref{counting} shows that each 
of these can be regarded as a 
$P$-invariant of weight $w$ of $\cRt^{n-4}$.  Reversing the steps, it
follows that $I(g)$ has the desired form.  
\stopthm

Theorem~\ref{Rinvariants} is proved by reduction to the following, which is
one of the main results of \cite{BEGr}. 

\begin{theorem}\label{TRinvariants}\cite{BEGr} If $n$ is odd, every 
$P$-invariant polynomial on $T\cRt$ is a 
Weyl invariant.  If $n$ is even, every even $P$-invariant polynomial  
on $T\cRt$ is a Weyl invariant.     
\end{theorem}
\noindent
{\it Proof of Theorem~\ref{Rinvariants} using Theorem~\ref{TRinvariants}.}  
We first extend the notion of weight to polynomials which are not
$P$-invariant.  
A polynomial on $\cTt$ is a polynomial in
the variables $(\Rt_{IJKL,\cM})_{|\cM|\geq 0}$, subject to the linear 
relations defining $\cTt$.  The degree of such a polynomial will refer to
its degree in the variables $(\Rt_{IJKL,\cM})$. 
A polynomial $Q$ on $\cTt$ will be said to have
weight $w$ if it satisfies $Q(p_a.(\Rt)) = \sigma_w(p_a)Q((\Rt))$ for
$a>0$.  The linear polynomial $\Rt_{IJKL,\cM}$ has $-w=\|IJKL\cM\|-2$.  
Since $\Rt_{IJKL,\cM}$ vanishes unless $\|IJKL\cM\|\geq 4$, one sees  
that a monomial in the variables $(\Rt_{IJKL,\cM})$ must vanish if its
weight $w$ and degree $d$ satisfy $d>-w/2$.
Observe that each of the relations 
(1)--(6) in Definition~\ref{Rt} defining $\cRt$ states the vanishing of
a polynomial of a specific weight, i.e. $\cRt$ is defined by polynomials 
homogeneous with respect to the action of the $p_a$. 

Let $Q$ be a $P$-invariant polynomial of weight $w$ on $\cRt$, assumed to
be even if $n$ is even.  We will
show by induction that for any $d\geq 0$, $Q$ can be written in the form 
$Q=Q_d + E_d$ on $\cRt$, where $Q_d$ is a linear combination of complete
contractions \eqref{ambcontrR} of weight $w$ and $E_d$ is a  
polynomial on $\cTt$ of weight $w$ which is the sum of monomials of degree  
$\geq d$.  The observation above implies that $E_d = 0$ for $d>-w/2$,
which gives the desired conclusion. 

The induction statement is clear for $d=0$.  Suppose it holds for $d$, so
that $Q=Q_d +E_d$ on $\cRt$.  Write $E_d = E'_d + E''_d$, where $E'_d$ and
$E''_d$ are polynomials on $\cTt$ of weight $w$, $E'_d$ is homogeneous of
degree $d$, and $E''_d$ is a linear combination of monomials of degree
$>d$.  Since $E_d$ is $P$-invariant when restricted to $\cRt$, 
it follows easily using the first part of Theorem~\ref{jetiso} that $E'_d$
is $P$-invariant of 
weight $w$ when restricted to $T\cRt$.  So by Theorem~\ref{TRinvariants},
there is a 
linear combination of complete contractions $Q'_d$ of weight $w$, which 
we may take to be homogeneous of degree $d$ as polynomials on $\cTt$, such
that $E'_d=  Q'_d$ on $T\cRt$.  Let $P_i$, $i=1, 2, \ldots$ be an
enumeration of the polynomials on $\cTt$ whose vanishing defines $\cRt$,
and let $p_i$ denote the linear part of $P_i$, so that $T\cRt$ is defined
by the vanishing of the $p_i$.  We conclude that we can write 
$E'_d -Q'_d = \sum U_i p_i$ as polynomials on $\cTt$, where the $U_i$ are 
homogeneous polynomials of degree $d-1$ and the sum is finite.  When
restricted to $\cRt$, we have $E'_d-Q'_d = \sum U_i (p_i - P_i)$.  
So if we set $Q_{d+1}= Q_d + Q'_d$ and 
$E_{d+1} = E''_d +\sum U_i (p_i - P_i)$, then we have 
$Q = Q_{d+1} + E_{d+1}$ on $\cRt$, where $Q_{d+1}$ is a linear combination
of complete contractions of weight $w$ and $E_{d+1}$ is a polynomial
containing only 
monomials of degree at least $d+1$.  This completes the induction step.
\stopthm

We remark that Theorem~\ref{TRinvariants} also holds
for odd invariants when $n\equiv 2 \mod 4$, by
combining results in \cite{BEGr} and \cite{BaiG}.
The argument given above thus shows that Theorem~\ref{Rinvariants} holds 
also for odd invariants when $n\equiv 2 \mod 4$.   
Of course, the jet isomorphism theorem of Chapter~\ref{jet} does not apply
to higher order jets in even dimensions,  
and $\cRt$ is no longer the correct space to study to understand conformal
invariants.

\end{document}